\documentclass[final,5p,times,twocolumn]{elsarticle}

\usepackage{amssymb}
\usepackage[applemac]{inputenc}
\usepackage{graphicx,amsmath}
\usepackage{ifthen}
\usepackage{endnotes}
\usepackage{multirow}
\usepackage{pdfpages}

\def\beq{\begin{equation}}
\def\eeq{\end{equation}}
\def\bea{\begin{eqnarray}}
\def\eea{\end{eqnarray}}

\makeatletter
\def\ps@pprintTitle{%
  \let\@oddhead\@empty
  \let\@evenhead\@empty
  \let\@oddfoot\@empty
  \let\@evenfoot\@oddfoot
}
\makeatother

\biboptions{square,comma,numbers,sort&compress}

\begin{document}

\begin{frontmatter}


\title{An assessement of global energy resource economic potentials}
\author{Jean-Fran\c{c}ois Mercure \corref{cor1} \fnref{fn1}}
\ead{jm801@cam.ac.uk}
\cortext[cor1]{Corresponding author: Jean-Fran\c{c}ois Mercure}
\fntext[fn1]{Tel: +44 (0) 1223337126, Fax: +44 (0) 1223337130}
\author{Pablo Salas}
\address{Cambridge Centre for Climate Change Mitigation Research (4CMR), Department of Land Economy, University of Cambridge, 19 Silver Street, Cambridge, CB3 1EP, United Kingdom}

\begin{abstract}

This paper presents an assessment of global economic energy potentials for all major natural energy resources. This work is based on both an extensive literature review and calculations using natural resource assessment data. Economic potentials are presented in the form of cost-supply curves, in terms of energy flows for renewable energy sources, or fixed amounts for fossil and nuclear resources, with strong emphasis on uncertainty, using a consistent methodology that allow direct comparisons to be made. In order to interpolate through available resource assessment data and associated uncertainty, a theoretical framework and a computational methodology are given based on statistical properties of different types of resources, justified empirically by the data, and used throughout. This work aims to provide a global database for natural energy resources ready to integrate into models of energy systems, enabling to introduce at the same time uncertainty over natural resource assessments. The supplementary material provides theoretical details and tables of data and parameters that enable this extensive database to be adapted to a variety of energy systems modelling frameworks.

\end{abstract}


\end{frontmatter}


\section{Introduction}

Energy policy decisions for the planning of new energy generation capacity designed to respond to future demand require information regarding the engineering feasibility and the requirements of such systems in terms of capital investments, but also, crucially, the availability of natural resources. Meanwhile, future energy systems are expected, in many contemporary policy frameworks, to evolve towards their gradual decarbonisation, in order to decrease anthropogenic interference with the climate system (see for instance Edenhofer $et$ $al.$ \cite{Edenhofer2010}). From an energy perspective, the decarbonisation of the sector involves a transfer of supply from traditional fossil fuel based technologies towards low GHG emission energy generation capacity such as renewable energy systems or nuclear reactors. Such a transfer requires changes in the technologies used through substitution processes, a subject extensively studied in the past (for instance \cite{Grubler1999, Marchetti1978, Mercure2012}). However, these transformations also require changes in the use of primary energy sources. Realistic energy scenarios of the future can only be designed in a way that does not exceed natural sources and flows of energy which are available in all regions of the world. Therefore, assessments of the potential of natural energy resources are essential to energy planning and policy.

Meanwhile, many models exist that generate scenarios for the evolution of global energy systems in the future (For a recent brief review of IMAGE/TIMER, MERGE, E3MG, POLES and REMIND see Edenhofer $et$ $al.$ \cite{Edenhofer2010}). Such models, in order to generate scenarios that are realisable, must take into account the limits to each type of natural resource. However, while some energy models do not currently take explicit account of limits to resource flows, most of them do not consider their associated uncertainty. Of particular interest in this work is the Energy-Environment-Economy Model at the Global level (E3MG), a macroeconometric model of the global economy, which calculates economic activity in 42 industrial sectors within 20 regions of the world \cite{E3MG, Barker2006, Kohler2006, Scrieciu2010, Dagoumas2010}. In a previous paper, a new sub-component designed for use in E3MG was introduced, modelling the trend of investor behaviour in the power sector facing the choice of various technologies, using their Levelised Cost Of Electricity production (LCOE), which includes the cost of natural resources, and the resulting technological substitutions \cite{Mercure2012, Mercure2011TWP}. In this model, the limitation of natural resources occurs through the use of cost-supply curves, requiring a complete set for each type of natural resource considered in order to properly constrain the model.

Many studies and reviews have been published that summarise what was known of global energy potentials at time of their publication \cite{WEA2001, IPCCAR4Ch4, BGR2010, IPCCSRREN2011,WEC2010, Krewitt2009, WWFEcofys2011}. These reports however do not provide any economic structure to energy potentials. Without underlying individual cost structures associated with energy resource extraction or clarifications the use of the concept of $economic$ potentials, natural resource limit values are ambiguous. As the key for strategic energy planning lies precisely with the cost structure of energy production, such assessments are of limited use for energy systems modelling and policy-making without additional assumptions over their cost distributions. 

Every energy source is limited, either in its total amount that can be consumed, for stock resources, or in the total energy flow it can produce at any one time, for renewable resources. Resources tend to be exploited in order of their cost of extraction. As consumption gradually progresses to higher and higher levels of exploitation of particular resources, additional units consumed tend to incur increases in production cost. Therefore, the economic potential is better defined as a $function$ of cost, rather than as a constant value. As the costs of production increase with the levels of use, developers increasingly seek alternatives, where they exist, through an evaluation of the opportunity cost. Thus, the economic potentials for all energy resources available within a particular market depend onto one another and cannot be determined individually without knowing the alternatives. Additionally, since every particular market for energy is composed in a particular way, economic potentials vary geographically. As resource use evolves and depletion progresses, all economic potential values in a market change with time, an effect that stems from changing individual costs of energy production, and and that has repercussions onto the rate of consumption of every natural resource. 

As discussed in previous work, the economics of energy resources when used, for instance, for electricity production can be modelled in a simplified manner using a complete ensemble of cost-supply curves, which express the cost of resources at various levels of exploitation  \cite{Mercure2012}. In such a framework, the marginal cost of electricity production for every individual natural energy resource using specific power technologies may be compared using a framework such as the LCOE at every level of natural resource use, in order to compare their profitability. In such a model, the consumption of particular resources is limited by increasing marginal costs, and the potential depletion of resources becomes naturally represented since the cost gradually diverges at a levels of exploitation closer and closer to the technical potential. Cost-supply curves for global resources of wind, solar and biomass energy have been calculated previously using the land use model IMAGE~2.2 \cite{HoogwijkThesis, Hoogwijk2004, Hoogwijk2005, deVries2007,Hoogwijk2009}, and are used in this way in the TIMER energy sub-model. Additionally, global cost curves for fossil fuels have been produced by Rogner \cite{Rogner1997}, an influential work which is unfortunately becoming increasingly outdated. However, no comprehensive global assessment of all major energy resources which provides an underlying cost structure currently exists in the literature.

Additionally, assessments of natural energy resources inherently possess high uncertainties, which must be taken into account in order to generate confidence levels in model outputs. For instance, the global bioenergy potential has been estimated to lie between around 310 to 660~EJ/y, by Hoogwijk $et$ $al.$ \cite{Hoogwijk2005}, between 0 and 650~EJ/y by Wolf $et$ $al.$ \cite{Wolf2003} and between 370 and 1550~EJ/y by Smeets $et$ $al.$ \cite{Smeets2007}. These particular uncertainties stem from those on future projections of food demand and levels of technology advancement in the agricultural sector used in these calculations. These types of uncertainties are present in all assessments of renewable energy resources. Similarly, uncertainty arises with knowledge on amounts of stock resources, such as uranium and fossil fuels, where lower and lower levels of confidence are associated with larger and larger quantities. These uncertainties originate directly from the cumulative amount of effort that has been deployed to explore  geological occurrences, and express the fact that it is never possible to know with certainty the detailed composition of the crust of the Earth. Using a review of literature, it is possible, and appropriate, to define ranges of energy potentials rather than specific values, therefore defining areas in the cost-quantity plane for where actual cost-supply curves are likely to lie. 

This work proposes a theoretical framework  and a computational methodology for building natural resource assessments readily useable in models of energy systems, by using a combination of cost-supply curves and a treatment of uncertainty. This methodology is then applied to produce a cost-quantity analysis for every major natural energy resource, those with a potential larger than 10~EJ/y. As part of this work, cost-supply curves were produced for 13 types of resources for every one of the 20 world regions specified in E3MG, and form a new sub-model for natural resource use and depletion. Underlying potentials were however defined for 190 individual countries, and can be aggregated to any other particular set of world regions. For the sake of presentation in this paper, the cost-supply curves were aggregated into global curves, providing a world view. Since data for all 190 countries could not realistically be provided here, tables of parameters are provided in the supplementary material that enable to reconstruct the cost-supply curves for a set of 14 world regions that were assumed the most representative of the requirements of the international modelling community. Additional resource specific information regarding theoretical derivations, additional methodology and justifications are also provided.

This work follows consistently a methodology that can be reused as presented by the modelling community. To this end, definitions of the concepts used are first given followed by a concise description of the approach, detailing the cost-quantity analysis of resources with the associated uncertainty. Following this, a theoretical characterisation of the statistical properties of resource occurrences is given in order to enable the use of functional forms for interpolating resource data. This methodology is then used to produce cost-supply curves for renewables and stock resources. A summary of all major energy resources is given in the last section.

\section{Definitions}

\begin{description}
\item[Renewable and stock energy resources:] Natural sources of energy that may be found in one of two forms: $stocks$, where  energy may be extracted from fixed amounts of geologically occurring materials with specific calorific contents; $renewable$ $flows$, where energy may be extracted from continuously producing  onshore or offshore surface areas with wind, solar irradiation, plant growth, river flows, waves, tides or various forms of heat flows. 
\item[Theoretical Potential:] Total quantity of energy stock or flow estimated to exist or stem from a particular natural process, disregarding its technical recoverability.  
\item[Technical potential:] Total quantity of energy stock or flow estimated to exist or stem from a particular natural process, recoverable using a specific technique, disregarding its level of technical difficulty and the associated costs. 
\item[Economic potential:] Quantity of energy stock or flow estimated to exist or stem from a particular natural process, recoverable at exploitation costs that are competitive compared to all other alternative ways of producing the same energy carrier. Since the cost value considered competitive at any one time changes continuously, the economic potential is expressed as a quantity of stock or flow of energy function of cost. However, this concept is used more conveniently in models when expressed as the inverse of this function, the cost-supply curve.
\item[Cost-supply curve:] Function of the cost of energy flow or stock from a particular resource given that a certain quantity is already in exploitation or has already been consumed (the $marginal$ $cost$). In this work, the cost-supply curve and the economic potential are used interchangeably.
\item[Uncertainty range:] Area in the cost-quantity or cost-flow planes in which actual real cost-supply curves have a 96\% probability of lying. This would correspond to two standard deviations, 2$\sigma$, if the distributions were normal, but they are in general skewed. Real values have a 2\% probability of occurring below the range, and a 2\% probability of lying above.\endnote{The 2$\sigma$ probability range correspond to erf$(\sqrt{2})$ = 95.45\%, yielding 2.28\% as a probability of values occurring above or below the range. The values of 96\% and 2\% are used instead for convenience.}
\item[Productivity:] Amount of energy stock or flow produced per unit of land or sea surface area (wind, solar, biomass, wave energy), bore depth or digging effort (oil, gas, coal, nuclear fuels, geothermal energy) or construction effort (hydroelectricity, tidal energy).
\item[Hierarchical resource distribution:] Statistical type of natural resource productivity distribution, with productivity values that strongly depend on the number of simultaneous positively contributing physical factors. Resource producing items of this statistical type within one kind of resource (windy areas, rivers, mines, wells, etc) possess widely different productivities which enable their ranking in order of resource quality.
\item[Distribution for nearly identical resources:] Statistical type of natural resource productivity distribution in which resource producing items possess nearly identical productivity values, which do not depend on the simultaneous occurrence of several factors. Producing items of this type within one kind of resource (for instance sunny or fertile plots of land) have nearly identical properties, cannot be ranked and can be exchanged for one another.

\end{description}

\section{Methodology}

\begin{figure*}[t]
	\begin{center}
		\includegraphics[width=1.4\columnwidth]{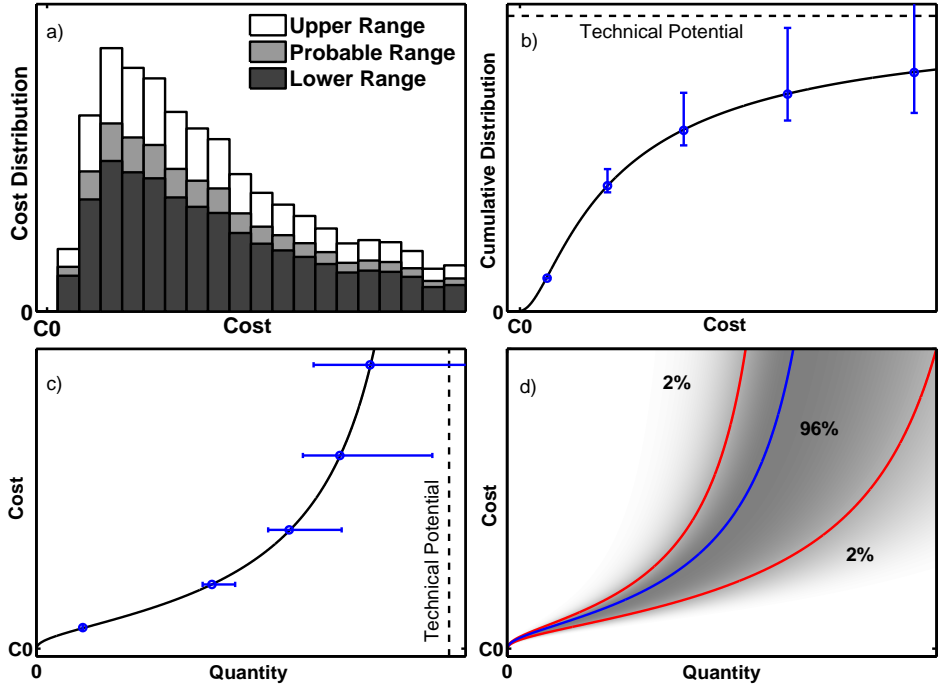}
	\end{center}
	\caption{$a)$ Sketch of a hypothetical distribution of cost ranked amounts of energy or energy flow units available in various cost ranges. The uncertainty over the amount available in each cost range is indicated with a colour shading: the top of the dark histogram represents the minimum amount which has a 98\% probability of being exceeded while the top of the white histogram indicates the maximum amount associated with a 2\% probability of being exceeded. The most likely quantity is intermediate, represented by the top of the grey histogram. $b)$~Cumulative distribution of energy resources, with uncertainty shown as vertical error bars. $c)$~The marginal cost, or cost of extracting an additional unit of energy stock or flow given that a certain quantity has already been exploited, commonly called the cost-supply curve, with uncertainty represented as horizontal error bars. $d)$~Cost-supply curve defined as a probability distribution, where the red curves indicate the limits of a 96\% confidence level region in the cost-quantity plane, while the blue curve corresponds to the most probable cost-supply curve. The assumption is therefore taken that there is a 2\% chance that the cost-supply curve lies below the upper boundary, and a 2\% chance that it lies below the lower boundary. }
	\label{fig:CCurve4All}
\end{figure*}

\subsection{Economic potentials}
Natural occurrences of energy resources are found in different forms, with varying productivity levels or require various levels of effort for their extraction, which enable their transformation into usable energy carriers with different levels of profitability. These variations together lead to particular distributions of costs for their utilisation. Naturally, resources with the best qualities for energy production, and thus lowest extraction costs, are likely to be  considered first by energy firms under financial constraints. Therefore, deriving economic potentials for energy resources involves the task of classifying and ranking different occurrences of specific resources in order of cost.

Information on energy resources is scarce and irregularly distributed, possibly inconsistent, and thus must be organised and classified in order to produce a consistent and complete set of economic potentials. Data may be patchy and incomplete, in which case assumptions are required in order to interpolate through missing parts. Such assumptions are taken in this work in the form of functional forms for the ranking of resources in terms of their cost of extraction. These are derived theoretically from basic statistical properties of resources. They have been carefully verified against several sets of data for specific types of natural resources which do not take any assumptions over the distribution of resources (wind, solar, two types of biomass resources as well as with uranium). They have been assumed to hold true for all other types of resources (for fossil, geothermal, hydroelectric and ocean resources). 

Methods of assessment differ and produce different results or ranges of results. In the absence of justified criteria onto which to base a choice of particular studies over all others, resource assessments must be considered equally, the collection of which can be used to generate uncertainty ranges. This allows to decouple this work from specific assumptions used in specific studies.

The methodology presented here builds upon the approach defined in our earlier work \cite{Mercure2012}. The economic potential of resources is defined using the cost-supply curve, which expresses the quantity of resource available for any cost value considered $economic$, or competitive with all other alternatives. Such curves are derived from cost rankings of resources and resulting distributions. The cost variable, however, stems from varying levels of technical difficulty for extracting resources, or alternatively, the productivity of energy producing resources such as plots of land, mines, oil wells, rivers, etc. Therefore, continuous distribution functions for the amounts of resources available in nature are defined in terms of their productivity. Two empirical forms for these distributions are defined and used throughout. Confidence ranges are derived from uncertainty analysis. The combination of both is used to construct probability densities for the location in the cost-quantity or cost-flow planes where the real cost-supply curves would be situated if it were possible to determine them with certainty.\endnote{Note that the use of uncertainty ranges in the cost-quantity/cost-flow plane relaxes the constraints of using specific functional forms, since it allows variations in the particular forms of the functional dependences within the ranges.} These probability distributions may be used as inputs to uncertainty analysis (such as Monte Carlo simulations) in energy systems modelling.

\subsection{Cost-supply curves}

The ranking of resources in terms of their productivity, required for building cost-supply curves, can be done using a set of histograms of the quantity of energy stocks of flows that can be obtained within various ranges of productivity values. The productivity variables may be converted into costs, which result in a new set of histograms representing the amount of energy that can be produced at costs within various cost ranges. This is shown in Figure~\ref{fig:CCurve4All} $a)$, with a typical distribution of energy resources, which decreases at low cost values due to the decreasing number of resources of exceptional quality, and to high cost values due to the decreasing productivity or recoverability of the resources. The shading is a representation of the confidence associated to their potential availability. The top of the dark grey distribution shows the lowest quantity of assured resources, assumed to be exceeded with a probability of 98\%. The top of the white histogram represents the upper range of speculative resources, those assumed to be exceeded with a probability of 2\%. The quantity which is the most likely to be available lies between these extremes, shown with the top of the  grey histogram. 

The amounts in each cost range thus possess an uncertainty.\endnote{No error bars are present for the cost variables, since an uncertainty on costs corresponds to an uncertainty on how to distribute energy quantities between existing cost ranges, which translates into an uncertainty in the quantity in each range.} In order to determine the quantity that can be obtained at or below certain cost values, the cumulative distribution function is calculated, shown in panel $b)$. This sum converges towards the technical potential. The uncertainty increases approximatively cumulatively with increasing cost values through the root of the cumulative sum of the squares of the individual uncertainty values, shown with blue error bars.

The marginal cost of resources, or the cost of extracting an additional unit of resource given that a certain quantity has already been used, corresponds to the inverse of the cumulative sum, shown in panel $c)$. Thus, the cost of additional units diverges when the number of units used approaches the technical potential, at the point of resource depletion. Using the uncertainty ranges, or error bars, to define two additional curves, assumed to define the upper and lower 96\% confidence limits, a probability density can be defined in cost-quantity space for the location where the real cost-supply curve would lie if it were possible to know it with certainty. This is shown in panel $d)$, where the red curves define the uncertainty area and the blue curve represents the most probable of all possible cost-supply curves, the mode of the distribution. Such probability densities are normally skewed, since the uncertainty over undiscovered resources lies at higher quantity values.\endnote{Thus, the most probable cost-supply curve is neither the mean or the median of the skewed distribution, it is the mode, or maximum.} Note that the uncertainty is assumed to vanish at the contemporary position in the cost-quantity plane, since current costs and levels of exploitation are well known.

\subsection{Distributions \label{sect:Distributions}}

\begin{figure}[h!]
	\begin{center}
		\includegraphics[angle = -90, width=1\columnwidth]{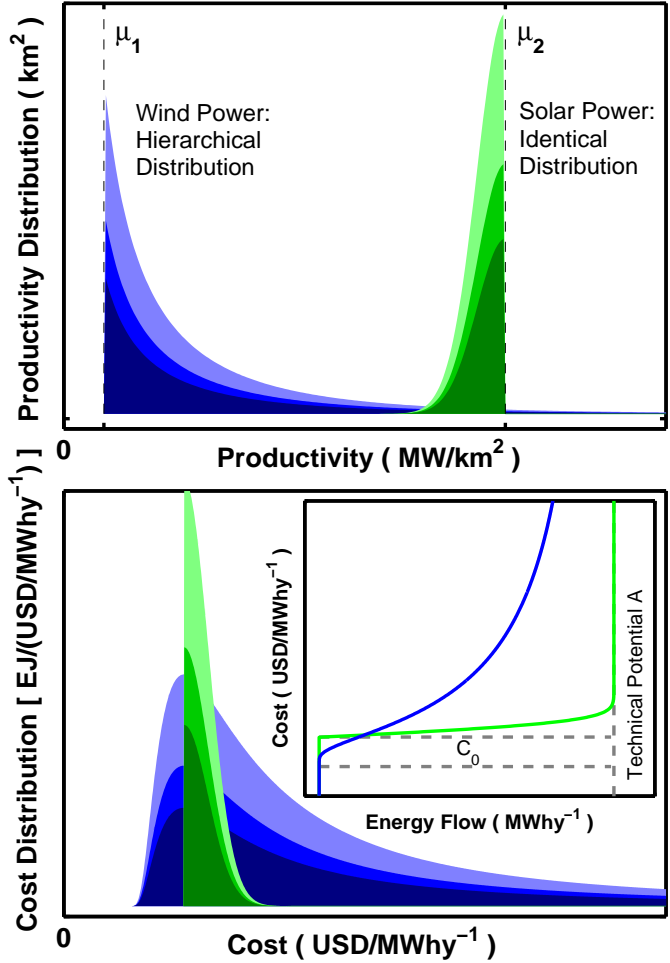}
	\end{center}
	\caption{Depiction of two types of natural resources, based on their statistical properties. $Top$. A typical sharp distribution for nearly identical resource types is shown in green, while the broad blue distribution is for hierarchical resources, from equations~(\ref{eq:DistProd1}) and (\ref{eq:DistProd2}). Both are expressed in terms of productivity. Different colour shadings represent uncertainty, as in panel $a)$ of Figure~\ref{fig:CCurve4All}. $Bottom$. Same distributions expressed as functions of cost, from equations~(\ref{eq:Dist1}) and (\ref{eq:Dist2}) through equation~(\ref{eq:Cost}). Associated cost-supply curves are given in the inset. Note that the technical potential was adjusted to be the same for both curves for visual clarity.}
	\label{fig:ResDistributions}
\end{figure}

Natural resources are scattered around the planet in different forms with different probabilities for the cost of their exploitation. Complex processes underlie the formation of these distributions, however, they may follow certain statistical trends, the nature of which stems from the nature of the resource. One particular property affects significantly these statistical distributions, whether the productivity of resource producing units (rivers, plots of land, wells, mines etc) tend to be very similar, which makes their ranking difficult, as opposed to resource types which have strong ordering. For example, solar energy can be produced using photovoltaic (PV) panels, and these panels may be installed with equal ease almost anywhere, and scarcely populated regions of similar solar irradiation will have large potentials for solar energy situated within very narrow productivity ranges where every plot of land has nearly the same productivity value. These are described by a particular type of statistical distribution for nearly identical resources. The properties of wind energy potentials follow a very different structure. Plots of land within a geographical region are unlikely to possess the same average wind speed, and can therefore be ranked in order of productivity, and thus of cost for the production of wind power. As opposed to solar energy sites, wind power production sites can be described by a hierarchical statistical distribution type. 

Resource distributions for nearly identical resources are sharply defined in a narrow range of productivity, but cut off at a maximum value, which corresponds to the best possible conditions for energy production. Meanwhile, resources with hierarchical distributions occur in large numbers in low productivity ranges, and in ever lower amounts as productivity increases. This property stems from the large number of positively contributing factors which are required simultaneously in order to have a resource producing unit with high productivity. Such a property results in a distribution that decreases exponentially with increasing productivity, but cut off below a certain low productivity value, where it is simply assumed that no energy can be obtained with any reasonable amount of effort. This is shown in Figure~\ref{fig:ResDistributions}, top panel, where a typical distribution for nearly identical resources is shown in green, and a hierarchical distribution is given in blue. Similarly to panel $a)$ of Figure~\ref{fig:CCurve4All}, the colour shading indicates uncertainty.

These resource distributions are well described by the following density functions:
\beq
f(\nu)d\nu =  \left\{ \begin{array}{cc} {A \over \sigma} e^{-{\nu \over \sigma}}d\nu  & \nu > \mu_1 \\ 
0 & \nu \leq \mu_1 \end{array} \right .,
\label{eq:DistProd1}
\eeq
\beq
g(\nu)d\nu =  \left\{ \begin{array}{cc} {A \over \sqrt{2\pi} \sigma} \nu e^{-{(\nu-\mu_2)\over 2 \sigma^2}^2}d\nu  & \nu < \mu_2 \\ 
0 & \nu \geq \mu_2 \end{array} \right .,
\label{eq:DistProd2}
\eeq
where $\nu$ is the productivity (see definition above), $A$ is the technical potential, $\sigma$ is the width of the distribution and $\mu_1$ is the minimum usable productivity, in the first case, and $\mu_2$ is the maximum productivity available in the second case. 

Costs are related to the productivity by an inverse relationship,
\beq
C = {C_{var} \over \nu} + C_0,
\label{eq:Cost}
\eeq
with which the distributions can be transformed into cost-quantity space. The scaling factor $C_{var}$ corresponds to a cost per unit of land or sea area (wind, solar, biomass, wave energy), digging effort (oil, gas, coal, nuclear fuels, geothermal energy) or construction effort (hydroelectricity, tidal energy), and the ratio $C_{var}/\nu$ has units of cost per unit of energy produced.  $C_0$ corresponds to the sum of costs which do not depend on the productivity, per unit of energy produced. The conversion of these productivity distributions into cost distributions is described in detail in section~S.2 of the supplementary material, and yields the following:
\beq
f(C)dC =  \left\{ \begin{array}{cc} {A B \over (C-C_0)^2} e^{-{B \over C-C_0}}dC  & C > C_0  \\ 
0 & C \leq C_0 \end{array} \right .,
\label{eq:Dist2}
\eeq
\beq
g(C)dC = \left\{ \begin{array}{cc} {A \over \sqrt{2\pi} B} e^{-{(C-C_0)^2\over 2 B^2}}dC & C > C_0\\
0 & C \leq C_0 \end{array} \right .,
\label{eq:Dist1}
\eeq
where $A$ is the technical potential, $B$ a scaling factor and $C_0$ a cost offset, the set of three parameters required to define every distribution given in this work. These functions are illustrated in Figure~\ref{fig:ResDistributions}, where the inset shows the associated cost-supply curves. It is observed that for a similar technical potential, the curve for nearly identical resources possesses less curvature up to very near the technical potential than those for hierarchical resources, a property that stems from their lack of ordering, and results in similar cost values for most of the resources.

These functional forms have been found to reproduce very closely the cost-supply curves calculated by Hoogwijk $et$ $al.$ using the land use model IMAGE, whose work does not assume any functional dependence on cost for its distributions  \cite{HoogwijkThesis, Hoogwijk2004, Hoogwijk2009}. Distributions were calculated by producing cost ranked histograms of calculated potential renewable energy flows (wind, solar and biomass) at every point of a 0.5$^{\circ}\times$0.5$^{\circ}$ grid of the earth's onshore land. Thus, their form originates purely from statistical properties of the aggregation and ranking of the resources modelled. Using least-squares non-linear fits, the cost-supply curves in their work were found to agree very well with one or the other of the functions given above, depending on the nature of the resources: solar energy and agricultural land are well represented by the distribution for nearly identical resources, while wind power and rest land are well represented by the hierarchical distribution. Additionally, the distribution for hierarchical resources was found to agree well with observed cost distributions of uranium as reported by the International Atomic Energy Agency (IAEA) \cite{IAEA2009}. Examples of non-linear fits of these functions to IMAGE data are given in section~S.2.5 of the supplementary material. 

No such global cost ranked data exist for the remaining types of natural resources that could enable fits of distributions distribution. Choices of distributions were therefore taken as assumptions, made based on the physical nature of the resources. Potential basins that could be created for hydroelectricity possess very individual characteristics, which makes them hierarchical. Geothermal resources, however, were treated as a hybrid mixture of the two, since good geothermal sources in active volcanic areas such as Iceland can be ranked, but large amounts of very similar sites can be found in non active areas. Stock resources however were treated slightly differently. Different oil and gas occurrence subtypes (e.g. conventional gas, shale gas, clathrates, etc) originate from different geological processes which have no strong relation to one another, and should therefore be treated independently. These resource subtypes are characterised with different costs of extraction. They were assigned one hierarchical distribution per subtype, and subsequently aggregated. This resulted in composite cost-supply curves with complex structures. Coal, uranium and thorium resources were considered to occur as a single subtype, supported by the fact that the data for uranium was found to follow well the hierarchical distribution.

\subsection{Uncertainty\label{sect:uncert}}

The methodology used in this work for treating uncertainty is fundamental to this analysis of economic potentials as it allows the incorporation all available information, even when sources are inconsistent or conflicting. Inconsistencies can be found between assessments for most individual natural resources, and are the result of the use of different approaches and assumptions, which can be determinant for the technical potential values derived. This is most obvious in resources such as wind power, solar energy and bioenergy, where the total amount of appropriate land depends highly onto competing activities, making the assumptions in the evaluation of the land suitability factors the main drivers of uncertainty. Other such assumptions are world population and the associated food demand, levels of technological development and changes in agricultural productivity. Resource assessments are uncertain by nature, since it is not possible to know with certainty the complete geological content of the crust of the earth, or to predict the weather and associated wind, sunshine and rainfall with perfect foresight. Thus the comparison of ever larger numbers of natural resource assessments is the key to define ranges of confidence, and these are as important as their associated most probable potential. 

This work uses a consistent methodology to define probability distributions for cost-supply curves. Three cost-supply curves are derived from resource assessment data, where two are used to define the 96\% confidence region in the cost-quantity plane, and one taken as the most probable of all possible curves. In all plots of this work, the most probable cost-supply curves are given in blue and the 96\% confidence limits are displayed in red. Uncertainty ranges are almost always asymmetric, since upper ranges are intrinsically characterised by smaller amounts of accumulated knowledge. 

Uncertainty is treated differently for renewable resources compared to stock resources. For renewables, cost-supply curves were obtained from the literature or calculated and taken as the most probable curves, while the 96\% confidence limits were obtained by scaling the technical potential to the limits of its uncertainty range, defined by an ensemble of different studies.\endnote{This is done in order to avoid inconsistencies where curves calculated independently, for instance by fitting data, could cross in some cases.} In the case of stock resources, all resource assessments provide classifications associated with cost ranges and various levels of confidence. In these cases, three cost-supply curves were calculated by assigning probabilities to uncertainty classifications, as in panel $a)$ of Figure \ref{fig:CCurve4All} (i.e. reserves were assumed to exist with a 98\% probability, while reserves plus all speculative resources were assumed to be available with a 2\% probability). Individual methodologies for all types of resources are described in the supplementary material.

\section{Renewable energy resources}

\begin{figure*}[p]
	\begin{minipage}[t]{1\columnwidth}
		\begin{center}
			\includegraphics[width=1\columnwidth]{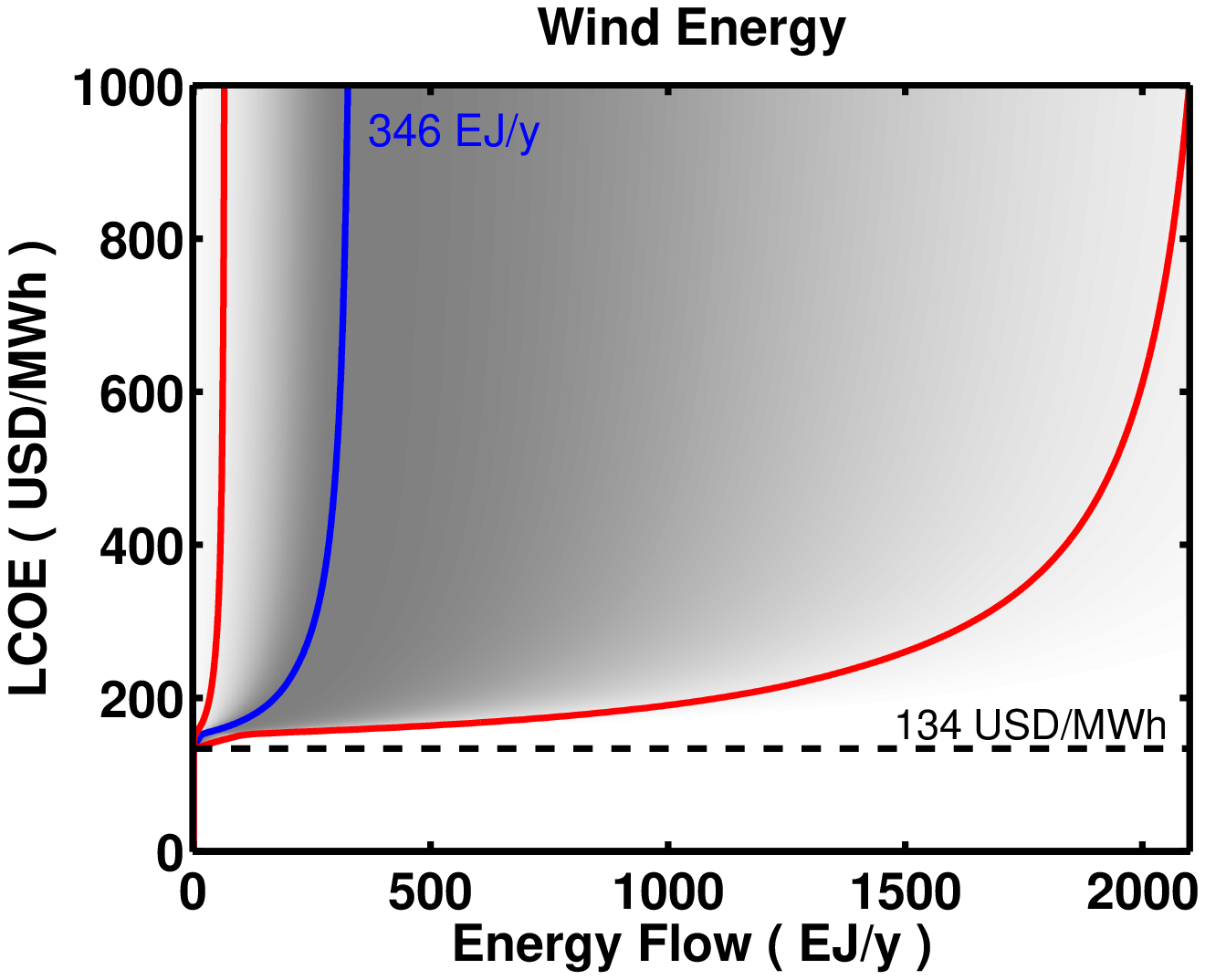}
			\includegraphics[width=1\columnwidth]{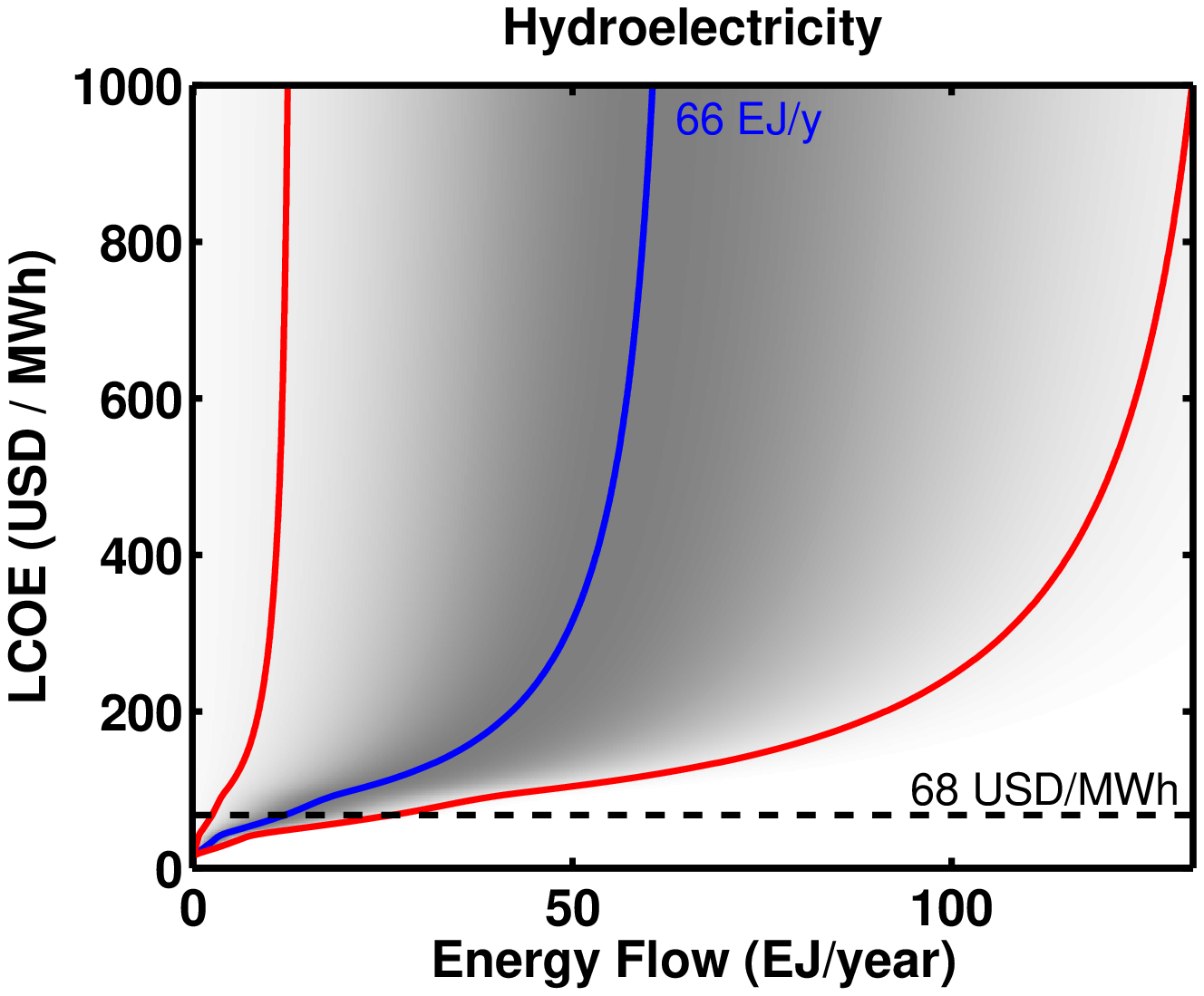}
			\includegraphics[width=1\columnwidth]{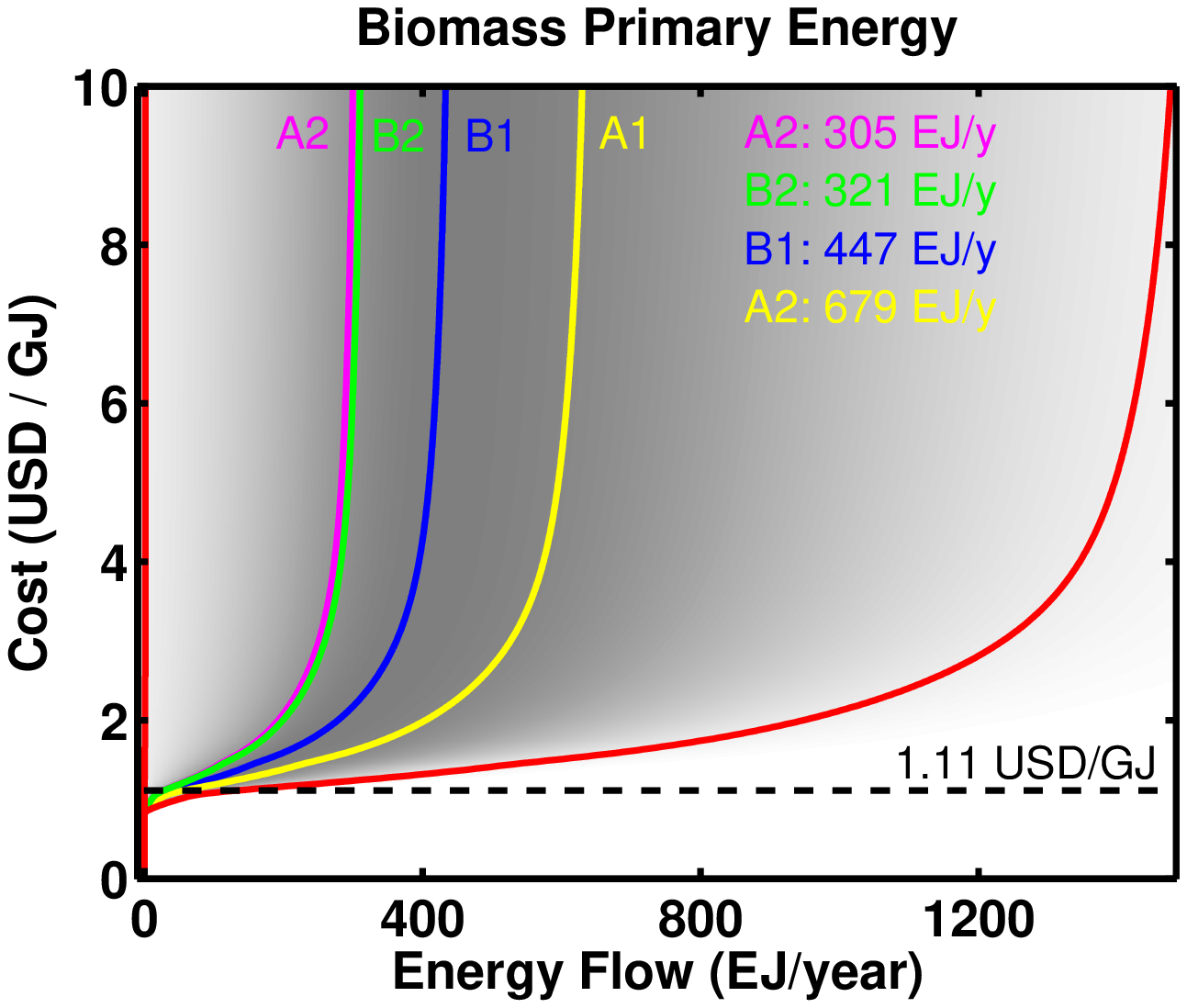}	
		\end{center}
	\end{minipage}
	\hfill
	\begin{minipage}[t]{1\columnwidth}
		\begin{center}
			\includegraphics[width=1\columnwidth]{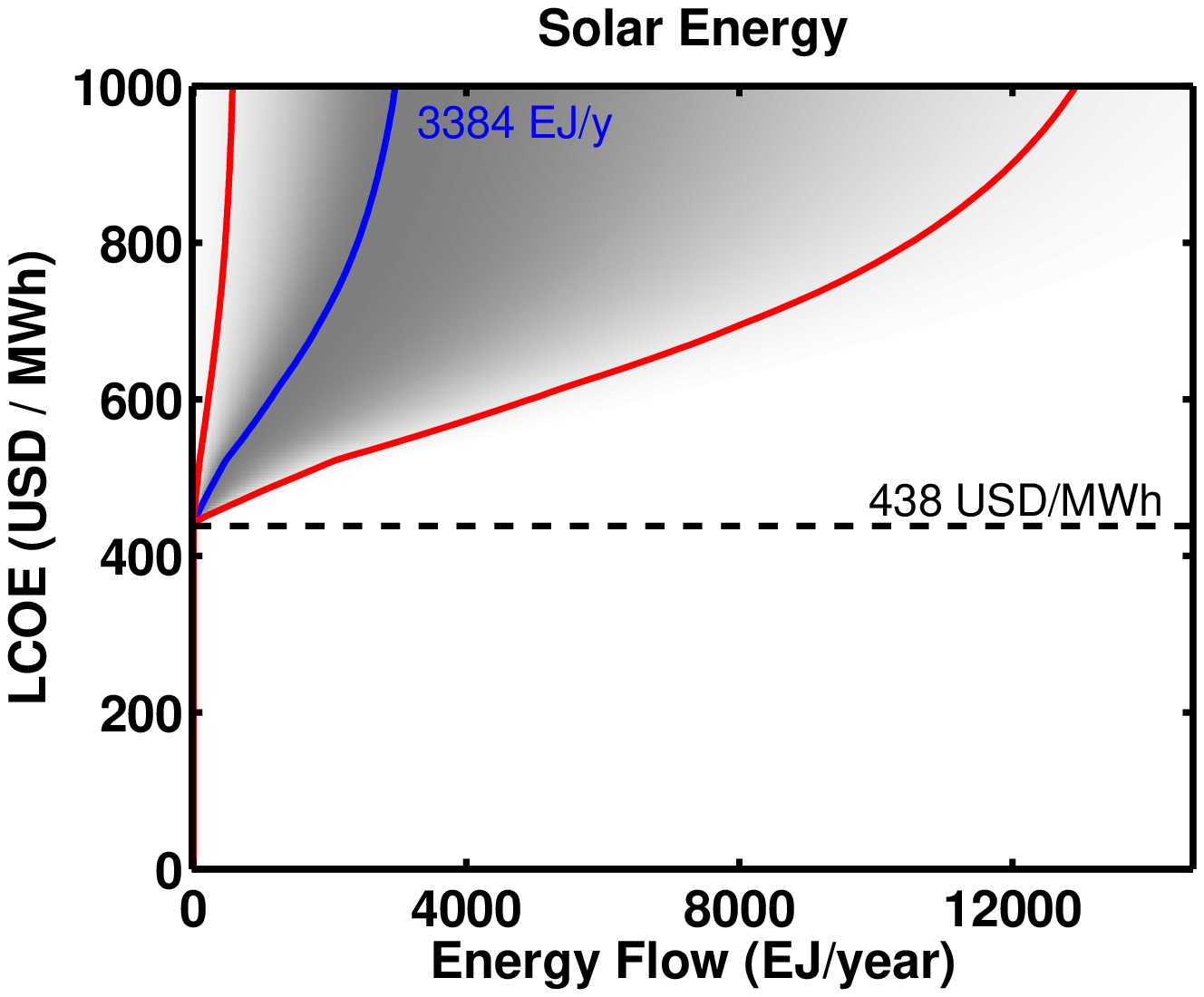}
			\includegraphics[width=1\columnwidth]{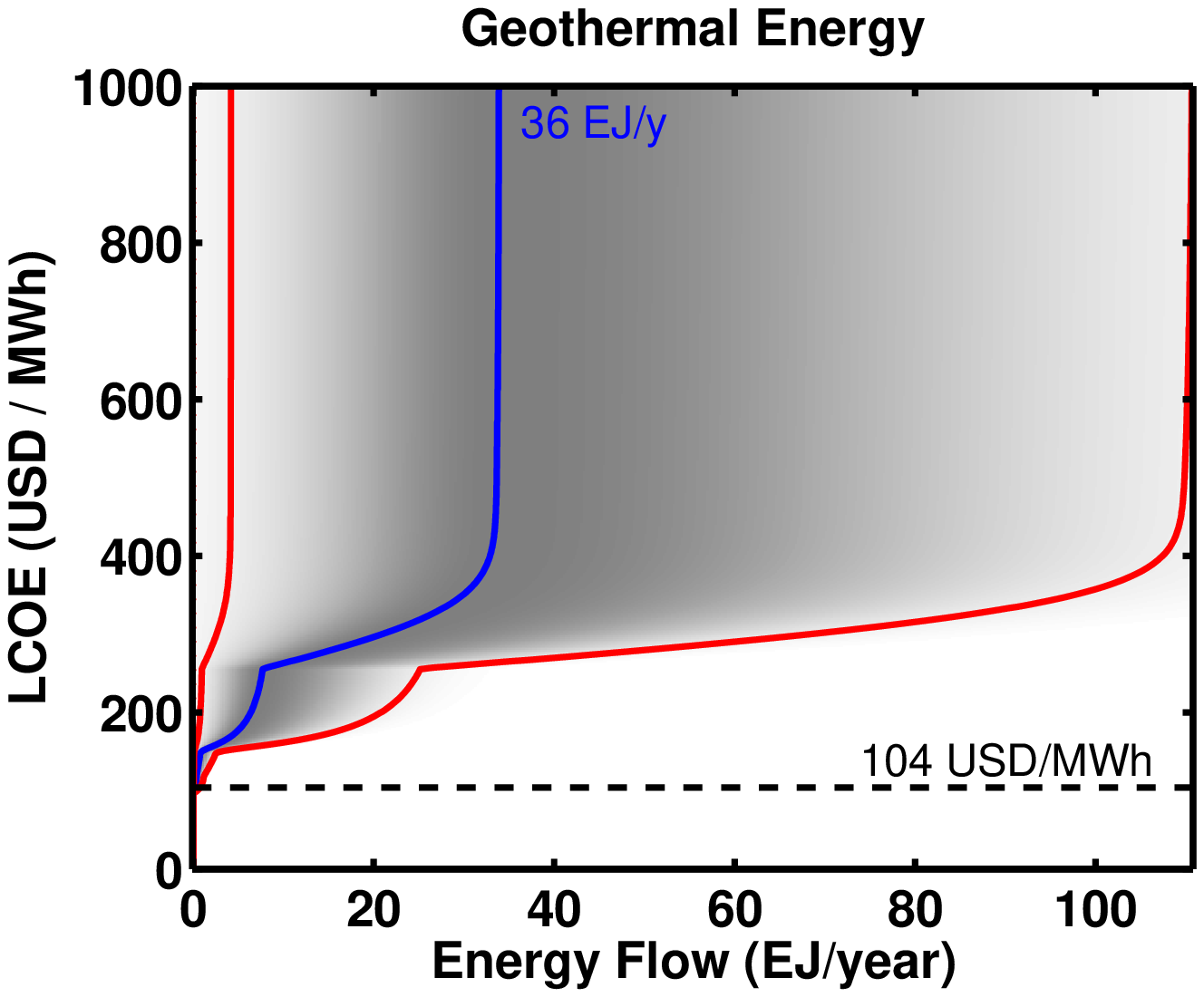}
			\includegraphics[width=1\columnwidth]{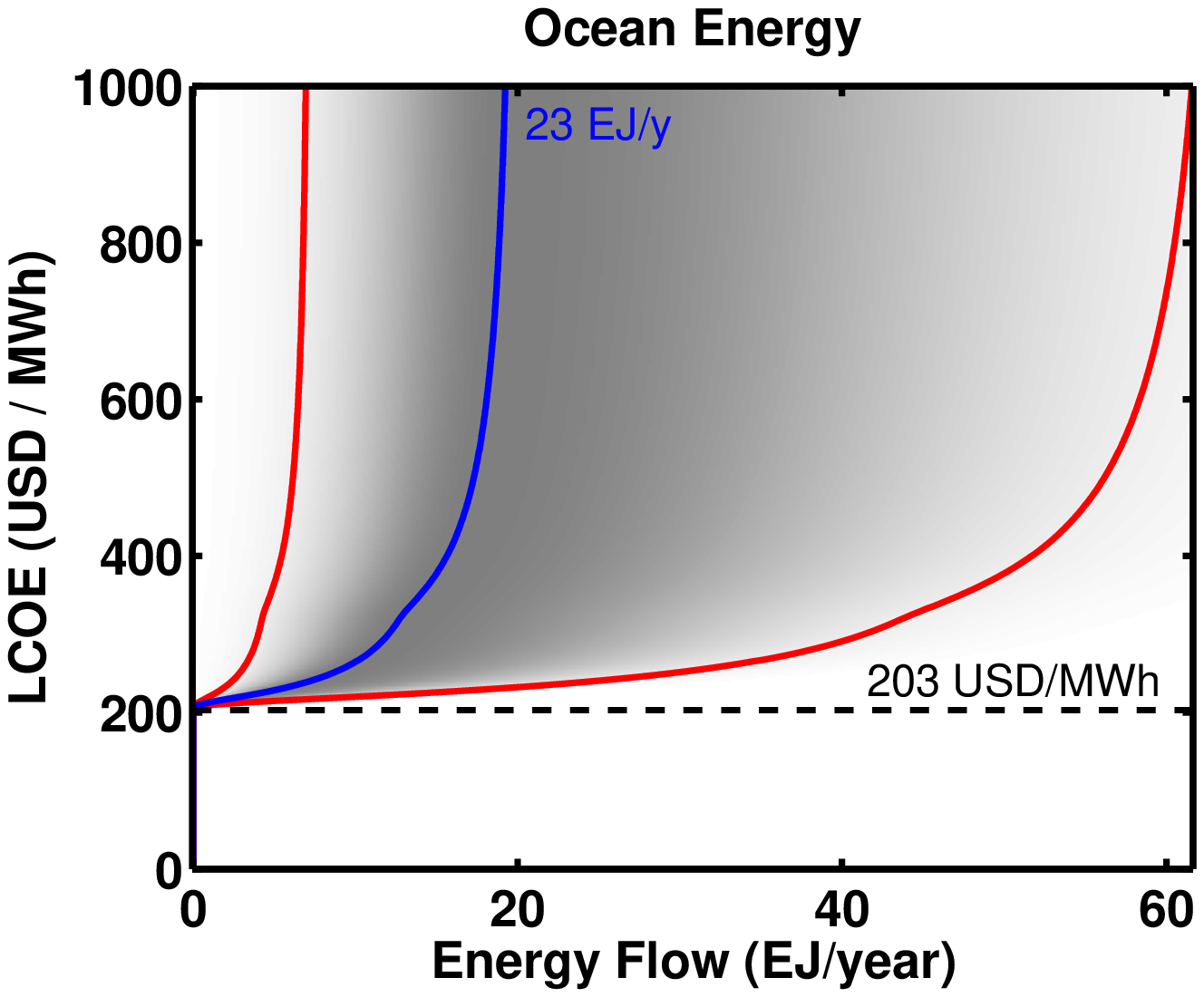}
		\end{center}
	\end{minipage}	
	\caption{Cost-supply curves for renewable resources: wind power, solar energy, hydroelectricity, geothermal power, biomass and ocean energy.}
	\label{fig:GlobalRenCSC}
\end{figure*}

\subsection{Wind energy \label{sect:wind}}

Wind speeds depend strongly on altitude as well as on landscape topologies, the climate and the type of land cover, or roughness. In general, wind speeds increase logarithmically with elevation at low altitude (see for instance S{\o}rensen \cite{Sorensen2011}), and, for a specific elevation and geographic location, occur statistically following a well defined Weibull probability distribution which decreases both towards zero and large wind speeds (for instance Grubb and Meyer \cite{Grubb1993}). Average wind speeds on sites useable for energy production, for instance in the United Kingdom, range between 5.1 and 9 m/s at 10m elevation, the lower boundary determined by technology and the upper limit by the decreasing supply of locations with large wind speeds (for instance in \cite{Sorensen2011}). The power density offered by land areas must be calculated using technical characteristics of particular turbines. For a particular site, one integrates over all wind speeds the product of the site wind speed probability distribution and the turbine power curve. However, a correlation exists between the yearly averaged wind speed and the number of full load hours \cite{Abed1997}. The minimum possible distance between a turbine and its neighbours in a wind farm are determined by losses produced by the wake of neighbouring turbines which results in lower wind speeds and increased turbulence, and scales with the turbine rotor diameter, limiting the density of energy that can be extracted per unit of land area (for instance in the work of Mackay \cite{Mackay2008}).\endnote{While larger turbines intercept a larger wind front, they are also spaced further apart in two spatial directions. Thus, while the power production of large wind farms scales with the square of the length of the blades, it scales inversely with the square of the distance between turbines. These two effects almost cancel each other out, except for the fact that taller turbines intercept higher wind speeds at higher altitudes.} Hoogwijk $et$ $al.$ assumed a maximum density of energy production of 4~MW/km$^2$ \cite{Hoogwijk2004}.

Various research groups have calculated the global distribution of wind power, resulting in a range of values between 70 and 2509~EJ/y for onshore wind power  \cite{Lu2009, Hoogwijk2004, Fellows2000, Archer2005, WEC1994, Grubb1993}, and about 57.4~EJ/y for offshore wind power \cite{Krewitt2009}. A global onshore value of 346~EJ/y has been derived Hoogwijk $et$ $al.$ \cite{Hoogwijk2004}, in whose work, used for the present analysis, estimations of average wind speeds were applied to points on a global onshore grid, as well as the land suitability for the installation of wind farms using the land use model IMAGE~2.2.\endnote{For details on IMAGE see \cite{IMAGE}.} Energy potentials obtained from yearly averaged wind speeds determined on every point of the grid were subsequently aggregated into various cost-supply curves for specific regions of the world, according to the land suitability factor of each point. Cost values were determined using a present value calculation including fixed and variable costs, capacity factors and energy densities associated to particular land areas. 

Wind farm sites follow very strongly the distribution for hierarchical resources (see for instance the exceptional fit of Figure~S.4.1 in the supplementary material), since good sites with average wind speeds exceptionally suitable for energy production are geographically scattered, and the majority of areas in any region of the world possess mediocre average wind speeds, allowing a strict ordering of resource units (see for instance the European Wind Atlas \cite{EuropeanWindAtlas1989}). The profitability of a wind farm venture depends strongly on the quality of the wind resource, determined through capacity factor and average turbulence values. For a fixed turbine investment cost, low capacity factors increase dramatically the cost per unit of electricity produced. 

Figure~\ref{fig:GlobalRenCSC} presents the global economic potential for wind power. It gives an aggregate cost-supply curve, using for the most probable curve the data calculated by Hoogwijk $et$ $al.$, which involves the most detailed methodology for determining the suitability factors, resulting in a technical potential of 346~EJ/y \cite{Hoogwijk2004}.  While Lu $et$ $al.$ estimated an optimistic technical potential of 2509~EJ/y by calculating wind potentials over the global onshore area and excluding low wind areas by restricting capacity factors to values above 20\% \cite{Lu2009}, Archer $et$ $al.$ calculated a potential of 2257~EJ/y in an assessment where the land included was restricted to class~3 wind energy sites but did not include alternate uses of the land, a value taken as the upper boundary of the uncertainty range \cite{Archer2005}. Meanwhile, an estimate of 70~EJ/y was obtained using an evaluation of the number of sites with average wind speeds above 5.1~m/s, but with an arbitrarily chosen value of 4\% of that land available for wind turbine installation, in order to account for alternative land uses, taken as the lower boundary of the uncertainty range \cite{WEC1994}. All other existing studies result in values within this range \cite{Grubb1993, Fellows2000}.\endnote{Fellows $et$ $al.$ concludes with an estimate for 2020 of 148~EJ/y, while Grubb and Meyer calculated a global potential of 191~EJ/y \cite{Fellows2000,Grubb1993}.} 

It is to be noted that using up a large fraction of that potential results in large areas becoming covered by wind farms. Since typical individual wind turbines currently have capacities of 3~MW but capacity factors of about 25-35\%, compared to 80\% for coal power stations, replacing one coal power station of 1~GW for wind energy requires 700 to 1100 turbines, covering an area of 500 to 1500~km$^2$, compared to about 1~km$^2$ for the original power station\endnote{The variation originates from both assumed ranges in capacity factors of 25-35\% and turbine densities of 2.2 to 4~MW/km$^2$ \cite{Mackay2008,Hoogwijk2004}.}. Even though agricultural land used by wind farms may still be cultivated, and therefore a competition for land with agriculture does  not directly occur, strong emissions reductions pathways based on substituting fossil fuels for renewables likely implies that, given the large number of wind turbines required, these would invade permanently traditional rural landscapes. 

In the case of offshore wind power, although distributions of wind speeds at offshore locations tend to have a higher median, the air flow possesses less turbulence and wind speeds vary less in time, all of which contribute to produce a higher power density in terms of geographical area, the total area where such turbines can be installed is small compared onshore areas, unless floating turbines become widely available \cite{Weinzettel2009}. The potential of offshore wind energy was evaluated to 22~EJ/y by Hoogwijk $et$ $al.$ and to 57.4 by Krewitt $et$ $al.$, both based on the work of Fellows $et$ $al.$  with power density values of 10 and 12~MW/km$^2$ \cite{HoogwijkGraus2008,Krewitt2009,Fellows2000}.  This potential could be at best only approximately six times lower than the most probable potential of onshore energy given in Figure~\ref{fig:GlobalRenCSC}, and costs per unit energy are significantly higher \cite{EWEA2009}. Due to the lack of reliable and consistent cost data, a cost-supply curve for offshore wind is not presented in this paper.

\subsection{Solar energy \label{sect:solar}}

Solar radiation over the Earth surface is of about 1.2$\times$10$^5$~TW, or 3.6$\times$10$^6$~EJ/y \cite{Crabtree2008}. The fraction of that energy that can be harvested with existing systems has been estimated by several studies, and ranges between about 1340 to 14800~EJ/y \cite{Hofman2002, HoogwijkThesis, deVries2007}. Even though the total generation of energy from solar technologies has been increasing steadily during the last two decades, they still represent very low percentages in their respective categories; solar heating systems account for 0.3\% of the total energy used for heating in 2008, and solar electricity generation represented only 0.06\% of the total electricity generation during the same year \cite{IEAWEO2010}. 

Solar energy can be harvested using either of two existing technologies, photovoltaic (PV) devices (see \cite{Avrutin2011}) or concentrated solar power (CSP) (see the International Energy Agency (IEA) Energy Technology Systems Analysis Programme (ETSAP) \cite{IEAETSAP2010CSP}). Single crystal silicon photovoltaic diodes currently have light conversion efficiencies of up to 25\%, while III-V semiconductor cells such as GaAs systems have efficiencies of up to 28\%, and solar cells using concentrated sunlight can convert light as efficiently as 43.5\% \cite{Green2011}. The resulting electricity generation energy density ranges between 5 and 100~MW/km$^2$, depending on the type of devices used and the geographical location. On the other hand, CSP technology uses a traditional steam turbine where water was heated using sunlight concentrated with various arrangements of mirrors, and have efficiencies of around 13-24\% and energy densities near 25~MW/km$^2$ \cite{IEAETSAP2010CSP,IPCCSRREN2011}. These systems are however restricted to high irradiance areas, and therefore have a lower global technical potential. Solar energy is well represented by a distribution for nearly identical resources, since within areas of similar irradiance and average cloud coverage, its cost does not depend strongly on the nature of the land, having identical productivity values that depend solely on the chosen technology, and the opportunity cost of the land is the limiting factor to its technical potential. 

Figure~\ref{fig:GlobalRenCSC} shows the global economic potential for solar energy, using PV as technology. It is an aggregation of curves determined in various world regions, based on the work of Hoogwijk \cite{HoogwijkThesis}, but rescaled to match the values of de Vries $et$ $al.$ \cite{deVries2007}. The global technical potential is of 3384~EJ/y. Regional cost-supply curves were drawn from an analysis performed using the land use model IMAGE~2.2 to determine land suitability factors and the opportunity cost of the land at every point of a global grid, while regional estimates of solar irradiation were used to determine the energy potential, which have been taken here for the most probable cost-supply curve. The lower boundary curve of the uncertainty range, with a technical potential of 1340~EJ/y, is the one evaluated by Hoogwijk $et$ $al.$, while the most probable curve was taken as the 2000 value of de Vries $et$ $al.$. The upper boundary curve of the uncertainty range, with a technical potential of 14778~EJ/y, is the 2050 projection of de Vries $et$ $al.$, and stems mostly from increases in land availability following future reductions in the amounts of land required for agriculture through improvements in productivity \cite{deVries2007}. 

The use of a large fraction of the technical potential for solar energy using PV systems signifies covering up large amounts of land with solar panels. Mackay calculates an average productivity value of 10~MW/km$^2$ for the United Kingdom, while Hoogwijk used values between 6 and 25~MW/km$^2$, depending on the geographical location, with capacity factors of up to 50\% \cite{Mackay2008,HoogwijkThesis}. When compared to coal power plants of capacity of 1~GW and capacity factor of 80\%, this implies that the replacement of one such power plant by solar panels requires an area between 50 and 500~km$^2$.
 
\subsection{Hydroelectricity}

Hydropower stems from water pressure gradients that are produced by the run-off of rainfall through landscape topographies, using dams to restrict water flow and accumulate water at elevation level higher than that given by the landscape, producing a potential for electricity production using turbines (see for instance \cite{IPCCSRRENHydro2011}). Hydroelectricity is the most deployed renewable electricity technology, with a global installed capacity of close to 1~TW, which produced around 2\% of the total primary energy supply in 2008 \cite{IEAWEO2010}. As a fraction of its total technical potential, it is also the most developed of all renewable resources, to the extent that around 23\% of the global hydroelectric technical potential is currently in use. However, its exploitation around the world is not even: 25\% of the European technical potential has already been developed, while Africa uses only 7\% of its hydroelectric resource \cite{IJHD2011}. 

Costs of hydroelectric systems are highly site-specific and were found to have varied between around 400 to 4500~USD2002/kW in an extensive global analysis done by Lako $et$ $al.$ \cite{Lako2003}. These values are influenced by many different factors, which include material and labour costs, but also critically the opportunity cost of the land. The latter refers to the consequences of flooding large areas of land, and the resulting displacement of communities and agricultural activities, and thus varies strongly from region to region. For this reason, the deployment of hydropower is often decided on political rather than financial grounds. Hydroelectric resources were assumed in this work to follow a hierarchical distribution, since available natural basins that can be flooded possess vastly different geographic characteristics that make them unique and produces strong ordering. 

Figure~\ref{fig:GlobalRenCSC} shows the global cost-supply curve for  hydroelectricity, where the dotted line indicates the cost at the current deployment of 12~EJ/y, 23\% of the modest value of the most probable technical potential of 66~EJ/y, calculated in this work from data gathered in the World Atlas and Industry Guide 2011 \cite{IJHD2011}. The high deployment to potential ratio is an indication that the remaining number of suitable sites for building dams is relatively limited. The intersection of this curve with the current total hydroelectricity generation value yields a cost of production of about 68~USD/MWh. However, the development of hydropower projects hardly follows an order based onto cost, but follows instead an order dictated by political considerations, which are out of the scope of this work. Therefore, this value is only indicative, and projects with LCOE values between 23 and 460 USD/MWh have been recently developed \cite{IEAProjCosts}. The use of this cost-supply curve in modelling involves the inevitable assumption of development following a cost order. In long term scenarios, this is reasonable, since the development of the limited number of remaining available sites involves either increasing opportunity costs in inhabited areas due to increasing local populations, or increasingly large transmission costs associated with increasing distances to uninhabited areas. The cost-supply curve was derived using the theoretical, technical and economic local potential values from \cite{IJHD2011}. Since the definition of the economic potential in IJHD is not given, the (asymmetric) range of the distribution of recent cost values in the work of Lako $et$ $al.$ was interpreted (in 2008 dollars) as what is currently assumed economic. The remaining technical potential (above the economic potential) was assumed to involve higher costs. The upper boundary curve of the uncertainty range was derived by considering the aggregated global theoretical potential of 144~EJ/y from the data in the World Atlas and Industry Guide, while the lower boundary curve of the uncertainty range was derived by assuming that no additional construction of hydroelectric dams occurs in the future, limiting future hydroelectric generation to 12~EJ/y.

\subsection{Geothermal energy}

Geothermal resources, stored beneath the Earth’s surface in the form of heat, are heat sources constantly replenished by the radioactive decay of isotopes of uranium, potassium and thorium (see for instance \cite{Macdonald1959, Wasserburg1964}). Although geothermal heat has been used since prehistory, and its utilization for electricity generation commenced at the beginning of the last century, its current deployment is small in comparison with other sources of energy.  It currently accounts only for 0.3\% and 4\% of the total electricity generation and heating production respectively \cite{Cataldi1993, Lund2005,IEAWEO2010}.

Geothermal resources are classified in four categories: hydrothermal (liquid and vapour dominated), hot dry rock (where fluids are not produced spontaneously), magma (molten rock in regions of recent volcanic activity) and geopressured (hot high-pressure brines containing dissolved methane) \cite{Mock1997}. The most commonly used type is hydrothermal, although high expectations exist regarding the development of Enhanced Geothermal Systems (EGS), oriented towards the hot dry rock type through hydraulic stimulation \cite{Tester2006}. According to estimations made by Aldrich $et$ $al.$ based on a report of the Electric Power Research Institute (EPRI), the estimated geothermal heat accumulated under the continental masses to a depth of 5~km depth is of approximately 1.46$\times10^8$~EJ, most of it assumed to be stored in rocks and water, with a proportion of 6:1 in favour of the former \cite{Aldrich1981,EPRI1978}.  Even though this is a vast amount of heat, only a small part is recoverable for productive purposes.\endnote{Note that the average replenishment of the geothermal heat underground is several orders of magnitude inferior to the stock of heat currently available: around 65~mW/m$^2$ at the continental level, producing an average thermal energy recharge rate of about 315~EJ/yr \cite{Pollack1993}. This value can be considered as the theoretical potential of geothermal energy if viewed in terms of sustainable extraction of geothermal resources over an extended period. However, the amount of time over which geothermal resources could be used at higher rates than this is likely to be more than one thousand years.}

While geothermal resources are available all over the world, their accessibility differs from site to site  according to various technical characteristics including the geological structure of the ground and the depth and type of heat reservoirs. In the vicinity of tectonic plate boundaries, narrow zones characterised by significant volcanic activity (so-called volcanic belts), geothermal gradients are particularly high, between 40 and 80~$^{\circ}$C per km of depth, enabling the extraction of high temperature geothermal resources.  On the other hand, areas with low volcanic activity are characterized by low and uniform geothermal gradients: around 25~$^{\circ}$C per km of depth \cite{EPRI1978, Aldrich1981}. The extraction of geothermal resources in active areas are highly site-specific, and thus were assumed to follow a distribution for hierarchical resources \cite{Pasqualetti1983}. Meanwhile,  geothermal gradients in the rest of land masses have very similar properties and costs, and were therefore assumed to follow a distribution for nearly identical resources. Cost-supply curves were produced for both types of land and both hydrothermal and EGS technologies, generating four curves which were subsequently aggregated in each world region.

Stefansson found a high correlation between the number of active volcanoes in a particular region and the estimate of the size of hydrothermal resources for electricity generation in the same region \cite{Stefansson2005}. Using the total number of volcanoes active in the world, discarding those located on the sea floor or in arctic regions, he estimated a global installable hydrothermal electricity producing capacity of approximately 200~GW (producing 6.0~EJ/y of electricity with 95\% capacity factor). Using this information, along with the statistical analysis between wet and dry systems developed by Goldstein $et$ $al.$ \cite{Goldstein2009}, Bertani estimated the total geothermal installable electricity producing capacity of 1200~GW (36~EJ), including hydrothermal and EGS technologies \cite{Bertani2010, Bertani2012}.

The global aggregation of curves yields the cost-supply curve presented in Figure~\ref{fig:GlobalRenCSC}. The associated global technical potential of 36~EJ/yr, involves a 95\% capacity factor. Cost values were obtained from the International Energy Agency (IEA) \cite{IEAGeo2010}. The lower and upper boundaries of the uncertainty range, of 4 and 114~EJ/y, are explained in section~S.3.4 of the the supplementary material.

\subsection{Bioenergy \label{sect:biomass}}

Bioenergy, energy derived from plants, is currently the most widely exploited renewable energy resource, with 51~EJ/y, 10\% of global annual primary energy use \cite{IEAWEO2010}. The combustion of biomass derived fuels is nearly carbon neutral if CO$_2$ uptake during plant growth is taken into account, minus losses occurring in transformation processes. Thus, biomass based technologies provide an important emissions mitigation potential. While biomass combustion using integrated gasification combined cycle technology (BIGCC) is expected to become the most efficient biomass based electricity production method \cite{Rhodes2005}, the combination of biomass and carbon capture and storage technology has been shown to produce negative CO$_2$ emissions \cite{Gough2011}, thus providing the potential for $reductions$ of atmospheric CO$_2$ concentrations, or for compensating other emissions. Moreover, the emissions factors of some power plants using conventional coal technologies are being be reduced by co-firing coal and biomass fuels. Meanwhile, liquid biofuels derived from biomass, such as ethanol and biodiesel, have the potential to replace oil-derived transport fuels with minimal changes in vehicle internal combustion engine technology and jet engines \cite{IPCCSRRENBiomass2011}. 

Biomass currently used for electricity and biofuel production largely originates from forestry and agricultural residues, and other forms of commercial or household mixed solid waste. Volumes of waste available could amount up to 100~EJ/y but are highly uncertain and not studied here\endnote{Waste amounts depend on efficiency of biomass use (such as in food or timber production) and therefore subject to significant changes depending on future policy. These are therefore difficult to model, as are their associated costs of production.} \cite{Hoogwijk2003,Smeets2007}. The larger share of bioenergy potential lies with the production of dedicated biomass crops. Global technical potentials for primary bioenergy range between 0 and 1550~EJ/y \cite{Wolf2003, Hoogwijk2005, Smeets2007, deVries2007,vanVuuren2009}. Bioenergy crops include perennial woody short rotation coppiced trees, such as willow, poplar or eucalyptus, perennial grasses such as miscanthus, elephant grass and switchgrass, starch rich crops such as wheat, corn, sugar beet and cane, and oil rich crops such as rapeseed and palm. Depending on their nature, they can be transformed into energy carriers by using, among many processes, combustion, gasification or anaerobic decay for electricity production, fermentation or the Fischer-Tropsch process for transport fuel production \cite{IEABioenergy2009}.

Biomass production for energy purposes makes use of agricultural land and thus may have a high opportunity cost. The technical potential that lies in agricultural land is large, but energy production from biomass is in direct competition for land with food production, a situation which has the potential to drive significant increases in world food prices \cite{Dornburg2010}. Following the approach which has been used recently by many authors (for instance in \cite{Hoogwijk2005, Hoogwijk2009,Smeets2007,Wolf2003}), the explicit assumption is taken in the present work that future bioenergy production uses no more than leftover land after the global food demand has been met, a premise that is difficult to justify in the absence of specific legislation and further investigations, but it avoids the complex problem of simulating food and biomass prices.\endnote{Bioenergy potentials could in principle be larger if global food demand is not met. However, it will not be lower if global food demand is indeed met. The problem of simulating food prices is complex as it involves modelling both local food markets, underreported in developing countries, and efficiency changes associated with changes in food prices.} Thus the bioenergy potential is obtained by subtracting from the total biomass potential the amount required by the food demand, based on population growth curves and dietary assumptions.

Hoogwijk $et$ $al.$ evaluated the use of the land at each point of a global grid yearly up to year 2100 using IMAGE~2.2, in which leftover agricultural land was termed ‘abandoned land’ \cite{Hoogwijk2005, Hoogwijk2009}. The reported cost-supply curves were observed in the present work to follow a distribution for nearly identical resource units using non-linear fits  of eq.~(\ref{eq:Dist1}) to their data. In addition to agricultural land, however, other types of geographical areas with lower productivities exist which can be used for particular bioenergy crops. These were labelled ‘rest land’ by Hoogwijk $et$ $al.$ and contribute a significant global technical potential. They were found to follow the distribution for hierarchical resources by using fits of eq.~(\ref{eq:Dist2}) to their data. Examples of such fits are described in section~S.2.5 of the supplementary material. Land use depends strongly on assumptions regarding world population, diet habits, global urbanisation and trade of agricultural products. The four main SRES scenarios, A1, A2, B1 and B2 (see \cite{IPCCSRES}), were taken as assumptions for all exogenous variables in these calculations, and results are presented for each. Large differences arise between scenarios, with technical potentials ranging between around 302~EJ/y for the A2 scenario to 676~EJ/y for the A1 scenario, which result in large uncertainties for values of the global biomass technical potential. This work has however been revisited more recently taking into account additional assumptions of future water scarcity and land degradation by van Vuuren $et$ $al.$, yielding lower estimates ranging from 120 to over 300 EJ \cite{vanVuuren2009}. Other ranges have been estimated using different methodologies, with more pessimistic projections of 0 to 648~EJ/y by Wolf $et$ $al.$, and optimistic values of 367 to 1548~EJ/y by Smeets $et$ $al.$ \cite{Wolf2003,Smeets2007}. The low end of the range given by Wolf $et$ $al.$ stems from high projected food demand and low agricultural productivity, while the high end is due to mostly vegetarian diets and high productivity. Meanwhile, the high end of the range of Smeets $et$ $al.$ originates from `super high' agricultural productivity, high availability of the land and landless animal production systems.

Figure~\ref{fig:GlobalRenCSC} presents the global economic potential of bioenergy in terms of primary energy before conversion to electricity or liquid biofuels (derived from the data in \cite{Hoogwijk2005, Hoogwijk2009, Smeets2007, Wolf2003}), using both abandoned and rest land, but includes also a small component from bagasse of 3~EJ/y. Four cost-supply curves are given, calculated by Hoogwijk $et$ $al.$ for the A1, A2, B1 and B2 SRES scenarios in 2050, shown as solid curves\endnote{Values projected for the year 2050 were used since amounts of land potentially available in the future are expected to increase due to increased agricultural efficiency.} \cite{Hoogwijk2009}. A value near zero was taken for the lower boundary of the uncertainty range, consistent with the low end of the range calculated by Wolf $et$ $al.$ while the high end of the range calculated by Smeets $et$ $al.$ was taken for the upper boundary of the uncertainty range \cite{Wolf2003, Smeets2007}. For a decarbonisation scenario, the cost-supply curve derived for the B1 SRES scenario was considered the most probable cost-supply curve, but for other types of scenarios, choices of curves consistent with particular working assumptions should be made. 

\subsection{Ocean energy sources}

The term ocean energy denotes renewable energy produced using seawater as a resource, where unlike for wind energy or hydroelectricity, not only the kinetic energy of seawater can be used to produce electricity, but also temperature gradients in the ocean and salinity differences near river mouths.  Using ocean energy as a general classification type, it can be divided into four main sources of energy \cite{IEAETSAP2010Marine, IPCCSRRENOcean2011}:

\begin{itemize}
\item {\bf Wave energy}, driven by transfers of energy from the wind to the surface of the ocean,
\item {\bf Tidal energy}, driven by the rise and fall of sea levels due to gravitational forces (tidal range) and the resulting water currents,
\item {\bf Ocean Thermal Energy}, driven by temperature gradients between upper and lower ocean layers,
\item {\bf Salinity Gradient energy}, derived from salinity gradients between ocean and fresh water at the mouths of rivers. 
\end{itemize}

Section~S.3.6 of the supplementary material provides a review of theoretical potentials for these sources, resulting in a total that could be as high as 600~EJ/y. Technical potentials however are much lower and uncertain, since the current development status for ocean energy technologies excluding tidal is preliminary, and cost data is in some cases unavailable. Specific geographical and configurational requirements for tidal and salinity gradient technologies involves, as it is the case for hydroelectricity, calculating the technical potential by summing the potential values from a large number of individual studies. Such studies have not been performed exhaustively on a wide scale yet. Meanwhile, wave and ocean thermal are based onto global extrapolations carried out using physical measurements. Global energy potentials calculated in various studies are given in table~\ref{tab:Ocean}, and additional details are given section~S.3.6 of the supplementary material.

\begin{table}[h]\footnotesize
\begin{center}
		\begin{tabular*}{1\columnwidth}{@{\extracolsep{\fill}} l|r r r l }
			\hline
			Technology & Min. & Mode & Max. & Study\\
					&EJ/y	&EJ/y	&EJ/y	&\\
			\hline
			\hline
			Wave	&6.3	&	&	&\cite{Sims2007}\\
			Energy	&	&18.9&	&\cite{WEC1994}\\
					&	&	&65	&\cite{WEA2001}\\
			\hline
			Tidal		&1.8	&	&	&\cite{Hammons1993}\\
			Energy	&	&3.6	&	&\cite{Hammons1993}\\
					&	&	&7.2	&\cite{WEC1994}\\
			\hline
			\bf Total	&\bf 8.1&\bf 22.5&\bf 72.2&\\
			\hline
			Ocean	&3.2	&	&	&\cite{Charlier1993}\\
			Thermal	&	&32	&	&\cite{Charlier1993}\\
			Energy	&	&	&85	&\cite{Nihous2007}\\
			\hline
			Salinity	&5.8	&	&	&\cite{Skramesto2009}\\
			Gradient	&	&7.2	&	&\cite{Krewitt2009}\\
			Energy	&	&	&83	&\cite{Cavanagh1993}\\
			\hline
			\bf Total		&\bf 17.1&\bf 62.8&\bf 240.2&\\
			\hline
		\end{tabular*}
	\caption{Technical potentials for different types of ocean energy used to define the cost-supply curve. The uncertainty ranges are defined using the $Min$ and $Max$ values, while $Mode$ represents the most probable value.}
	\label{tab:Ocean}
\end{center}
\end{table}
Energy potentials for ocean thermal and salinity gradient energy are theoretical and highly uncertain, and no reliable cost estimates were found. These types of resources were therefore not included in the present calculations for the most probable cost-supply curve, due to the risk of generating misleadingly optimistic potentials given the lack of reliable information.\endnote{The costs for ocean thermal and salinity gradient energy systems are likely to be much higher than those of tidal and wave energy. This would result in a piecewise cost-supply curve featuring an additional step at high costs.} Wave and tidal systems are better established. Therefore, using cost values obtained from the IEA and the ETSAP, a cost-supply curve for ocean energy based on an aggregation of separate cost-supply curves for wave and tidal energy was produced, and is given in Figure~\ref{fig:GlobalRenCSC}, with a small technical potential of 22.5~EJ/y  \cite{IEAETSAP2010Marine,IEAProjCosts}. The  lower and upper boundaries of the uncertainty range were obtained from the extremal values  of 8 and 72~EJ/yr respectively given in table~\ref{tab:Ocean}.

\section{Stock resources \label{sect:Stock}}
\begin{figure*}[p]
	\begin{minipage}[t]{1\columnwidth}
		\begin{center}
			\includegraphics[width=1\columnwidth]{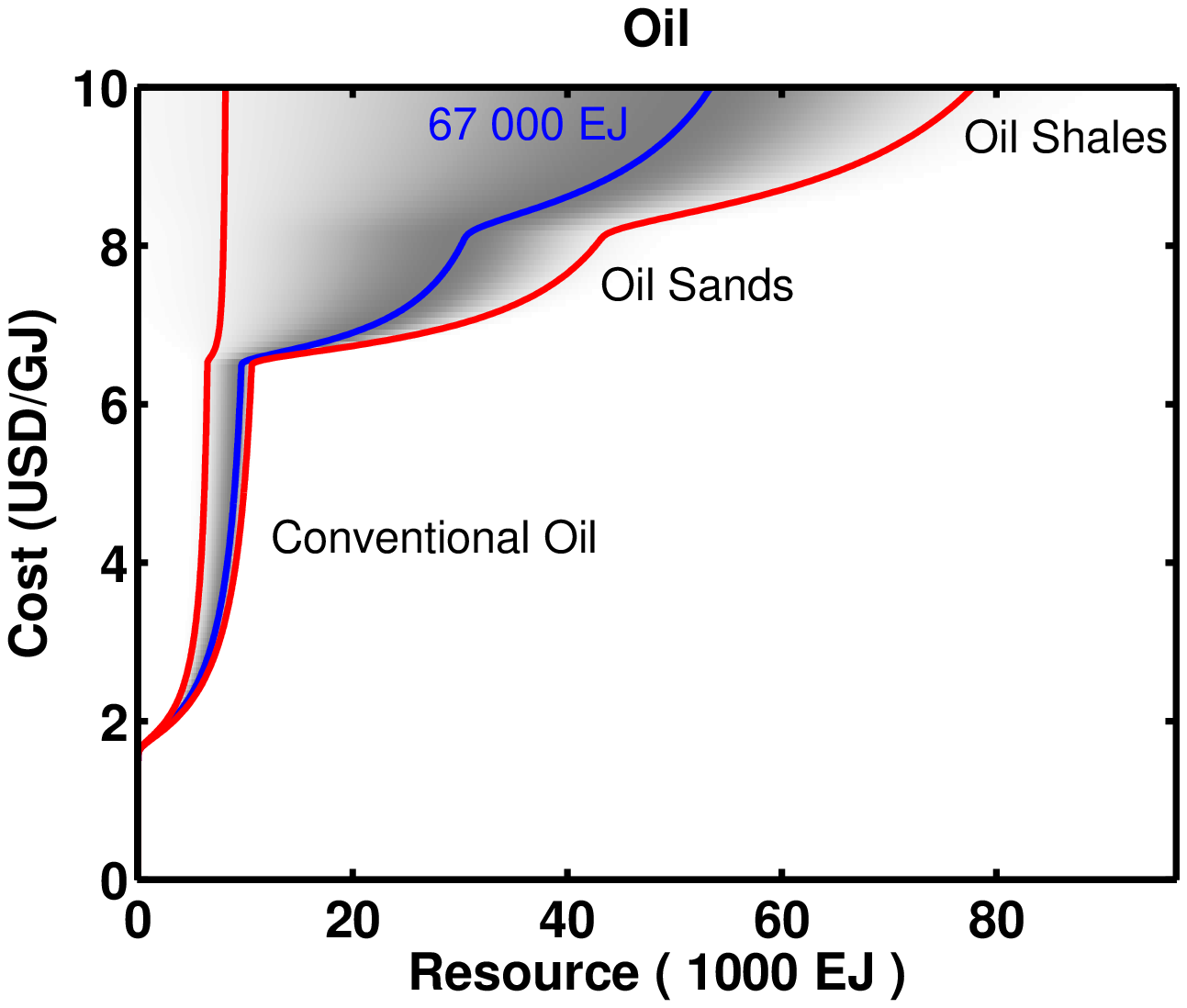}
			\includegraphics[width=1\columnwidth]{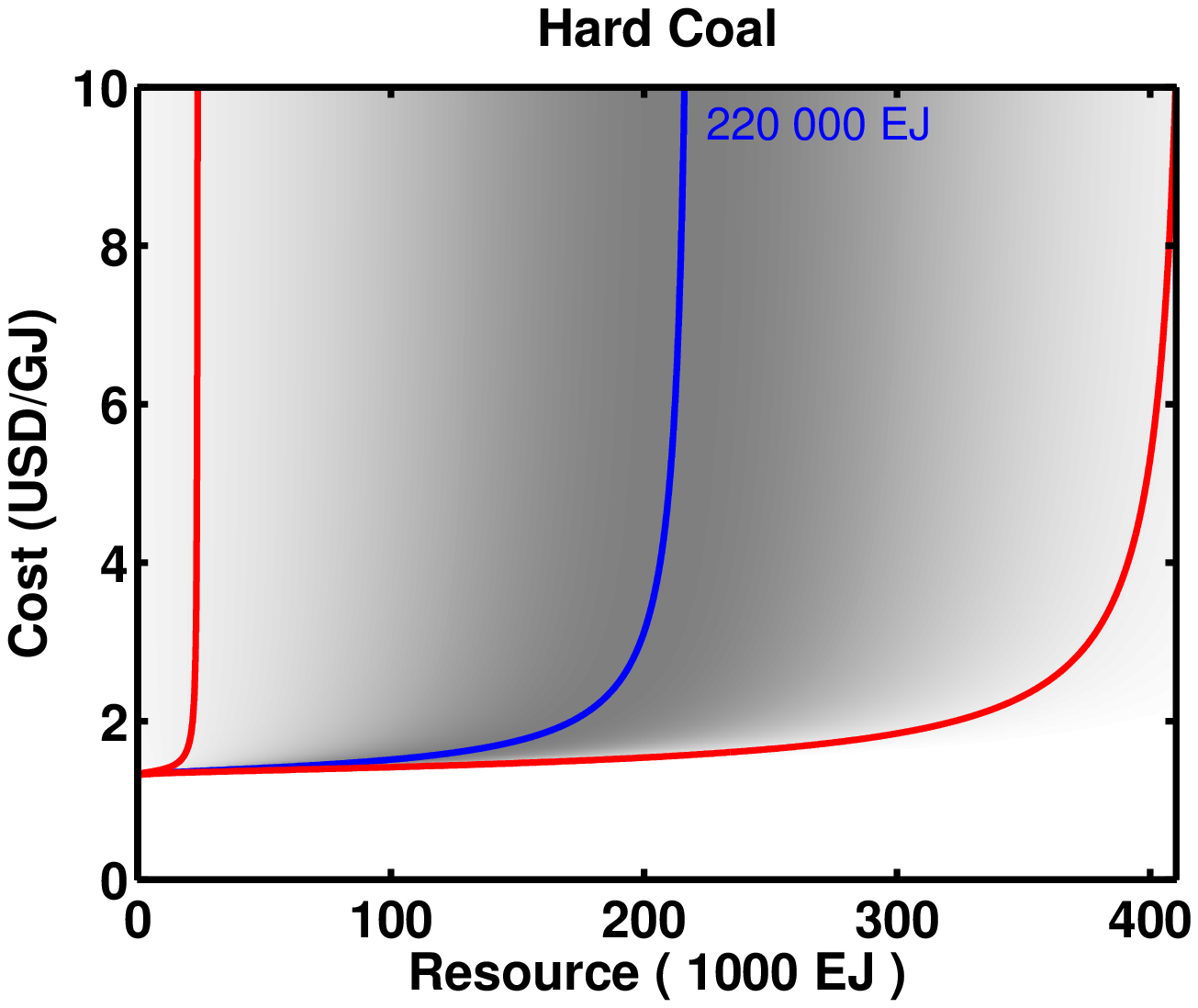}
			\includegraphics[width=1\columnwidth]{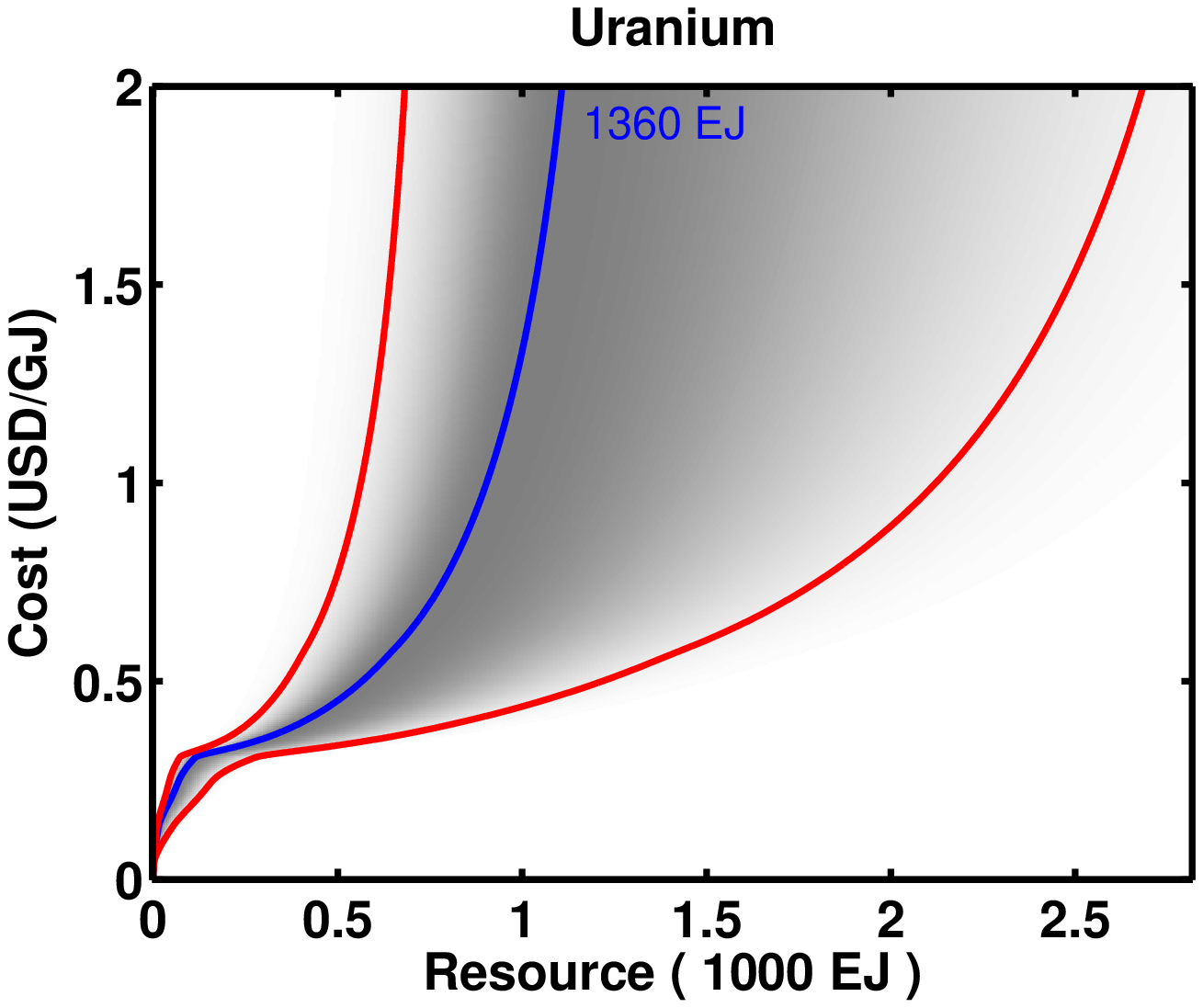}
		\end{center}
	\end{minipage}
	\hfill
	\begin{minipage}[t]{1\columnwidth}
		\begin{center}
			\includegraphics[width=1\columnwidth]{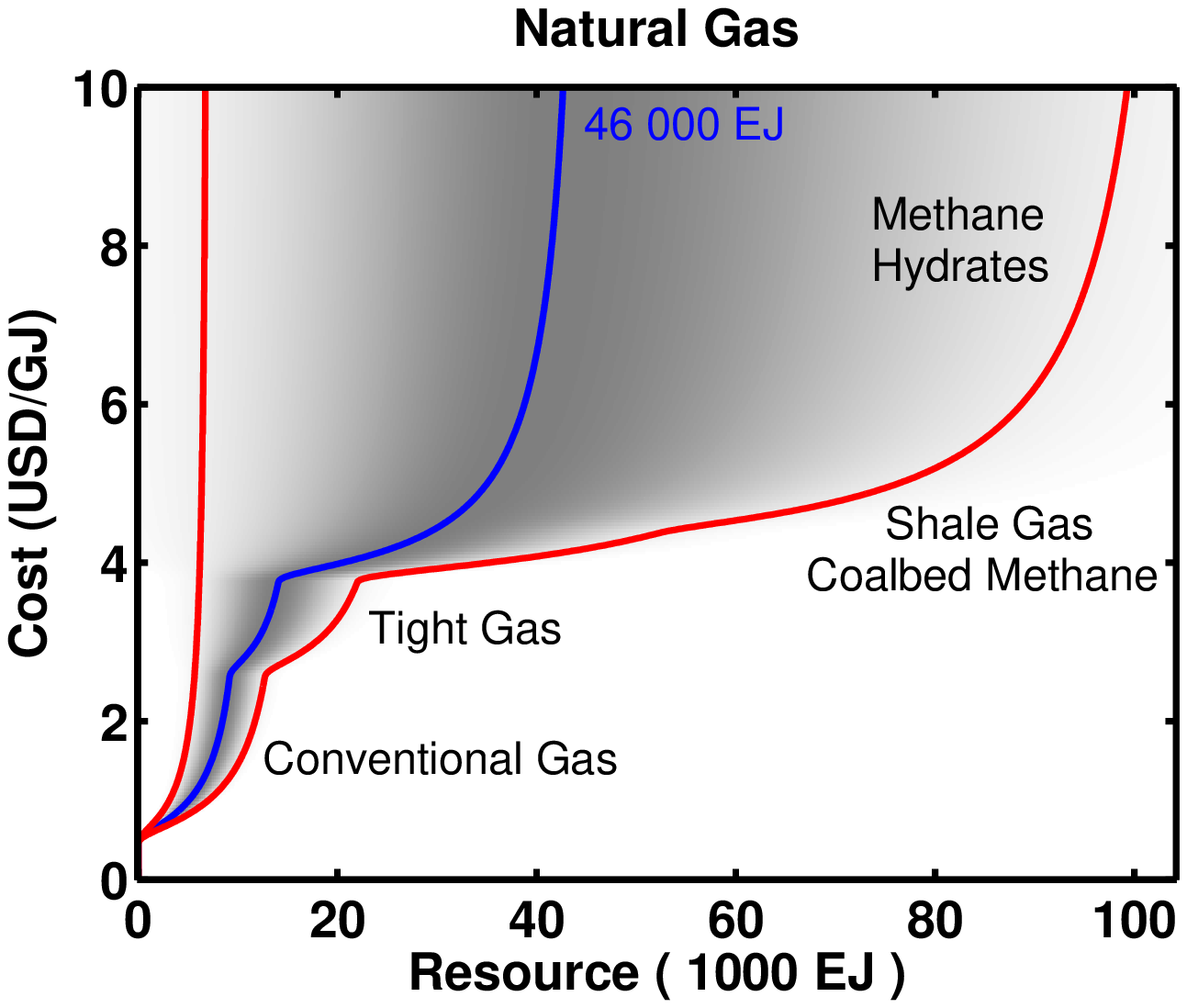}
			\includegraphics[width=1\columnwidth]{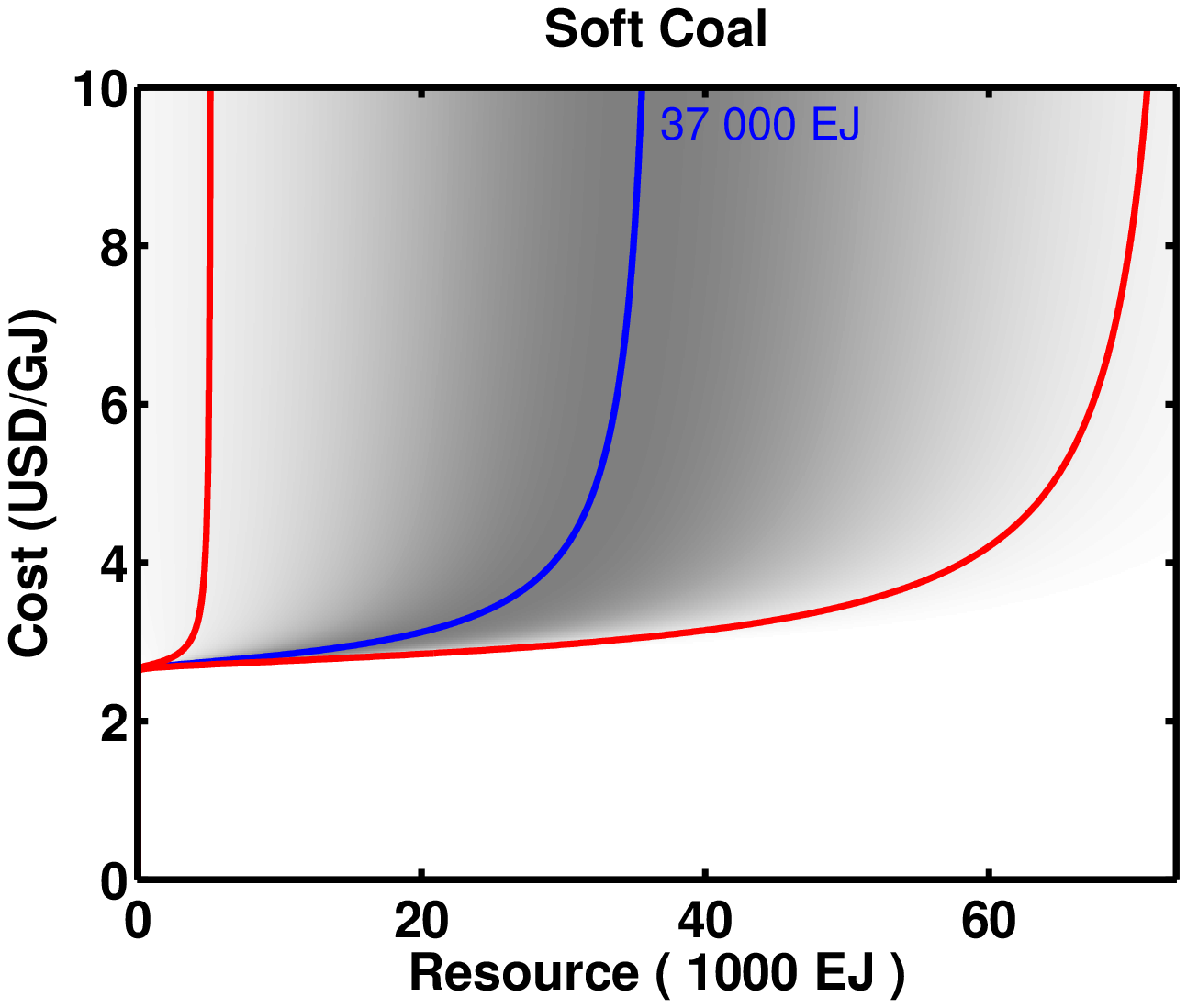}
			\includegraphics[width=1\columnwidth]{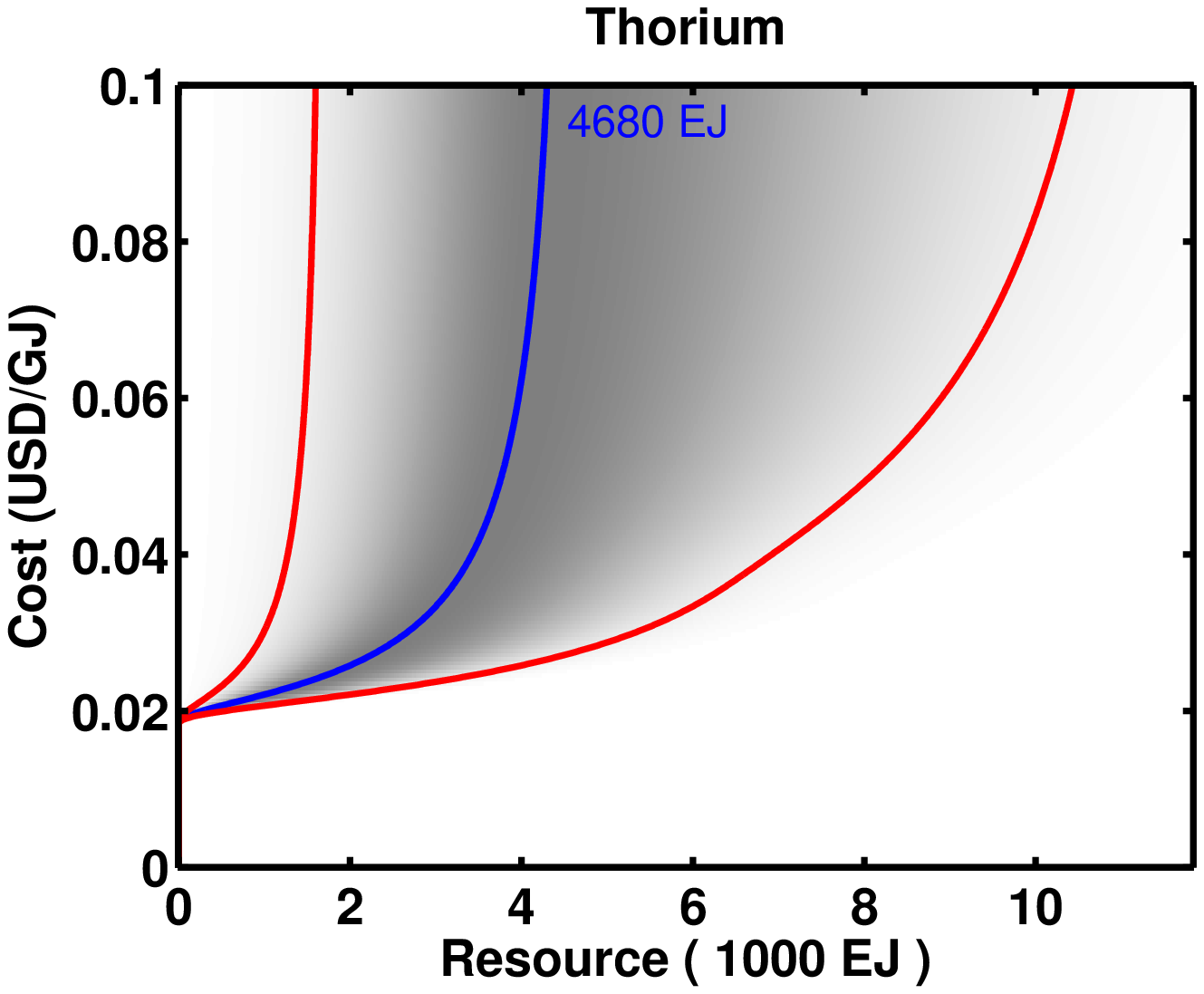}
		\end{center}
	\end{minipage}	
	\caption{Cost-supply curves for fossil and nuclear resources, including oil, gas, hard coal, soft coal, uranium and thorium. Hard coal includes anthracite and bituminous coal, defined as coal with a calorific content above 16~500~kJ/kg. Soft coal corresponds to sub-bituminous coal and lignite, and includes all coal with a calorific content lower than 16~500~kJ/kg.}
	\label{fig:GlobalStockCSC}
\end{figure*}

\subsection{Fossil fuels}

As it occurs with all types of exhaustible natural resources, fossil fuel resources and reserves are known to continuously expand, even though they are gradually consumed. This is due to periodic resource discoveries and improvements in the methods of extraction. Therefore, what is considered economical to extract changes every year. Reserves are distinct from resources, the former referring to the resources that are known to exist with almost complete certainty and to be economical to extract, while the latter refers to those which are thought to exist with various degrees of confidence, and those currently thought too expensive to extract. As technological improvements and additional knowledge affect the economics of different methods of extraction, there is a flow from resources towards reserves, and thus reserves expand \cite{McKelvey1972, Rogner1997}. Meanwhile, discoveries continuously add to resources. As prospecting activities for hydrocarbon sources remains very active, this makes the production of cost-supply curves more difficult than for renewables, since is at best a snapshot in time of what is known to exist and recoverable with current technology.  

In order to assess global energy potentials, it is nevertheless necessary to explore cost-supply curves for fossil fuels, even if they are derived from current knowledge, and therefore expected to change in the future. It is unlikely that fossil fuel resources turn out smaller than what is currently expected to exist. On the contrary, it is probable that they turn out significantly larger as methods of extraction are devised for types of occurrences which were until recently not thought possible to use, such as gas hydrates or oil shales. Therefore, the cost-supply curve uncertainty ranges are highly asymmetric. The associated extraction costs, which increase as low cost conventional sources are depleted, nevertheless decrease due to technological improvements, and it is therefore not immediately obvious whether costs are likely to go up or down in the future. 

Global cost-supply curves have been calculated previously by Rogner \cite{Rogner1997}. These results have been used extensively by the energy modelling community; however they are becoming increasingly outdated. This section provides an update to the work of Rogner, but using an approach emphasising uncertainty, and thus, following the spirit of the current treatment of renewable resources, in opposition to the approach of Rogner, the results of this section should be interpreted as ranges rather than specific values. 

The economic potentials of fossil fuels are given in Figure~\ref{fig:GlobalStockCSC}, showing in order liquid hydrocarbons, gaseous hydrocarbons, hard and soft coal, the last two being classified using their calorific content.\endnote{Hard coal includes anthracite and bituminous coal, while soft coal includes sub-bituminous coal and lignite, the last two having lower calorific contents than the first two. The limiting calorific value used to separate the two categories is of 16~500~kJ/kg \cite{BGR2010}.} For oil and gas, different types of occurrences considered in this assessment are indicated with text. These are associated with independent distributions of the hierarchical type, aggregated to produce composite curves. Due to the wide use and global diffusion of fossil fuel extraction technology, extraction cost ranges were assumed to be the same for all regions of the world. In the case of coal, less information was found over differing types of mines and associated costs, and single distributions were used, where costs were assumed to vary little with the amount extracted. This is unlikely to matter in the long run given the very large scale of the resource, and limited expectations of its depletion.

Uncertainty ranges were determined using resource classifications, and in some cases where this is unavailable, their nature. Oil occurrences obtained from the World Energy Council and the German Federal Institute for Geosciences and Natural Resources are classified as either reserves or resources, with the exception of oil shales, which are as a whole considered resources only \cite{WEC2010,BGR2010}. Four types of oil resources were considered, conventional (crude) oil, oil shales, oil sands and extra-heavy oil. Cost ranges were obtained from the IEA \cite{IEAWEO2008}. Gas occurrences follow a similar trend, but a larger number of types of resources were considered: conventional gas \cite{BGR2010}, shale gas \cite{EIA2011}, tight gas \cite{BGR2010,WEA2001}, coalbed methane \cite{Boyer1998} and methane hydrates \cite{Boswell2011}. The respective cost ranges were obtained from the ETSAP \cite{IEAETSAP2010}. Large amounts of methane are known to exist dissolved in aquifers \cite{WEA2001}, but were not included due to the lack of reliable data. Information for coal was derived from a mixture of data \cite{WEC2010,BGR2010}. Complete details on the methodology underlying these curves, as well as region specific data tables can be found in the supplementary material.

\subsection{Fissile materials}

Five sources of fissile materials for nuclear reactors are known to exist. These are enumerated in order of cost. The first comes from stocks of highly concentrated $^{235}$U (uranium) or $^{239}$Pu (plutonium) originating from decommissioned nuclear arsenals diluted with $^{238}$U. The second source is lightly enriched $^{235}$U/$^{238}$U produced from mined natural deposits. The third originates from U and Pu recovered from spent fuel (using the PUREX process \cite{Bonche2005}). The fourth source is thorium (Th) using the $^{232}$Th/$^{233}$U fuel cycle. The fifth source is U which occurs in very low concentrations in  seawater. Producing a cost-supply curve involves creating a scenario for the nuclear sector, and requires careful consideration of uncertainty. Additionally, if ingenious use of fast reactors is invoked for the future, fuel efficiencies of up to 50 times larger could be obtained, altering dramatically these expectations. 

In order to construct a cost-supply curve for U and Pu, the nuclear industry was assumed to continue to use current methods and thermal reactors, and therefore, only fuels originating from naturally mined U and from nuclear arsenals were considered. Many authors stress that deposits of Th worldwide are three times larger than those of U \cite{Bonche2005, Sinha2006, AbuKhader2009, Suess1956}. However, less efforts at prospecting for Th ore have been carried out and as a consequence, the current reasonably assured reserves of Th, in tonnes of natural Th, are lower than those of U \cite{IAEA2009}, a situation which is likely to change if interest in Th grows. The nuclear fuel cycle for Th being more efficient than that of thermal reactors based on U, it leads to larger amounts of energy per tonne of natural Th and thus leads to lower fuel costs per unit of energy. Costs only include the extraction costs given by the IAEA, without the inclusion of enrichment or transformation components.

Detailed resource data from the IAEA for naturally occurring U and Th were used to construct two cost-supply curves and associated uncertainty ranges \cite{IAEA2009}. The data are classified into four levels of certainty and four cost ranges. Such resources generally increase naturally in size with increasing costs of extraction, as well as with uncertainty, an effect produced by the hierarchical ordering of natural resource consumption and by the decreasing amount of effort which has been spent on prospecting activities for resources more and more difficult to exploit. For the conversion of resources from tonnages to energy values, an average conversion efficiency for thermal U reactors of 159~TJ/t was used, determined from the 2008 electricity production of 2611.1~TWh from a global fleet capacity of 273.7~GW, with a capacity factor of 80\%, which used 59~065~t of natural uranium \cite{IAEA2009}. Meanwhile, the burnup rate for Th reactors was derived from the value of 24~000~MWd/t reached by the experimental Indian model \cite{Sinha2006}, equal to about 2100~TJ/t. Panels $e)$ and $f)$ of Figure~\ref{fig:GlobalStockCSC} present the resulting global economic potentials for U  and Th. Uncertainty ranges for U were obtained by considering only reasonably assured reserves (RAR) for the lower boundary of the uncertainty range, RAR and inferred reserves for the most probable cost-supply curve, and all of the RAR, inferred, prognosticated and speculative resources for the upper boundary of the uncertainty range.

The uncertainty ranges are highly asymmetric due to the tendency of the size of speculative resources to increase with the level of uncertainty. It is observed that in terms of energy, reserves of Th are larger than those of U, and that these are also less expensive per unit of energy, due to the higher burnup rate of the Th system. It must be emphasised that U resources could, in principle, be used with much higher burnup rates, were fast reactors to be deployed globally. The resources of U do not include seawater U, as data over these are scarce and highly speculative.\endnote{U is present dissolved at very low concentrations in seawater (3-4 ppb, from \cite{IAEA2009}), giving rise nevertheless to large amounts of U given the size of the terrestrial body of seawater. Water turnover due to currents is very slow however, making it highly speculative whether a significant portion of seawater U can be recovered, and the costs involved are very high \cite{Bonche2005}.} Finally, it is to be noted that the fuel component of the levelised cost for nuclear reactors is very small compared to the investment costs, which results in a very small influence of the fuel costs onto the decision-making, unless nuclear resources become depleted.

\section{Summary of energy resources}

\begin{table}[hbt]\footnotesize
\begin{center}
		\begin{tabular*}{1\columnwidth}{@{\extracolsep{\fill}} l l l|r|r r r r }
			\hline
			Resource &  &  & Use &\multicolumn{4}{ l }{Technical Potential}\\
			\hline
			Name &Type& Dist. & EJ/y & L & M &U&Units\\
			\hline
			\hline
			Wind		&Flow	&Hierarch.	&.72		&72	&346	&2257	&EJ/y	\\
			Solar	&Flow	&Identical		&.04		&1340&3384	&14778	&EJ/y\\
			Hydro	&Flow	&Hierarch.	&12		&12	&66		&148	&EJ/y\\
			Geotherm.&Flow	&Hybrid		&0.23	&4	&36		&111	&EJ/y\\
			Biomass	&Flow	&Hybrid		&51		&0	&447	&1548	&EJ/y\\
			Ocean	&Flow	&Hierarch.	&.002	&8	&23		&72		&EJ/y\\
			Oil		&Stock	&Hierarch.	&170	&9	&67		&98		&10$^3$EJ\\
			Gas		&Stock	&Hierarch.	&109	&7	&46		&106	&10$^3$EJ\\
			Hard Coal	&Stock	&Hierarch.	&\multirow{2}{*}{139}&24	&220	&419	&10$^3$EJ\\
			Soft Coal	&Stock	&Hierarch.	&	&5	&37		&75		&10$^3$EJ\\
			Uranium	&Stock	&Hierarch.	&30	&0.83&1.36	&3.43	&10$^3$EJ\\
			Thorium	&Stock	&Hierarch.	&-	&1.74	&4.68&12.27	&10$^3$EJ\\
			\hline
		\end{tabular*}
	\caption{Summary table for all energy resources. $Stock/Flow$ indicate whether resources are renewable flows or stocks. $Hierarch./Identical/Hybrid$ identifies the type of statistical distribution assigned. $Use$ refers to current yearly consumption of these resources. $L$ indicates the lower boundary of the uncertainty range. $M$ indicates the most probable technical potential. $U$ indicates the upper boundary of the uncertainty range.}
	\label{tab:Summary}
\end{center}
\end{table}

Table~\ref{tab:Summary} provides a summary of all types of global energy resources, classified by type (renewable flows or stocks), to which a type of statistical distribution it is assigned (for hierarchical or nearly identical resources, or a hybrid mixture of both), along with technical potential values. The potential values are given with their lower and upper boundaries of the uncertainty ranges. For comparison, current consumption of these resources is given based on data from the IEA \cite{IEAWEO2010}. Note that biomass resources are expressed in terms of primary energy, which become smaller when converted into electricity according to the efficiency of transformation.

\section{Conclusion}

This paper presents an assessment of global economic energy potentials for all major energy resources, those with a potential larger than 10~EJ/y. These were given in the form of cost-supply curves, adding an economic structure to energy potentials, and therefore providing them an unambiguous definition. Additionally, these were provided using a probabilistic construction that allows a simple representation of uncertainty. The curves were calculated using assumptions over the cost distribution of resources using functional forms based on statistical properties of resource types. The set of energy potentials include six types of renewable energy sources, wind, solar, hydroelectric, geothermal, biomass and ocean energy, as well as four types of fossil fuels, oil, gas, hard and soft coal, and two nuclear materials, uranium and thorium. While the potentials for renewable resources were determined predominantly based onto an extensive review of the literature, potentials for stock resources were determined directly using resource and reserve assessment data. 

The cost-supply curves calculated in this work were produced in order to be used by the global energy modelling community, for the purpose of constraining models of the energy sector in order to produce realistic scenarios of future energy use. It is hoped by the authors that this work will supersede outdated studies currently used and provide a consistently calculated update for all types of energy resources. In particular, the large set of regional cost-supply curves underlying the aggregate curves presented in this paper form the core of a new model for natural resource use and depletion for the global E3 model E3MG, to be used through the family of technology models FTT. Other regional aggregates can be provided by the authors.

Resource assessments, however, change continuously as new information becomes available, and as new and more sophisticated studies are carried out. Therefore, while the potentials presented here may become outdated in the future, however the methodology presented will not. The simple and robust theoretical framework and computational methodology presented here can be very useful for reporting natural resource assessments, and enable direct use in models of energy systems. It would therefore be appropriate to reuse this approach with new data as they become available, either by rescaling the curves given (the values for the $A$ parameters) or by recalculating new sets of parameters. In either cases, this approach provides a consistent and general framework for limiting all types of natural resources in energy models.

\section*{Acknowledgements}
The authors would like to acknowledge Dr T. S. Barker and T. Hanaoka for guidance and support, as well as A. Anger, H. Pollitt, P. Summerton and P. Bruseghini for highly informative discussions. This work was supported by the Three Guineas Trust,  Conicyt (Comisi\'on Nacional de Investigaci\'on Científica y Tecnol\'ogica, Gobierno de Chile) and the Ministerio de Energie, Gobierno de Chile. 

\theendnotes



\newpage
\onecolumn
\setcounter{section}{0}
\renewcommand*\thesection{S.\arabic{section}}
\renewcommand*\thesubsection{S.\arabic{section}.\arabic{subsection}}
\part*{Supplementary Material}
\section{Introduction and use of the supplementary material}

This part of the work aims to complement the main paper by providing all details that could be required by researchers interested in either:
\begin{itemize}
\item Verifying the methodology or reproducing the results
\item Rebuilding this database for use in a particular model of energy systems
\end{itemize}
A high amount of care was put into summarising compactly all relevant information in this part of the work in order to make this possible. Mathematical details underlying the calculations given in the main text are provided. Additionally, tables of data and parameters are given for a chosen set of  world regions, which may not necessarily correspond to the particular divisions used by other research groups. It is however impractical to provide larger tables involving all countries of the world, even though such tables exist underlying this work  (190 countries exist in our database). For more information, the authors may be contacted at the address provided.

Natural resource assessments are performed continuously and what is known of global natural resources changes continuously. Therefore the cost-supply curves in this work may become outdated. This will happen due to three processes: firstly the total amount of resources might change (parameter $A$), secondly the scaling of costs may change through inflation (parameter $B$) and the costs of technology may reduce due to learning-by-doing (parameter $C_0$). However, the structure of the cost-supply curves, or shape, will not change. These parameters should be simply scaled to new values and the results will still be valid.

\section{Derivation of distribution functions and cost-supply curves \label{sect:SDist}}
\numberwithin{table}{section}
\numberwithin{figure}{section}

\subsection{Distribution function for the hierarchical type of resources}

Hierarchical resources have an exponential energy distribution in productivity space:
\beq
f(\nu)d\nu =  \left\{ \begin{array}{cc} {A \over \sigma} e^{-{\nu \over \sigma}}d\nu  & \nu > \mu \\ 
0 & \nu \leq \mu \end{array} \right .,
\eeq
where $A$ is the technical potential and $\sigma$ is the half width of the function.
This function is required in cost space, and the equation connecting cost $C$ to productivity $\nu$ is 
\beq
C = {C_{var}\over \nu} + C_{0}.
\eeq
where $C_{var}$ corresponds to costs per unit of effort or per unit of resource producing items such as the rent of the land (in \$/km$^2$), bore hole depth (\$/km), dam size and type, mine depth, etc. The ratio of $C_{var}$ to $\nu$ has units of \$/GJ. $C_{0}$ is the sum of fixed costs, such as capital investments, transport or transformation costs, etc, in \$/GJ. 

The density productivity interval must  be transformed into a cost interval:
\beq
dC = -{C_{var}\over \nu^2}d\nu, \quad \Rightarrow d\nu = -{C_{var} \over (C-C_{0})^2}dC.
\eeq
Using the cost scaling parameter $B = C_{var}/\sigma$, the distribution becomes:
\beq
f(C)dC =  \left\{ \begin{array}{cc} {A B \over (C-C_{0})^2} e^{-{B \over C-C_{0}}}dC  & C > C_{fixed} \\ 
0 & C \leq C_{0} \end{array} \right ..
\eeq

\subsection{Distribution function for nearly identical resources}

In the case of nearly identical resources, there is no convenient exact analytical form that can be derived from eq.~\ref{eq:DistProd2} of the main paper. However, the form given in eq.~\ref{eq:Dist2} works extremely well when compared to data (see section~\ref{sect:SFits}), and can be derived from eq.~\ref{eq:DistProd2} through a simple approximation. 

Nearly identical energy producing resources, such as land plots, are assumed truly identical, and therefore have a potential situated at a single value of productivity $\nu = \mu$,
\beq
n(\nu)d\nu = N\delta(\nu-\mu)d\nu,
\eeq
where $n$ is a density of energy producing land area, while $N$ is the total energy producing land area (in km$^2$) and the function $\delta(\nu)$ is the Dirac delta function.\footnote{The Dirac delta function is defined such that $\int_{-\infty}^{\infty}\delta(x-a)f(x)dx = f(a)$.} Without any additional reductions in productivity, the total amount of energy that can be obtained from these land resources would be their area times the productivity:
\beq
A =  \int_0^\infty \nu N\delta(\nu-\mu)d\nu = N\mu.
\eeq
Unit land areas have a suitability factor, however, that reduces their productivity below the maximum value of $\mu$ by a small amount $\epsilon$ with a probability $P$. The probability for the reduction in productivity is assumed to be normally distributed around zero but positive $\epsilon$, with standard deviation much smaller than the average productivity $\sigma << \mu$:
\beq
P(\epsilon) d\epsilon = \left\{ \begin{array}{cc}{2\over \sqrt{2 \pi} \sigma } e^{-{\epsilon^2 \over 2 \sigma^2}}d\epsilon & 0 \leq \epsilon \leq \mu\\
0 & otherwise \end{array} \right .,
\eeq
where the reduction in productivity $\epsilon$ must be less than the maximum value $\mu$. The distribution of resources must be calculated by summing over all reduction values $\epsilon$ given their probability $P(\epsilon)$:
\beq
n'(\nu)d\nu = \int n(\nu) P(\epsilon)d\nu d\epsilon = \int_0^{\mu} {2 N \over \sqrt{2 \pi} \sigma} e^{-{\epsilon^2 \over 2 \sigma^2}}\delta(\nu-\mu+\epsilon)d\nu d\epsilon.
\eeq
This is a convolution and can be seen as a sum of several Dirac Delta functions centred at slightly reduced values of productivity, $\mu-\epsilon$, with probability $P(\epsilon)$, instead of one Dirac Delta function centred at $\mu$ with probability 1. The total amount of energy that can be obtained from each plot of land corresponds to its area times its productivity. Thus, the productivity distribution of energy production potential leads to eq.~\ref{eq:DistProd2} of the main paper:
\beq
g(\nu) d\nu = \nu n'(\nu)d\nu = \left\{ \begin{array}{cc}{2 N\over \sqrt{2 \pi} \sigma } \nu  e^{-{(\nu-\mu)^2 \over 2 \sigma^2}}d\nu & \nu \leq \mu\\
0 & \nu > \mu \end{array} \right ..
\eeq
This function is required in cost space, and the equation connecting cost to productivity is 
\beq
C = {C_{var}\over \nu} + C_{fixed},
\eeq
where $C_{var}$ corresponds to the rent of the land (in \$/km$^2$), while $C_{fixed}$ is the sum of fixed costs (in \$/GJ). The productivity is situated very near the value of $\mu$, since the variations of productivity are small and $\sigma << \mu$. $\nu$ can be rewritten as a small variation around $\mu$, i.e. $\mu - \Delta$:
\beq
C = {C_{var}\over \mu - \Delta} + C_{fixed} = {C_{var}/ \mu \over 1 - \Delta/\mu} + C_{fixed} \simeq {C_{var}\over \mu}\left(1+{\Delta \over \mu}\right) + C_{fixed} = {C_{var}\over \mu^2}\left(2\mu - \nu \right) + C_{fixed},
\eeq
which is the crucial approximation, and using $C_0 = C_{var}/\mu + C_{fixed}$, 
\beq
\nu  = (C_0 - C){\mu^2\over C_{var}} + \mu,
\eeq
where $C_0$ is defined as the sum of fixed costs plus $C_{var} / \mu$, the total cost at the maximum productivity value. Since $dC = -C_{var}/ \mu d\nu$, the density can be rewritten as
\beq
g(C)dC = \left\{ \begin{array}{cc}{2 N\over \sqrt{2 \pi} B } \left( (C_0-C){\mu^2\over C_{var}} + \mu \right) e^{-{(C-C_0)^2 \over 2 B^2}}dC & C > C_0\\
0 & C \leq C_0 \end{array} \right .,
\eeq
where the cost scaling parameter $B$ is defined as $C_{var} \sigma / \mu^2$. This can be rewritten further as 
\beq
g(C)dC = \left\{ \begin{array}{cc}{2 N\over \sqrt{2 \pi} B } \mu \left( {C_0-C \over B} {\sigma \over \mu} + 1 \right) e^{-{(C-C_0)^2 \over 2 B^2}}dC & C > C_0\\
0 & C \leq C_0 \end{array} \right ..
\eeq
The value of $C$ cannot be below $C_0$ by definition, but also, the factor $\exp[-(C-C_0)^2/2B^2]$ decreases rapidly to zero as $C-C_0$ becomes larger than $B$. Therefore, the value of the term $(C_0-C)/B$ is less than one wherever a significant potential of energy exists. Since $\sigma$ is much smaller than $\mu$, this results with
\beq
{\sigma \over \mu}{C-C_0 \over B} << 1,
\eeq
and therefore the distribution becomes
\beq
g(C)dC = \left\{ \begin{array}{cc}{2 A\over \sqrt{2 \pi} B } e^{-{(C-C_0)^2 \over 2 B^2}}dC & C > C_0\\
0 & C \leq C_0 \end{array} \right .,
\eeq
where now $A = N \mu$ is approximately the total energy potential. 

This is the strict region of validity of the expression given in \ref{eq:Dist2} of the paper. In numerical terms, it is empirically found that the rigidity of these rules can be relaxed and the validity extended. For example, the distribution in productivity space can actually have a tail towards higher values or have a large value for $\sigma$, and this does not significantly alter the goodness of fit of the function in cost space. 

\subsection{Cost-supply curve expressions}

From the distributions $f(C)dC$ and $g(C)dC$, cumulative distributions $N(c)$ can be derived. For hierarchical resources, this results in 
\beq
N(C) = A\,e^{\left(-{B \over C-C_0}\right)},
\eeq
while for nearly identical resources this is
\beq
N(C) = A\,\text{erf}\left({(C-C_0)\over \sqrt{2}B}\right),
\eeq
where `erf' is the error function.

The cost-supply curves are the inverse of these functions, which respectively give
\beq
C(N) = {-B \over \ln\left({N\over A}\right)} + C_0
\eeq
and
\beq
C(N) = \sqrt{2}B\, \text{inverf}\left({N \over A}\right) + C_0,
\eeq
where `inverf' is the inverse error function.

\subsection{Parameterisation formulas}

In order to define a particular distribution or cost-supply curve, the parameters $A$, $B$ and $C_0$ are required, while the data available usually involves the total technical potential of the resource and a fraction of this considered to exist at costs situated between two values, as well as the current level of use of the resource. Assuming that two points of the curve are known, this may be expressed as two quantities $Q_1$ and $Q_2$ which occur at two cost values $C_1$ and $C_2$. These quantities can be expressed as fractions of the total technical potential, $\delta_1$ and $\delta_2$, the latter corresponding to the parameter $A$. In the case of the distribution for hierarchical resources, the values of $B$ and $C_0$ are the following:
\beq
C_0 = {C_2\ln \delta_2 - C_1\ln \delta_1 \over \ln \delta_2 - \ln\delta_1},
\label{eq:Parametrisation1}
\eeq
\beq
B = - (C_1 - C_0)\log\delta_1.
\eeq
In the case of nearly identical resources, this becomes
\beq
B = {C_2-C_1 \over \sqrt{2}( \text{inverf}\,\delta_1 - \text{inverf}\,\delta_2)},
\eeq
\beq
C_0 = \sqrt{2}B\, \text{inverf}\,\delta_1 + C_1
\label{eq:Parametrisation2}
\eeq

\subsection{Demonstrating the validity of the functional forms using IMAGE data\label{sect:SFits}}

\begin{figure}[h]
	\begin{minipage}[t]{.5\columnwidth}
		\begin{center}
			\includegraphics[width=1\columnwidth]{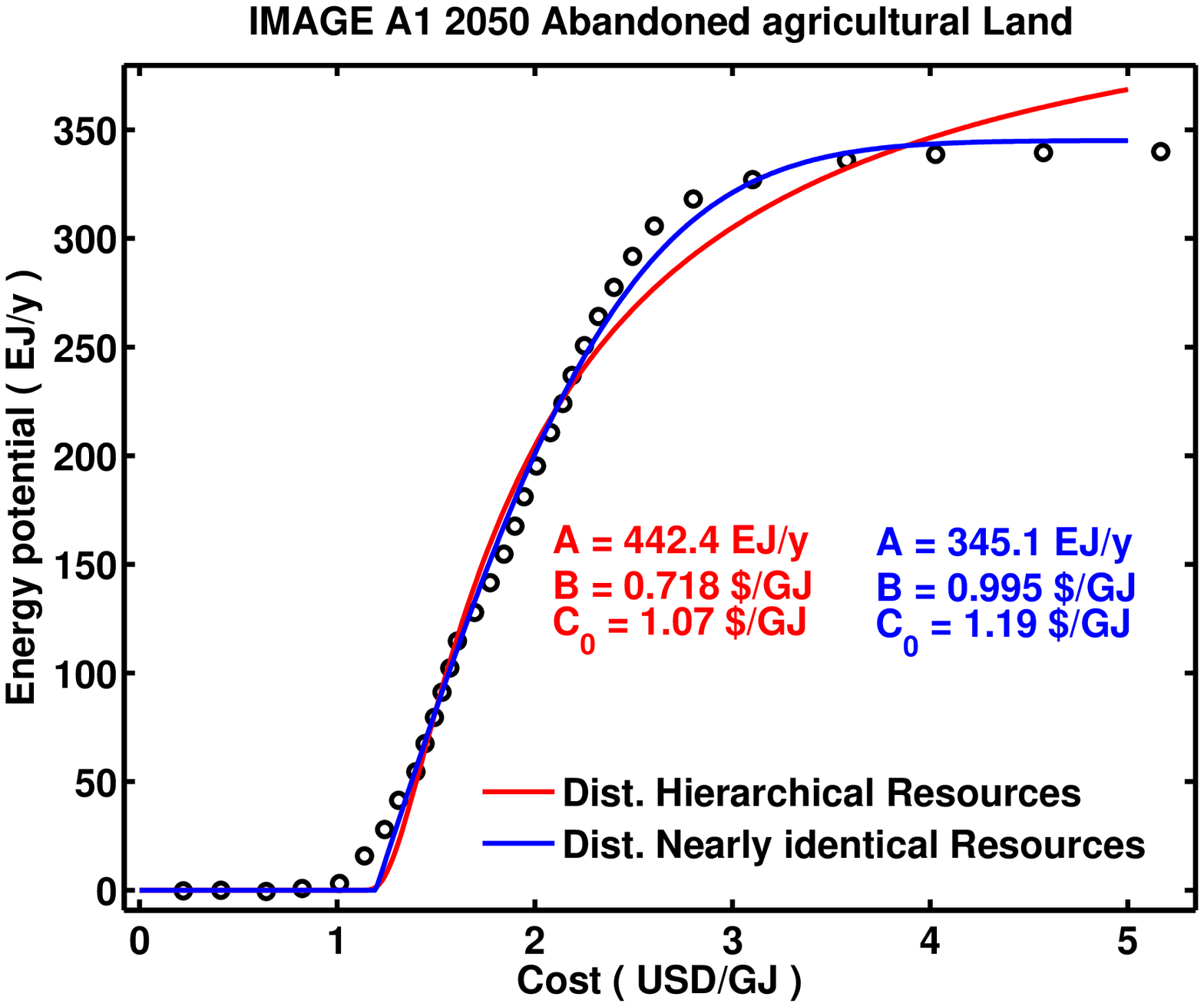}
			\includegraphics[width=1\columnwidth]{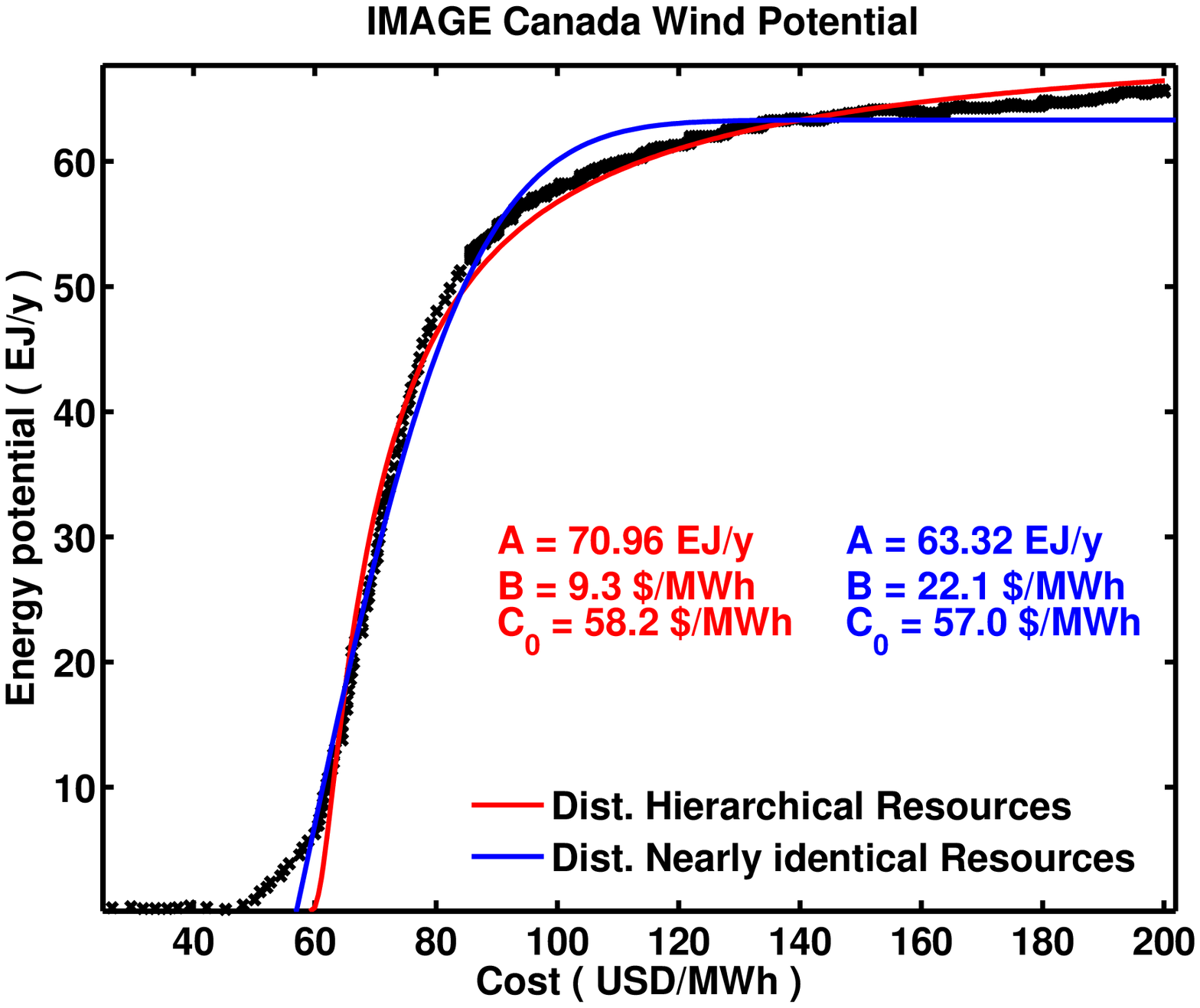}
		\end{center}
	\end{minipage}
	\hfill
	\begin{minipage}[t]{.5\columnwidth}
		\begin{center}
			\includegraphics[width=1\columnwidth]{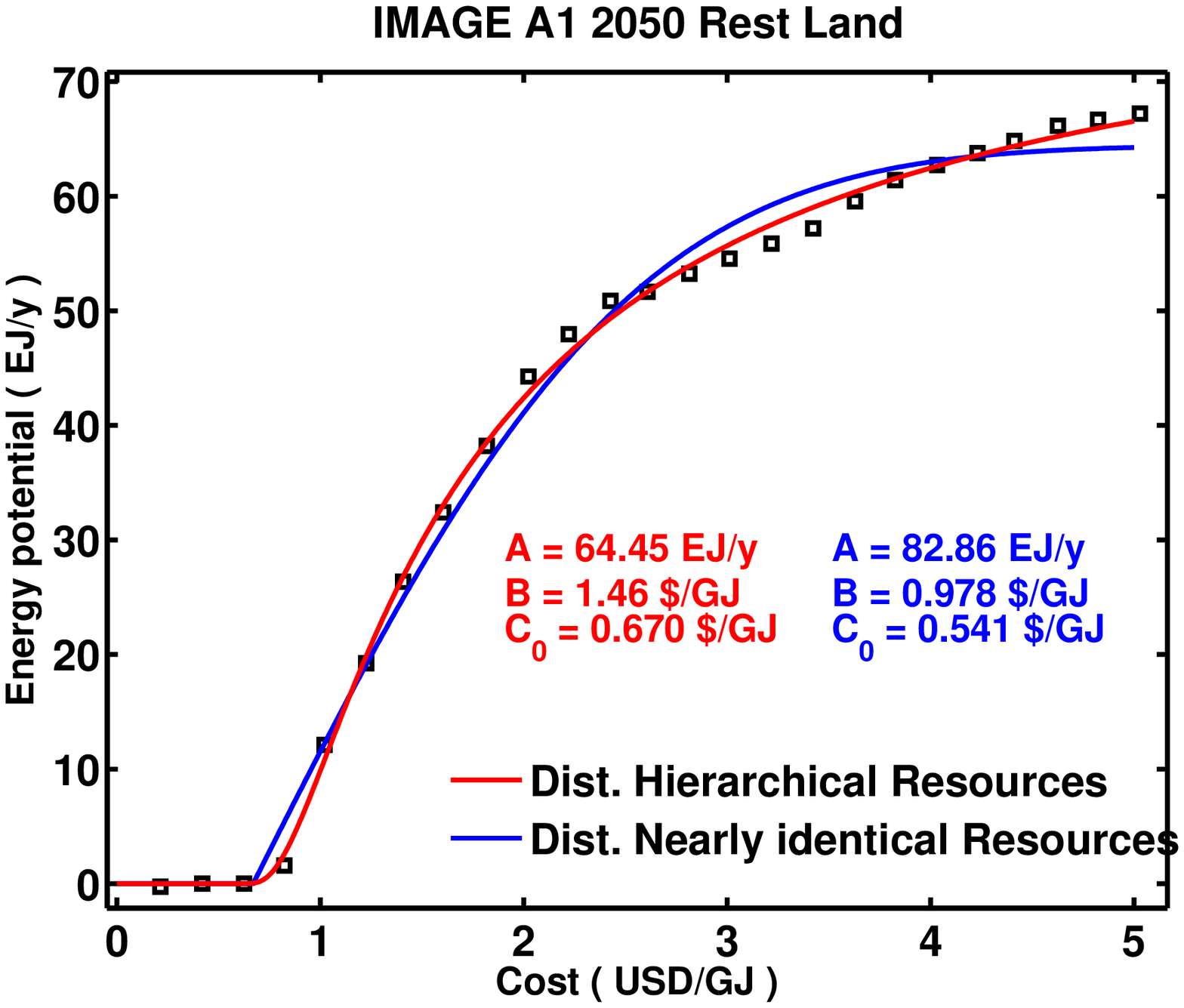}
			\includegraphics[width=1\columnwidth]{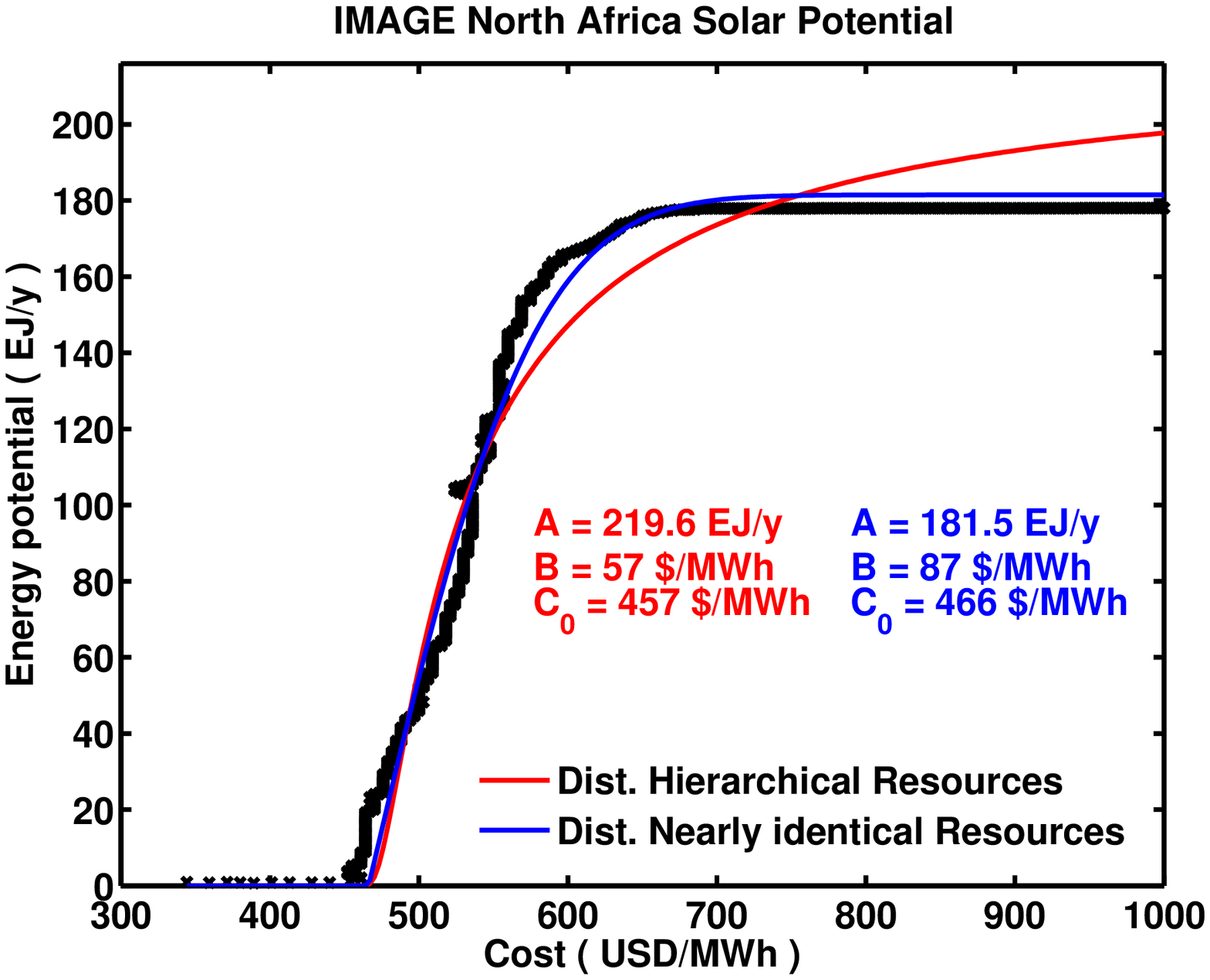}
		\end{center}
	\end{minipage}	
	\caption{Curve fits using non-linear least-squares of the two types of cumulative distribution with data from various studies of renewable energy potentials previously reported, calculated using the model IMAGE (reproduced from \cite{Hoogwijk2009, Hoogwijk2004, HoogwijkThesis}). The goodness of these fits are a good indication for which type of distribution represents best each type of resource. It can observed that data for abandoned agricultural land is well described by the cumulative distribution for nearly identical resources ($top$ $left$), while the data for rest land is described by the cumulative distribution of the hierarchical type ($top$ $right$). Data for wind energy is well described by hierarchical resource distribution ($bottom$ $left$), while the data for solar energy is well described by a distribution for nearly identical resources ($bottom$ $right$).}
	\label{fig:DistFits}
\end{figure}

Examples of the use of the analytical forms of the distributions are presented in figure \ref{fig:DistFits} for biomass, solar and wind energy. The data are taken from land use simulations performed using IMAGE by Hoogwijk $et$ $al.$, which provide the only sources of cost-supply curves calculated outside of this project that do not already use assumptions on the analytical form of the resource distribution \citep{HoogwijkThesis, Hoogwijk2009, Hoogwijk2004}. Since IMAGE simulates the use of the land on each point of a global grid, and since these cost-supply curves were calculated by building histograms of the number of grid points with productivities situated within various ranges, their form stems purely from the statistical nature of the data. These are thus appropriate for testing the functions given above.

Non-linear least-squares fits were performed with both analytical forms for each data set. In every case, only one of the two functions given above is appropriate, while the other is not. Fits are moreover of exceptional quality. For instance, wind resources are the best example of resources of the hierarchically ordered type, which stems from the exponentially increasing number of simultaneous factors required to produce ever higher productivities. The data is found to follow almost exactly the distribution for hierarchical resources (note that the deviation at low cost values stems from the aggregation of a region with a different cost structure into the region for Canada). Meanwhile, solar resources represent the best example of nearly identical resources, since in regions of similar irradiation, all sun-facing areas are equivalent. The data is found to follow closely the distribution for nearly identical resources. Biomass resources from abandoned agricultural land are nearly identical. This stems from the similar nature of local areas of agricultural land (i.e. large plains, deltas, similar irradiation, etc). Land plots with lower productivity are used for other activities. Rest land, however, is the category of land which would not be used for agriculture, and can be of various nature, but includes mainly savannah, shrubland and grassland or steppe. These can be ordered, and can be seen to follow the distribution for hierarchical resources. 

\section{Cost-supply curve calculation methodology per resource type}
\subsection{Definition of world regions \label{sect:WorldRegions}}

Cost-supply curves were calculated in this work for every E3MG world regions from aggregations of data defined for 190 countries. However, the region definition in E3MG is very specific and does not correspond closely to that of most other global models, and tables provided here for E3MG regions would be of limited use to the global modelling community. For accuracy, data for 190 countries would be required to be provided here, but is not possible for space considerations. For the convenience of potential users, the results are provided in tables with a definition of regions resembling that of other models such as IMAGE, AIM, etc. Any other aggregation of data, in table or curve form, can be supplied by the authors upon request. Table \ref{tab:Regions} gives the list of regions used here with most countries that belong to them.

\subsection{Wind and solar energy\label{sect:WindSolar}}

Following the justification of section~\ref{sect:SFits}, wind resources were modelled using a distribution of hierarchical type, while solar resources were modelled using a distribution for nearly identical resources. In both cases simulations performed by Hoogwijk $et$ $al.$ (wind energy), Hoogwijk (PhD thesis, wind, solar and biomass energy) and de Vries $et$ $al.$ (both) using IMAGE 2.2 were used, published in the form of data tables featuring both technical potentials and interpolations through cost supply curves at specific cost values for a list of 17 world regions \citep{Hoogwijk2004, HoogwijkThesis, deVries2007}. These values were used to find the distribution parameters $A$, $B$ and $C_0$ for every one of their regions. In the case of wind, $A$ values were thus obtained without additional processing. In the case of solar, $A$ values were obtained from de Vries $et$ $al.$ while $B$ and $C_0$ values from Hoogwijk \citep{deVries2007, HoogwijkThesis}. In both cases, $B$ and $C_0$ values were obtained using equations \ref{eq:Parametrisation1} to \ref{eq:Parametrisation2}. 

However, the regional aggregation in the aforementioned work does not match exactly the one chosen in the current study (or the one used in E3MG), detailed in section~\ref{sect:WorldRegions}. In order to obtain curves for this set of world regions, energy potentials from IMAGE~2.2 regions were disaggregated into 190 countries, and subsequently re-aggregated. This required additional assumptions in particular cases where $A$ values were required to be divided between underlying countries\footnote{Note that this is mostly true for E3MG regions; the regions used for this paper are very similar to those used by Hoogwijk $et$ $al.$ \citep{Hoogwijk2004}}. In the case of wind energy, the division of $A$ values was done proportionally to the cube of the yearly averaged wind speed\footnote{Wind energy scales with the cube of the average wind speed averaged over time (see for instance \citep{Sorensen2011}).}, times the amount of land suitable in each country for these energy production activities (the land area times the suitability factor provided by Hoogwijk $et$ $al.$, assumed the same for all countries member of a region) \citep{RREX, 3TIERWind2011, CIA2011}. In the case of solar energy, $A$ values were divided proportionally to the insolation averaged over countries times the amount of land suitable in each country \citep{RREX, 3TIERSolar2011, CIA2011}. Assuming an identical shape for the cost supply curves (identical values of $B$ and $C_0$) for every country within a particular IMAGE region, and using $A$ values thus divided, cost supply curves for the 190 countries were built. Given this set of curves, the re-aggregation of curves into new world regions was performed by summing the energy potential values at each cost (i.e. a sum along the horizontal axis of the cost-supply curve, called a horizontal sum henceforth). These aggregated curves do not correspond anymore to pure distributions of either type, but do not differ significantly from pure forms in any of the regions chosen for this work. Thus, new values for $A$, $B$ and $C_0$ for this work's regional definition were re-estimated using equations \ref{eq:Parametrisation1} to \ref{eq:Parametrisation2}, for the sake of simple presentation in this work (avoiding listing parameters for 190 countries, or providing aggregate curves defined on large numbers of cost data points). For E3MG, data curves evaluated on 1000 cost data points are used directly instead.\footnote{Exact analytical forms for cost-supply curves correspond to the inverse the cumulative distribution. When the cumulative distribution involves the sum of several distributions, an analytical form for the cost-supply curve does not exist.}

Cost values with which cost-supply curves were calculated using equations \ref{eq:Parametrisation1} to \ref{eq:Parametrisation2} were also obtained from Hoogwijk $et$ $al.$, but were rescaled to 2008 prices \citep{HoogwijkThesis, Hoogwijk2004}. This procedure, however, generates costs of energy production slightly different (wind) or higher (solar) than recent estimates available from the International Energy Agency, due to small errors (wind) or significant learning-by-doing cost reductions (solar) stemming from economies of scale with large expansion of electricity generation capacities that occurred between 2004 and 2008 \citep{IEAProjCosts}. The curves were therefore recalibrated with a constant offset to match recent values. The results are provided in table~\ref{tab:RenParam} for this work's list of world regions.

\subsection{Hydropower}

Hydroelectric resources, highly site dependent, were modelled using the distribution for hierarchical resources. Hydroelectric potentials and current annual electricity generation values were obtained from the last available technical report of the International Journal on Hydropower and Dams (IJHD), while the costs were obtained using an extensive study of 250 recent projects by Lako $et$ $al.$ from which statistics were derived \citep{IJHD2011, Lako2003}. These statistics were performed for the countries studied in Lako $et$ $al.$, and were used as proxies for regions not studied in their work, or where no information on recent hydroelectric developments was found. Some countries do not have recently reported hydroelectric projects onto which to base cost values. 

Recent developments have hardly followed an order of cost, since they were scattered between 500 and 4000 2003USD/kW. In order to use a cost-supply curve, it can only be assumed that future developments actually will approximately follow a cost order. Although only approximately true, this is reasonable, since development costs $will$ significantly increase when more and more usable sites are depleted, irrespective of the particular order in which they were built, and only difficult or distant river basins remain. This is important since the costs of hydroelectricity are currently not high in comparison to alternatives, but the resources are limited, and therefore the development of hydroelectric resources must be limited through an increasing cost in models of power systems such as FTT:Power. 

As can be seen in the current hydroelectricity generation data compared to the data for hydroelectric potentials from IJHD, a significant portion of the technical potential of every region is already developed \citep{IJHD2011,WEC2010}. The cost values delimiting the technical and economic potentials amongst remaining potential hydroelectric sites are not given by IJHD. Since the distribution of costs is not symmetrical, the assumption was taken that the amount of resources considered economical lies at costs between the local average cost $\mu$ minus its standard deviation $\delta$, $\mu - \delta$, and plus twice its standard deviation $\mu + 2 \delta$. This puts the upper cost limit to around 5000~2008USD/kW. Thus, sites within the technical potential with costs higher than this are considered currently uneconomical. Given this definition, a cost-supply curve for each region was calculated. Parameters for each  regional cost-supply curve are given in table \ref{tab:RenParam}. The global cost-supply curve of figure~\ref{fig:GlobalRenCSC} of the main text is an aggregation (a horizontal sum) of these regional curves.

\subsection{Geothermal energy\label{sect:SupplGeoth}}

Geothermal resources were divided into two groups, occurring in either ``in belt'' or ``out of belt'' land areas, referring to the so-called volcanic belt. ``In belt'' areas are located in volcanically active zones with high geothermal gradients (temperature gradients with bore depth from the surface of the ground). Given the particular characteristics of geothermal active areas in terms of their heat storage and underground temperature variation, the extraction of geothermal resources in those places are highly site-specific, and were therefore modelled using a hierarchical distribution. ``Out of belt'' areas corresponds to the rest of the continental masses, with sites that are characterised by smaller geothermal gradients, and that are almost identical to one another within large geographical areas. ``Out of belt'' resources were thus modelled using a distribution for nearly identical resources. The ratio of ``in belt'' to ``out of  belt'' land area values were obtained from the 1978 report of the EPRI \cite{EPRI1978}, enabling to divide reported technical potentials into two $A$ values for each distribution type. Geothermal resources were moreover calculated for both hydrothermal and EGS dry rock technologies, yielding four sets of parameters. Each cost-supply curve in each region was obtained by aggregating four curves. 

Technical potentials for different world regions were obtained from Bertani \cite{Bertani2012}. Given the differences between their regional aggregation and this work, the same methodology was used as for wind and solar energy in order to disaggregate the regional technical potentials between the same 190 countries. The proportion of the regional technical potentials assigned to every country within a particular region was assumed to be proportional to the total amount of geothermal energy stored up to five kilometres of depth in each country \citep{Aldrich1981}. 

Cost values for geothermal electricity production were taken from the IEA \cite{IEAGeo2010}. It was assumed that 90\% of the resources `in belt' were situated within these ranges of costs. However, resources `out of belt' follow the distribution for nearly identical resources, but face higher costs due to lower geothermal gradients and less productivity per unit investment. Since no additional cost information was available in this regard, these resources were assumed to lie in the upper half of the cost range given by \citep{IEAGeo2010}.\footnote{In `out of belt' areas, the same technologies are involved, either for hydrothermal or EGS, as for `in belt' areas. However, the resources are nearly identical over large areas and of equally low quality in comparison to `in belt' areas. Significantly higher productivities are found in volcanic areas.} Table \ref{tab:GeoParam} gives the parameters that can be used to reproduce these cost-supply curves using both types of distributions.

The lower boundary curve of the uncertainty range assumes a technical potential of 4~EJ/y based on differing assumptions for both technologies. In the case of hydrothermal technology, a conservative potential estimate of 70~GW (2~EJ/y) was derived by limiting the calculation to well known sites that have been already characterised by direct involvement or informed calculations \citep{Bertani2010, Bertani2012}. Meanwhile, the limited amount of accumulated experience with EGS technology creates uncertainties in the evaluation of the technical potential through variations in the efficiency of extraction \citep{Tester2006}, leading Bertani to estimate a lower limit of 70~GW (2~EJ/y) \citep{Bertani2012}. The upper boundary of the uncertainty range involves yet another set of assumptions for both hydrothermal and EGS technologies. In the case of hydrothermal, according to Stefansson, undiscovered or additional resources could exist which would be five to ten times higher than identified resources, increasing the potential to 1000-2000~GW (57~EJ/y) \citep{Stefansson2005}. In the case of EGS, the technical potential is calculated by an extrapolation of resources in the United States to the global level using the proportion between the heat stored at depths of less than 10~km in the United Stated with the known EGS primary energy potential in the same region, estimated as 1$\times10^6$~EJ of heat stored per 2.61~EJ/yr of EGS primary energy potential \citep{Tester2006}. Using the estimation of the heat stored at depths less than 10~km on the global scale of 403$\times10^6$~EJ \citep{Rowley1982}, this estimation results in a global technical potential of 57~EJ/yr, yielding a total of 114~EJ/y.

\subsection{Bioenergy}

Four bioenergy cost-supply curves are given in figure~\ref{fig:GlobalRenCSC} of the main text, for each of the SRES scenarios A1, A2, B1 and B2, based on simulations performed using IMAGE 2.2 \citep{Hoogwijk2005, Hoogwijk2009, HoogwijkThesis} (see \citep{IPCCSRES} for information on SRES scenarios). The primary biomass energy sources considered in the cost supply curve are abandoned agricultural land, rest land and bagasse, where the first is the largest source in all scenarios. Following the justification of section~\ref{sect:SFits}, abandoned agricultural land was modelled using distributions for nearly identical resources, while rest land was modelled using distributions for hierarchical resources. The remaining type of primary biomass resources, bagasse (from \citep{WEC2010}), contributes very small fractions of the total potentials, and its technical potentials were simply added to the potentials from abandoned agricultural land and rest land, for every region in every scenario. Cost values, however, are only given by Hoogwijk $et$ $al.$ for the total amount of biomass resources in each region, not individually for abandoned agricultural and rest land \citep{Hoogwijk2009}. Therefore, the right distribution to use had to be determined, by deciding which of the two represented best the data. This corresponds to finding the $dominant$ distribution. Therefore, the appropriate $type$ of distribution was determined by visual inspection for each region. Potentials for abandoned agricultural land are for most regions much larger than those for rest land, and therefore most regions were modelled using distributions for nearly identical resources. These distributions were disaggregated into 190 countries, following the methodology described in section~\ref{sect:WindSolar}, proportionally to country land areas times their suitability factor. Table \ref{tab:BiomassParam} provides values that can be used to parameterise biomass cost-supply curves for the world regions used in this work, with the appropriate type of distribution used indicated in the last column.

\subsection{Ocean energy \label{Sect:SOcean}}

Given the vast extent of oceans, the calculation of theoretical potentials for ocean energy sources produces large values. For instance, using a global wind-wave model, Mork $et$ $al.$ estimated a potential for wave energy between 2986 and 3703~GWe (94 to 117~EJ/yr), while Charlier and Justus estimated a global theoretical tidal power potential between 1000 and 3000~GWe (32 to 95~EJ/yr) using a capacity factor of 100\% \citep{Mork2010, Charlier1993}. In the case of ocean thermal energy, Pelc and Fujita estimated a theoretical potential of approximately 10~TW (315~EJ/yr) using a capacity factor of 100\%, while for salinity gradient energy, Cavanagh $et$ $al.$ calculated a value of 2.6~TW (82~EJ/yr) using a capacity factor of 100\% \citep{Pelc2002, Cavanagh1993}. Using values from these particular studies, the total theoretical potential for ocean energy would be as high 523 and 609~EJ/yr. However, more reliable and conservative potentials have also been evaluated, given below. These values are modest in comparison. As indicated in the main text, cost-supply curves were calculated for wave and tidal systems only. For presentation in this work only, these two cost-supply curves were combined into a single one for ocean energy. Parameters for regional ocean cost-supply curves are given in table \ref{tab:RenParam}.

\subsection*{Wave Energy}
In the case of wave energy, WEC estimated a maximum global installable capacity of 2~TW by limiting developments to technically favourable locations near coastlines \citep{WEC1994}. Using this value, and assuming a single capacity factor value of 32\%, Krewitt $et$ $al.$ estimated a technical potential for wave energy of 20~EJ/yr, while UNDP estimated a technical potential of 65~EJ/yr using the same value but assuming a capacity factor of 100\% instead \citep{Krewitt2009, WEA2001}.  Following a more conservative approach restricted to shorelines exceeding a power production of 30~kW/m (resulting in around  2\% of global coastlines), Sims $et$ $al.$ estimated a technical potential of 500~GW (6.3~EJ/y using a capacity factor of 40\%) \citep{Sims2007}. It is clear however that using single capacity factor values is not appropriate. It is likely that, as it is the case for wind, the cost variation of a wave energy cost-supply curve should stem from capacity factor variations which stem from the local quality of the resource. Such data is however not currently available as it is for wind resources. Capacity factor distributions are likely to follow roughly those of wind energy, which vary between 15 and 35\%, since both resources are closely related (wind capacity factor distributions were obtained by extracting the capacity factor from Hoogwijk $et$ $al.$ data \citep{Hoogwijk2004}). The assumption was therefore taken in this work that wave energy resources are captured using a single technology, with an investment cost given by ETSAP of 6600~USD/kW, with a capacity factor that varies between 35\% (where the resource quality is highest, and the cost of electricity production is lowest per unit energy produced) to a low value of 15\% (below which sites are not economically useable) \citep{IEAETSAP2010Marine}. Using a maximum global capacity of 2~TW, the cost-supply curves were calculated with a hierarchical distribution, assuming that 90\% of the wave resources are available at capacity factors within the range 15-35\%.
The disaggregation into 190 countries was performed according to the lengths of their respective coastlines, using data from the Central Intelligence Agency (CIA) \cite{CIA2011}. 

\subsection*{Tidal Energy}
Accounting for most of the global installed capacity of ocean energy systems, tidal energy is the only technology that has reached a commercial scale, with approximately 523~MW installed at the end of 2010 \citep{IEAOcean2011}. WEC made a rough estimation of the technical potential of tidal energy of about 2000~TWh/y (7.2~EJ/y), 10\% considered economical \citep{WEC1994, Rodier1992,WPCWEC1986}.  In a more detailed study, Hammons presented a global but non-exhaustive list of potential tidal sites that could be considered for development, including projected  installed capacities and approximate annual outputs \citep{Hammons1993}. The total output from these sites would be of almost 400~TWh/y (1.4 EJ/y). Hammons furthermore extrapolated that the inclusion of additional sites around the world not studied specifically in his work would result in a global technical potential for tidal energy likely to range between 500 and 1000~TWh/yr (1.8 – 3.6~EJ/y). The cost-supply curve for tidal energy was calculated using the range of cost values given in ETSAP of 5000 to 6500 USD/kW \citep{IEAETSAP2010Marine}. Existing capacity was assumed to have been built at costs below that range, while the sites reviewed by Hammons (400~TWh/y) were assumed to be associated with costs within the range \citep{Hammons1993}. Additional sites were assumed to have costs above the range. 

\subsection*{Ocean Thermal and Salinity}
The state of development of ocean thermal and salinity gradient energy technologies is currently experimental and therefore large uncertainties accompany calculations of associated energy potentials \citep{Sims2007}. Upper limits in the form of theoretical potentials have been calculated. Nihous estimated a theoretical potential for ocean thermal energy of 2.7~TW (85~EJ/yr or 23~652~TWh/yr), which corresponds to the maximum amount of energy resources that could be extracted without disrupting significantly the temperature of the upper layers of the ocean in an steady state regime, using a one-dimensional model of oceanic temperature gradients \citep{Nihous2007}. Using a similar method, Charlier and Justus produced a more conservative estimation for the theoretical potential of 1000~GWe (32~EJ/yr), assuming a capacity factor of 100\% \citep{Charlier1993}. However, according to von Arx, such a level of heat extraction would imply a decrease in the ocean surface layer temperature of approximately 1$^{\circ}$C \citep{VonArx1974}. In order to avoid such a decrease, Charlier and Justus recommend a reduced estimate based on 10~TW of usable heat replenishment rate, corresponding to 100~GWe (3.2~EJ/y) \citep{Charlier1993}. 

In the case of Salinity Gradient, based on average discharge and low flow discharge values, Skramesto estimated the theoretical potential in the range of 1600 - 1700~TWh/yr (5.8 - 6.1~EJ/yr) \citep{Skramesto2009}. Using  a global discharge rate of fresh water to seas of 44~500~km$^3$ per year, Krewitt $et$ $al.$ estimated a theoretical potential of 2000~TWh (7.2~EJ/yr), value very similar to the estimate of Skramesto \citep{Krewitt2009, Skramesto2009}.

\subsection{Oil \label{sect:SOil}}
Oil resources (Table \ref{tab:Oil}) were considered in four types of occurrences, crude oil, oil shales, extra-heavy oil and oil sands \citep{BGR2010,WEC2010}. Cost information was obtained from the IEA \cite{IEAWEO2008}. The data were aggregated into this work's world regions. For each type of occurrence for each region, a hierarchical distribution was parameterised by assuming that 1\% of the resources have extraction cost below the lower bound, while 90\% have a cost of extraction below the upper bound. The distributions were summed for each region in order to calculate regional cost-supply curves, and all distributions were summed in order to determine the global cost-supply curves given in figure \ref{fig:GlobalStockCSC} of the main text.

The curve for the lower boundary of the uncertainty range was defined by assuming that only crude oil, extra-heavy oil and oil sands reserves are available, and that the rest is either unusable or does not exist. The most probable cost-supply curve was calculated assuming that crude oil, oil sands and extra-heavy reserves and resources are available, as well as oil shale resources, but no additional amounts. The curve for the upper boundary of the uncertainty range assumes that all reserves, resources and additional amounts are available, and that an additional amount of oil shales is discovered, evaluated at 50\% of the current resources. This was done in order to compensate for the absence of speculative resources and lack of detailed information available for oil shale resources, which are likely to become larger if additional exploration is carried out, and will occur if (but only if) interest in oil shales intensifies.\footnote{For instance, if strong decarbonisation policies are implemented globally, oil shales are not likely to be explored much further, since they currently involve large processing costs and only small scale exploitation.}

\subsection{Natural gas}
Gas occurrences were considered in five forms, of which four unconventional: conventional gas, shale gas, tight gas, coalbed methane and methane hydrates \citep{BGR2010, EIA2011, WEA2001, Boyer1998, Boswell2011}. The associated cost ranges were obtained from the ETSAP \cite{IEAETSAP2010}. Of the unconventional forms, only shale gas has seen exploitation larger than experimental. An additional source of methane exists, which is thought very large, aquifer gas \citep{WEA2001}. However, its potential being very speculative, no reliable information over costs of extraction was found, and thus these were not considered in the present study. Similarly, methane hydrates provide a very large source of natural gas, however, these resources occurring under the sea, and the methods of extraction very experimental, the costs of exploitation are very large, and due to large amounts of shale gas available at lower costs, it is unclear whether the world will see wide-scale exploitation of methane hydrates. All resources were distributed into this work's world regions, except for the methane hydrates, for which it is not clear whether they are situated within territorial waters or not, and thus were assigned to an international category. Regional cost-supply curves were calculated with the same method as oil, and the global cost-supply curve is an aggregation of all regions.

The curve for the lower boundary of the uncertainty range includes conventional gas reserves only. The most probable cost-supply curve includes conventional, shale and tight gas reserves, along with half the conventional, shale and tight gas resources. It moreover includes half of the coalbed methane reserves. The curve for the upper boundary of the uncertainty range includes all reserves and resources, including methane hydrates. 

\subsection{Coal}
Although coal is a very common commodity and well known resource, information over its natural occurrences is not very detailed. Coal information was available from two sources, BGE and WEC, and table \ref{tab:Coal} was constructed using a mixture of both \citep{BGR2010, WEC2010}. Where information was inconsistent, the larger amounts were kept (such inconsistencies were not frequent nor very large). Since BGR does not report coal resources in the complete classification (proven, probable and possible reserves or resources) for all countries, some elements of the table are nil \citep{BGR2010}. This situation is likely to be due to the large amounts of coal available with conventional mining techniques, and therefore most of the resources are considered reserves, and occurrences with lower productivity or higher costs are not reported. Coal formations occur in different forms which have different calorific contents. For similar mining and transport costs, the costs of coal in terms of energy produced are higher for lower grade coals. Coal resources were divided into two categories, hard coal, including anthracite and bituminous coal which posseses higher calorific contents of between 16~500 and 35~000 kJ/t, and soft coal, including sub-bituminous coal and lignite, with calorific contents between about 11~000 and 16~500~kJ/t \citep{BGR2010}. Note that these classifications are not strictly well defined in geological terms, and that coal occurrences exist that have intermediate properties. This stems from different geophysical processes taking place during the slow formation of these hydrocarbons. 

The curves for the lower boundary of the uncertainty ranges for soft and hard coal include proven reserves only. The most probable curves include proven and probable reserves, and half of the proven and probable resources. The curves for the upper boundary of the uncertainty ranges include all proven, probable and possible amounts for both reserves and resources.

\subsection{Uranium}

Information for uranium occurrences is available from a survey of the IAEA \cite{IAEA2009}, which compiles data provided by all member countries. They are highly detailed with four classifications and four cost ranges, as seen in table \ref{tab:U}. As the amounts in each row of the table are cumulative (i.e. the data in one cell is inclusive of the sum of the cells to the left), with the associated cost values they correspond to cumulative distributions. Assuming that they should follow hierarchical distributions, the associated cumulative distribution may be fitted using a non-linear least squares method. Although a fit of a function with three parameters over four data points is hardly a reliable method to determine a best fit with any level of certainty, it nevertheless produced the best curves that could be interpolated between points, as was determined by close inspection of each fit. Note that additional data points were be defined in order to constrain the fits better, such as additional values at higher costs by repeating the last data point, assumed equal to the technical potential in the saturation region, and at (0,0). The resulting fits were found to follow the data very closely. Values for military stocks of U are uncertain, since the information that is publicly available is scarce and incomplete, and were omitted. These are very small in comparison to natural stocks \citep{IAEA2009}.

The curve for the lower boundary of the uncertainty range includes only RAR (reasonably assured reserves) in all cost ranges. The most probable curve includes RAR and inferred reserves. The curve for the upper boundary of the uncertainty range includes all four classifications of resources in all cost ranges. The speculative resources in the unassigned cost range were not included in the non-linear fitting procedure, since their distribution into existing or higher cost ranges is ambiguous. They were therefore added to the technical potential of the upper boundary of the uncertainty range. Finally, it was assumed that sea water U is too costly and uncertain to include in the present study. Given the large amount of sea water on the planet, this resources is thought very large even though the concentration of U in sea water is very low. However, due to the very slow mixing process of  sea water, it is misleading to consider the global body of sea water as a usable source of U, as it could be rapidly depleted locally, providing small amounts, without access to the remaining resources situated far offshore \citep{Bonche2005}. 

\subsection{Thorium}

Thorium (Th) deposits in the Earth's crust around the world are expected to be three times larger than those of U, as determined from isotope lifetimes and the composition of the accreted material which formed the planet \citep{Bonche2005, Sinha2006, AbuKhader2009, Suess1956}\footnote{Such arguments, which originates from astrophysics and the study of stars and nebulae, along with nuclear lifetimes of isotopes which were formed through nuclear reactions within stars or during supernova events, enable the prediction of relative amounts on Earth for all isotopes of elements of the periodic table (see for instance \cite{Suess1956}). }. However, known resources of Th are much smaller and less detailed than those of U, a situation which is the result of the relatively small interest that has been given to Th in comparison to U. Therefore, it is to be expected that Th resources increase in size significantly if at some time in the future interest grows in Th based nuclear reactors. Although the Th nuclear fuel cycle has been demonstrated several decades ago, it has not been used commercially as it involves more safety hazards related to radiation than the U fuel cycle \citep{Bonche2005, Sinha2006}. The Th nuclear fuel cycle is more efficient than that of U and therefore, involves less mass of Th per unit of electricity produced. For similar mining costs, Th resources are less expensive per unit of energy, however, the processing of Th into $^{233}$U for fuel preparation has not been performed at an industrial scale. Therefore, the cost variable in the Th cost-supply curve is highly uncertain. 

Data for Th resources were obtained from the IAEA \cite{IAEA2009}. These are provided with much less detail than for U, with four uncertainty ranges but only one cost category. Consequently, the strategy of curve-fitting cannot be applied here, and one distribution of the hierarchical type per uncertainty category per world region was parametrised in the same way as for fossil resources. The cost axis was transformed from a cost per unit of mass to a cost per unit of energy using the efficiency of the Indian experimental model reported by Sinha and Kakodkar of 2100 TJ/t \citep{Sinha2006}. The curve for the lower boundary of the uncertainty range includes RAR only. The most probable curve includes RAR and inferred reserves. The curve for the upper boundary of the uncertainty range includes all four categories.





\section*{References}

\bibliographystyle{elsarticle-num}
\bibliography{CamRefs}
\newpage
\section{Data tables and figures}

\begin{table}[h]\footnotesize
\begin{center}
		\begin{center}
			\includegraphics[width=1\columnwidth]{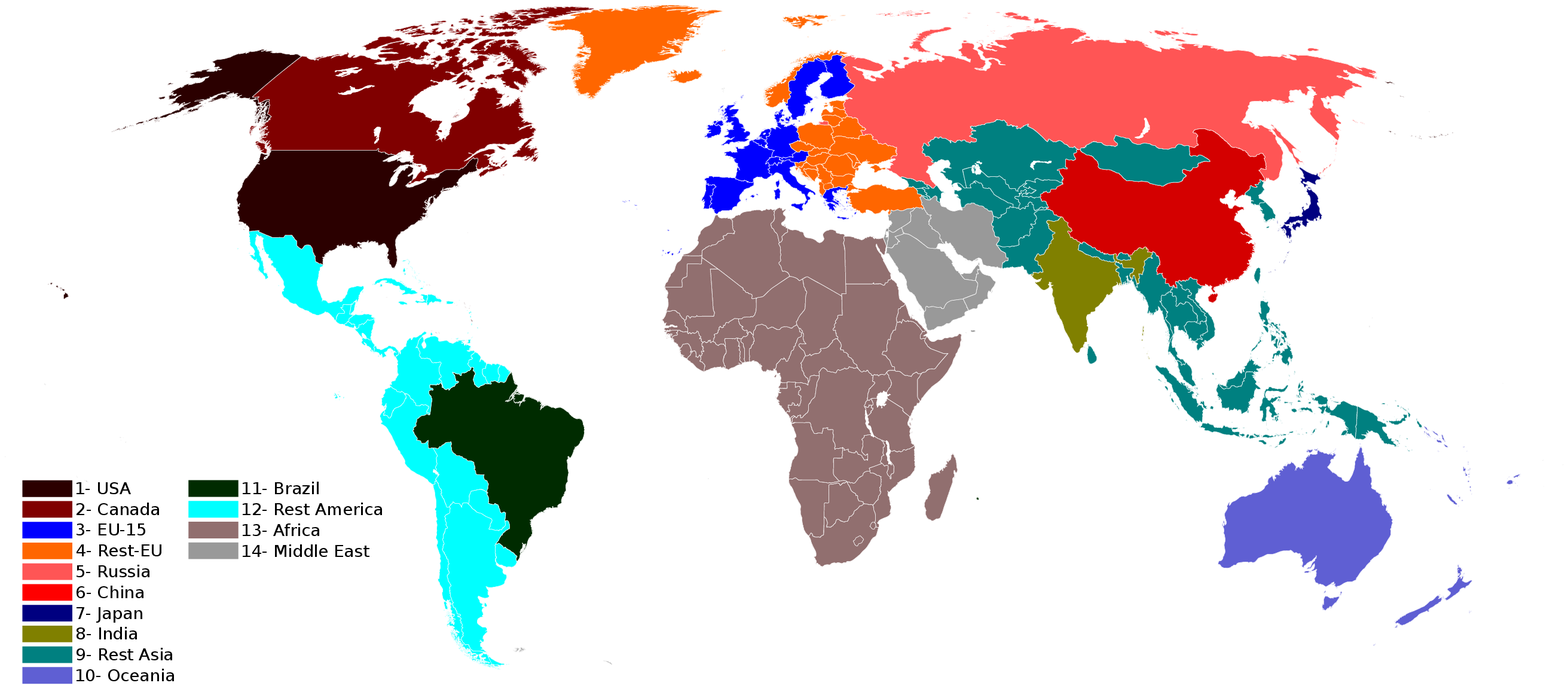}
		\end{center}
		\begin{tabular*}{1\columnwidth}{@{\extracolsep{\fill}} l|l }
			\hline
			Region & Member countries\\
			\hline
			\hline
			USA 			& USA\\
			Canada  		& Canada\\
			EU-15 		& Austria, Belgium, Danemark, Finland, France, Germany, Greece, Ireland, Italy, Luxemburg, Netherlands, Portugal, Spain, \\
						& Sweden, United Kingdom\\
			Rest Europe 	& Albania, Belarus, Bosnia-Herzegovina, Bulgaria, Croatia, Cyprus, Czech Rep., Estonia, Hungary, Iceland, Latvia, Lithuania, \\								& Macedonia, Malta, Moldova, Montenegro, Norway, Poland, Romania, Serbia, Slovakia, Slovenia, Switzerland, Turkey, Ukraine\\
			Russia 		& Russia\\
			China 		& China\\
			Japan	 	& Japan\\
			India 		& India\\
			Rest Asia 		& Afghanistan, Armenia, Azerbaijan, Bangladesh, Bhutan, Brunei, Cambodia, Georgia, Hong Kong, Indonesia, Kazakhtan, Korea,\\ 							& Kyrghizstan, Laos, Malaysia, Mongolia, Myanmar, Nepal, Pakistan, Philippines, Singapore, Taiwan, Tajikistan, Thailand, \\
						& Turkmenistan, Uzbeksitan, Viet Nam\\
			Oceania 		& Australia, New Zealand, Papua New Guinea, pacific islands\\
			Brazil	 	& Brazil\\
			Rest America 	& Mexico, Central America, South America excluding Brazil\\
			Africa	 	& Africa\\
			Middle East	& Barhain, Iran, Iraq, Israel, Jordan, Kuwait, Lebanon, Oman, Palestine, Qatar, Saudi Arabia, Syria, UAE, Yemen\\
			\hline
		\end{tabular*}
	\caption{Definition of world regions for this paper with member countries for each.}
	\label{tab:Regions}
\end{center}
\end{table}

\begin{table}[p]\footnotesize
\begin{center}
		\begin{tabular*}{1\columnwidth}{@{\extracolsep{\fill}} l|r r r|r r r|r r r }
			\hline
			\multicolumn{10}{ l }{Renewable cost-supply curve parameters}\\
			\hline
			 & \multicolumn{3}{ l|}{Hydro}  & \multicolumn{3}{ l|}{Wind} & \multicolumn{3}{ l }{Solar}\\
			\hline
			  & \multicolumn{3}{ l|}{Hierarchical}  & \multicolumn{3}{ l|}{Hierarchical} & \multicolumn{3}{ l }{Nearly identical}\\
			Region & $A$ & $B$ & $C_0$ & $A$ & $B$ & $C_0$ & $A$ & $B$ & $C_0$ \\
			Name& PJ/y & \$/MWh & \$/MWh & PJ/y & \$/MWh & \$/MWh & PJ/y & \$/MWh & \$/MWh\\
			\hline
			\hline
			USA 			&5 746	&198.40	&1.40	&75 600	&30.19	&145.53	&262 800	&350.03	&620.44	\\
			Canada  		&4 217	&38.82	&78.06	&43 505	&10.17	&149.28	&51 571	&1542.58	&839.26	\\
			EU-15 		&3 066	&22.01	&88.25	&12 009	&65.77	&113.60	&61 471	&817.28	&616.21	\\
			R. Eur. 		&4 283 	&54.13	&38.85	&32 183	&19.38	&142.14	&57 355	&593.61	&718.00	\\
			Russia 		&6 595	&82.84	&75.89	&38 774	&65.10	&129.73	&384 448	&428.97	&611.15	\\
			China 		&11 502	&26.37	&35.65	&6 057	&175.01	&161.04	&175 509	&241.71	&695.79	\\
			Japan	 	&846	&25.24	&82.31	&360	&109.24	&161.04	&3 600	&185.69	&840.10	\\
			India 		&2 788	&103.23	&0.00	&2 018	&109.62	&208.48	&120 525	&140.01	&483.33	\\
			R. Asia 		&7 227	&160.05	&0.00	&18 125	&83.64	&123.77	&282 441	&369.67	&498.17	\\
			Oceania 		&779	&71.64	&0.00	&50 410	&59.23	&138.60	&430 421	&141.75	&521.29	\\
			Brazil	 	&4 738	&18.45	&12.00	&13 248	&23.33	&127.38	&113 386	&310.45	&527.52	\\
			R. Amer.	 	&7 414	&127.92	&0.00	&22 752	&30.14	&132.96	&160 214	&268.27	&527.44	\\
			Africa	 	&5 767	&69.16	&64.42	&23 106	&143.77	&161.04	&974 259	&157.16	&456.31	\\
			Mid. East		&1 094	&304.10	&0.00	&7 200	&109.62	&208.48	&306 000	&172.80	&437.87	\\
			\hline
			Total			&\multicolumn{3}{ l|}{66 061}	&\multicolumn{3}{ l|}{345 348}	&\multicolumn{3}{ l }{3 384 000}\\
			\hline
		\end{tabular*}
	\caption{Table of cost-supply curve parameters for each region for hydro, wind and solar power.  }
	\label{tab:RenParam}
\end{center}
\begin{center}
		\begin{tabular*}{1\columnwidth}{@{\extracolsep{\fill}} l|r r r|r r r }
			\hline
			\multicolumn{7}{ l }{Renewable cost-supply curve parameters}\\
			\hline
			 & \multicolumn{3}{ l|}{Wave}  & \multicolumn{3}{ l }{Tidal}\\
			\hline
			  & \multicolumn{3}{ l|}{Hierarchical}  & \multicolumn{3}{ l }{Hierarchical}\\
			Region & $A$ & $B$ & $C_0$ & $A$ & $B$ & $C_0$\\
			Name& PJ/y & \$/MWh & \$/MWh & PJ/y & \$/MWh & \$/MWh\\
			\hline
			\hline
			USA		&496	&32.46	&199.44	&145	&89.18	&303.33\\
			Canada	&5030	&32.46	&199.44	&757	&89.18	&303.33\\
			EU-15	&1525	&32.46	&199.44	&287	&89.18	&303.33\\
			R. Eur.	&2442	&32.46	&199.44	&333	&89.18	&303.33\\
			Russia	&937	&32.46	&199.44	&743	&89.18	&303.33\\
			China	&361	&32.46	&199.44	&49		&89.18	&303.33\\
			Japan	&741	&32.46	&199.44	&101	&89.18	&303.33\\
			India		&174	&32.46	&199.44	&89		&89.18	&303.33\\
			R. Asia	&2843	&32.46	&199.44	&396	&89.18	&303.33\\
			Oceania	&1536	&32.46	&199.44	&238	&89.18	&303.33\\
			Brazil	&186	&32.46	&199.44	&25		&89.18	&303.33\\
			R. Amer.	&1303	&32.46	&199.44	&253	&89.18	&303.33\\
			Africa	&1027	&32.46	&199.44	&140	&89.18	&303.33\\
			Mid. East	&309	&32.46	&199.44	&42		&89.18	&303.33\\
			\hline
			Total		&18 910	&		&		&3600	&		&\\
			\hline
		\end{tabular*}
	\caption{Table of cost-supply curve parameters for each region for wave and tidal energy.  }
	\label{tab:RenParam}
\end{center}
\end{table}
 
\begin{table}[p]\footnotesize
\begin{center}
		\begin{tabular*}{1\columnwidth}{@{\extracolsep{\fill}} l|r r r r|r r r r}
			\hline
			\multicolumn{9}{ l }{Primary biomass cost-supply curve parameters}\\
			\hline
			  & \multicolumn{4}{ l|}{A1}  & \multicolumn{4}{ l }{A2} \\
			\hline
			Region & $A$ & $B$ & $C_0$ &Type & $A$ & $B$ & $C_0$ & Type\\
			Name & PJ/y & \$/MWh & \$/MWh && PJ/y & \$/MWh & \$/MWh &\\
			\hline
			\hline
			USA 			&53 082	&3.93	&3.60	&2	&33 082	&6.42	&3.10	&2	\\
			Canada  		&18 000	&1.64	&3.60	&2	&12 000	&1.50	&3.60	&2	\\
			EU-15 		&11 305	&1.62	&6.15	&2	&11 305	&0.93	&6.18	&2	\\
			Rest Europe 	&11 855 	&0.43	&6.28	&2	&10 781	&1.45	&3.44	&2	\\
			Russia 		&126 017	&1.73	&3.60	&2	&67474	&1.96	&3.16	&2	\\
			China 		&107 100	&11.26	&5.04	&1	&23 322	&9.92	&7.20	&2	\\
			Japan	 	&115	&0.69	&7.20	&2	&44		&9.98	&14.40	&2	\\
			India 		&25 890	&3.37	&3.00	&2	&13 756	&2.32	&2.86	&2	\\
			Rest Asia 		&12 921	&0.91	&3.48	&2	&8 605	&0.35	&3.52	&2	\\
			Oceania 		&55 625	&2.01	&3.14	&2	&34 477	&3.31	&2.52	&2	\\
			Brazil	 	&77 751	&1.33	&6.54	&2	&22 070	&5.12	&3.81	&2	\\
			Rest America 	&27 629	&2.06	&5.45	&2	&7 310	&3.58	&3.44	&2	\\
			Africa	 	&139 246	&4.32	&1.68	&2	&53 240	&6.78	&-0.86	&2	\\
			Middle East	&13 011	&25.38	&6.96	&1	&8 011	&13.08	&7.20	&2	\\
			\hline
			Total			&679 548	&	&	&	&305 477	&	&	&	\\
			\hline
			\hline
			  & \multicolumn{4}{ l|}{B1}  & \multicolumn{4}{ l }{B2} \\
			\hline
			Region & $A$ & $B$ & $C_0$ &Type & $A$ & $B$ & $C_0$ & Type\\
			Name & PJ/y & \$/MWh & \$/MWh && PJ/y & \$/MWh & \$/MWh &\\
			\hline
			\hline
			USA 			&36 082	&0.85	&5.00	&2	&49 082	&2.07	&3.60	&2	\\
			Canada  		&14 000	&0.84	&3.60	&2	&13 000	&0.94	&3.60	&2	\\
			EU-15 		&7 268	&3.29	&5.93	&1	&12 921	&4.77	&4.48	&1	\\
			Rest Europe 	&9 844  	&2.42	&3.59	&1	&12 178	&0.45	&5.70	&2	\\
			Russia 		&87 319	&0.99	&3.60	&2	&77 396	&0.94	&3.60	&2	\\
			China 		&77 179	&1.57	&3.60	&2	&46 261	&5.51	&7.20	&2	\\
			Japan	 	&115	&0.69	&7.20	&2	&225	&0.69	&7.20	&2	\\
			India 		&13 756	&3.35	&2.92	&2	&6 289	&5.24	&3.60	&2	\\
			Rest Asia 		&5 075	&0.72	&3.49	&2	&5 349	&1.80	&3.85	&2	\\
			Oceania 		&35 279	&0.96	&2.81	&2	&30 329	&0.87	&3.09	&2	\\
			Brazil	 	&56 539	&4.76	&4.44	&1	&38 863	&2.29	&6.03	&2	\\
			Rest America 	&18 841	&3.86	&3.00	&2	&10 517	&12.20	&3.60	&1	\\
			Africa	 	&81 240	&3.36	&2.50	&2	&15 237	&4.47	&1.17	&2	\\
			Middle East	&4 011	&4.99	&7.20	&2	&3 011	&5.49	&7.20	&2	\\
			\hline
			Total			&446 548	&	&	&	&320 658	&	&	&	\\
			\hline
			
		\end{tabular*}
	\caption{Table of cost-supply curve parameters for biomass primary energy resources for four SRES scenarios A1, A2, B1 and B2.}
	\label{tab:BiomassParam}
\end{center}
\begin{center}
		\begin{tabular*}{1\columnwidth}{@{\extracolsep{\fill}} l|r r r|r r r|r r r|r r r }
			\hline
			\multicolumn{10}{ l }{Geothermal energy cost-supply curve parameters}\\
			\hline
			  & \multicolumn{6}{ l|}{Direct Use}  & \multicolumn{6}{ l }{Electricity}\\
			\hline
			  & \multicolumn{3}{ l|}{In belt}  &\multicolumn{3}{ l|}{Out of belt}  &\multicolumn{3}{ l|}{In belt}  & \multicolumn{3}{ l }{Out of belt}\\
			  & \multicolumn{3}{ l|}{Hierarchical}  & \multicolumn{3}{ l|}{Nearly identical} & \multicolumn{3}{ l|}{Hierarchical} & \multicolumn{3}{ l }{Nearly identical}\\
			Region & $A$ & $B$ & $C_0$ & $A$ & $B$ & $C_0$ & $A$ & $B$ & $C_0$ & $A$ & $B$ & $C_0$ \\
			Name& PJ/y & \$/MWh & \$/MWh & PJ/y & \$/MWh & \$/MWh & PJ/y & \$/MWh & \$/MWh & PJ/y & \$/MWh & \$/MWh \\
			\hline
			\hline
			USA		&74		&7.22	&94.34	&138	&81.44	&118.44	&1 361	&20.62	&144.96	&2 528	&63.62	&255.71\\
			Canada	&9		&7.22	&94.34	&83		&81.44	&118.44	&253	&20.62	&144.96	&2 280	&63.62	&255.71\\
			EU-15	&1		&7.22	&94.34	&89		&81.44	&118.44	&14		&20.62	&144.96	&1 458	&63.62	&255.71\\
			R. Eur.	&12		&7.22	&94.34	&27		&81.44	&118.44	&123	&20.62	&144.96	&524	&63.62	&255.71\\
			Russia	&9		&7.22	&94.34	&177	&81.44	&118.44	&275	&20.62	&144.96	&5 226	&63.62	&255.71\\
			China	&44		&7.22	&94.34	&103	&81.44	&118.44	&809	&20.62	&144.96	&1 887	&63.62	&255.71\\
			Japan	&15		&7.22	&94.34	&0		&81.44	&118.44	&144	&20.62	&144.96	&0		&63.62	&255.71\\
			India		&2		&7.22	&94.34	&37		&81.44	&118.44	&39		&20.62	&144.96	&743	&63.62	&255.71\\
			R. Asia	&115	&7.22	&94.34	&90		&81.44	&118.44	&1 498	&20.62	&144.96	&1 482	&63.62	&255.71\\
			Oceania	&20		&7.22	&94.34	&80		&81.44	&118.44	&222	&20.62	&144.96	&1 647	&63.62	&255.71\\
			Brazil	&5		&7.22	&94.34	&96		&81.44	&118.44	&102	&20.62	&144.96	&1 942	&63.62	&255.71\\
			R. Amer.	&143	&7.22	&94.34	&91		&81.44	&118.44	&2 030	&20.62	&144.96	&1 552	&63.62	&255.71\\
			Africa	&54		&7.22	&94.34	&210	&81.44	&118.44	&761	&20.62	&144.96	&4 179	&63.62	&255.71\\
			Mid East.	&11		&7.22	&94.34	&49		&81.44	&118.44	&188	&20.62	&144.96	&901	&63.62	&255.71\\
			\hline
			Total			&\multicolumn{3}{ l|}{514}	&\multicolumn{3}{ l|}{1 269} &\multicolumn{3}{ l|}{7 818}&\multicolumn{3}{ l }{26 349}\\
			\hline
		\end{tabular*}
	\caption{Table of cost-supply parameters for geothermal energy, for both direct use of heat and electricity production.}
	\label{tab:GeoParam}
\end{center}
\end{table}
\newpage

\begin{table}[p]\footnotesize
\begin{center}
		\begin{tabular*}{1\columnwidth}{@{\extracolsep{\fill}} l|r r|r|r r r|r r r }
			\hline
			\multicolumn{1}{ l }{Oil} &\multicolumn{1}{ l }{Mtoe}&\multicolumn{8}{ l }{}\\
			\hline
			  & \multicolumn{2}{ l|}{Crude Oil} & Oil Shales & \multicolumn{3}{ l|}{Oil Sands} &\multicolumn{3}{ l }{Extra Heavy Oil}\\
			Region  & \multicolumn{2}{ l|}{\cite{BGR2010}} & \cite{WEC2010} & \multicolumn{3}{ l|}{\cite{WEC2010}} &\multicolumn{3}{ l }{\cite{WEC2010}}\\
			Name & Reserves & Resources & Resources & Reserves & Resources & Additional & Reserves & Resources & Additional\\
			\hline
			\hline
			USA 			&3 863	&10 000	&536 931	&0		&5 429	&2 388	&3		&379	&4	\\
			Canada  		&667	&2 400	&2 192	&24 909	&227 189	&355 828	&0		&0		&0	\\
			EU-15 		&1 193	&1 545	&13 248	&31		&276	&0		&24		&1 928	&0	\\
			Rest Europe 	&1 155	&3 530	&4 411	&0		&1		&0		&5		&51		&0	\\
			Russia 		&10 436	&16 400	&35 470	&4 147	&39 034	&7 505	&1		&25		&0	\\
			China 		&2 018	&2 300	&47 600	&0		&233	&0		&110	&1 168	&0	\\
			Japan	 	&6		&10		&0		&0		&0		&0		&0		&0		&0	\\
			India 		&792	&400	&0		&0		&0		&0		&0		&0		&0	\\
			Rest Asia 		&9 218	&10 535	&3 988	&6 203	&55 945	&0		&18		&1 163	&0	\\
			Oceania 		&595	&1 100	&4 534	&0		&0		&0		&0		&0		&0	\\
			Brazil	 	&2 450	&5 000	&11 734	&0		&0		&0		&0		&0		&0	\\
			Rest America 	&8 996	&9 588	&60		&0		&136	&0		&8 476	&270 637	&27 704	\\
			Africa	 	&17 277	&15 485	&23 317	&263	&2 364	&6 778	&7		&66		&0	\\
			Middle East	&102 366	&21 170	&5 792	&0		&0		&0		&0		&0		&0	\\
			\hline
			Total &  161 031 & 99 463 & 689 277 & 35 552 & 330 607 & 372 500 & 8 644 & 275 416 & 27 707\\
			\hline
			\multicolumn{10}{ l }{} \\
			\hline
			\multicolumn{1}{ l }{Costs} &\multicolumn{1}{ l }{USD2008 / boe}&\multicolumn{1}{ l }{\cite{IEAWEO2008}}&\multicolumn{7}{ l }{}\\
			\hline
			  & \multicolumn{2}{ l|}{Crude Oil} & Oil Shales & \multicolumn{3}{ l|}{Oil Sands} &\multicolumn{3}{ l }{Extra Heavy Oil}\\
			 & Reserves & Resources & Resources & Reserves & Resources & Additional & Reserves & Resources & Additional\\
			 \hline
			Upper & 10 & 10 & 50 & 40 & 40 & 40 & 10 & 40 & 10\\
			Lower & 40 & 100 & 100 & 50 & 70 & 70 & 50 & 70 & 70\\
			\hline
		\end{tabular*}
	\caption{Oil resources by world region in units of Mtoe (million tonnes of oil).}
	\label{tab:Oil}
\end{center}

\begin{center}
		\begin{tabular*}{1\columnwidth}{@{\extracolsep{\fill}} l|r r|r r|r r|r r|r r }
			\hline
			\multicolumn{1}{ l }{Gas} &\multicolumn{1}{ l }{10$^9$m$^3$}&\multicolumn{9}{ l }{}\\
			\hline
			  & \multicolumn{2}{ l|}{Conv. gas} & \multicolumn{2}{ l|}{Shale gas} & \multicolumn{2}{ l|}{Tight gas} &\multicolumn{2}{ l|}{Coalbed Methane}&\multicolumn{2}{ l }{Methane Hydrates}\\
			Region & \multicolumn{2}{ l|}{\cite{BGR2010}} & \multicolumn{2}{ l|}{\cite{EIA2011}} & \multicolumn{2}{ l|}{\cite{BGR2010}} &\multicolumn{2}{ l|}{\cite{Boyer1998}}&\multicolumn{2}{ l }{\cite{Boswell2011}}\\
			Name & Reserves & Resources & Reserves & Resources & Reserves & Resources & Reserves & Resources & Reserves & Resources \\
			\hline
			\hline
			USA 			&7 080	&20 000	&17 000	&45 600	&1000	&210 000	&9 700	&2 000	&0	&0	\\
			Canada  		&1 754	&7 000	&10 988	&42 198	&0		&7 000	&5 700	&70 800	&0	&0	\\
			EU-15 		&2 338	&2 530	&7 024	&28 547	&0		&7 000	&2 802	&0		&0	&0	\\
			R. Eur. 		&3 889  	&7 510	&9 374	&39 734	&0		&0		&1 908	&0		&0	&0	\\
			Russia 		&47 578	&105 000	&538	&2 152	&0		&45 000	&17 000	&96 300	&0	&0	\\
			China 		&2 455	&10 000	&36 109	&144 463	&0		&9 000	&30 000	&5 100	&0	&0	\\
			Japan	 	&21		&5		&0		&0		&0		&0		&0		&0		&0	&0	\\
			India 		&1 115	&900	&1 784	&8 213	&0		&1 000	&800	&0		&0	&0	\\
			Rest Asia 		&23 951	&22 805	&1 444	&5 834	&0		&0		&1 100	&0		&0	&0	\\
			Oceania 		&3 553	&2 450	&11 215	&39 111	&0		&1 000	&8 500	&5 700	&0	&0	\\
			Brazil	 	&365	&2 000	&7 533	&25 658	&0		&6 000	&0		&0		&0	&0	\\
			R. Amer.	 	&7 704	&8 858	&47 579	&170 745	&0		&0		&0		&0		&0	&0	\\
			Africa	 	&14 753	&16 155	&29 482	&112 206	&0		&0		&800	&0		&0	&0	\\
			Mid. East		&75 358	&35 370	&0		&0		&0		&0		&0		&0		&0	&0	\\
			International	&0		&0		&0		&0		&0		&38 000	&0		&0		&300 000	&300 000	\\
			\hline
			Total &  191 914 & 240 583 & 180 070 & 1 0748 62 & 1 000 & 324 000 & 78 310 & 179 900 & 300 000 & 300 000\\		
			\hline
			\multicolumn{11}{ l }{} \\
			\hline
			\multicolumn{1}{ l }{Costs} &\multicolumn{2}{ l }{USD2008 / GJ}&\multicolumn{2}{ l }{\cite{IEAETSAP2010}}&\multicolumn{6}{ l }{}\\
			\hline
			  & \multicolumn{2}{ l|}{Conv. gas} & \multicolumn{2}{ l|}{Shale gas} & \multicolumn{2}{ l|}{Tight gas} &\multicolumn{2}{ l|}{Coalbed Methane}&\multicolumn{2}{ l }{Methane Hydrates}\\
			 & Reserves & Resources & Reserves & Resources & Reserves & Resources & Reserves & Resources & Reserves & Resources \\
			\hline
			Upper & 0.5 & 0.5 & 3.8 & 3.8 & 2.6 & 2.6 & 3.8 & 3.8 & 4.4 & 4.4\\
			Lower & 5.7 & 5.7 & 8.6 & 8.6 & 7.6 & 7.6 & 7.6 & 7.6 & 8.6 & 8.6\\
			\hline

		\end{tabular*}
	\caption{Natural gas resources by world region in units of Gm$^3$ (billion cube meters).}
	\label{tab:Gas}
\end{center}
\end{table}

\begin{table}[p]\footnotesize
\begin{center}
		\begin{tabular*}{1\columnwidth}{@{\extracolsep{\fill}} l|r|r r|r|r r }
			\hline
			\multicolumn{1}{ l }{Coal} &\multicolumn{1}{ l }{Mt}&\multicolumn{5}{ l }{}\\
			\hline
			  & \multicolumn{6}{ l }{Hard coal}\\
			\hline
			  & \multicolumn{3}{ l|}{Reserves}  & \multicolumn{3}{ l }{Resources}\\
			\hline
			Region & Proven & Probable & Possible & Proven & Probable & Possible\\
			Name  & \cite{WEC2010} & \cite{BGR2010} & & \cite{WEC2010} & \cite{BGR2010} \\
			\hline
			\hline
			USA 			&226 694	&0		&0		&6 691 942	&0		&0		\\
			Canada  		&4 346	&0		&0		&187 606		&0		&0		\\
			EU-15 		&84 721	&0		&0		&278 420		&0		&0		\\
			Rest Europe 	&24 534	&752	&1 862	&254 658		&7 428	&11 422	\\
			Russia 		&68 655	&0		&0		&2 730 810	&0		&0		\\
			China 		&180 600	&0		&0		&681 600		&0		&0		\\
			Japan	 	&340	&0		&0		&4 603		&1 988	&7 375	\\
			India 		&56 100	&0		&0		&105 820		&123 470	&37 920	\\
			Rest Asia 		&54 678	&0		&0		&289 048		&0		&0		\\
			Oceania 		&44 627	&0		&0		&1 620 675	&0		&0		\\
			Brazil	 	&1 547	&0		&0		&6 212		&0		&0		\\
			Rest America 	&9 960	&4572	&4 237	&20 496		&0		&0		\\
			Africa	 	&32 546	&0		&0		&58 150		&0		&0		\\
			Middle East	&1 203	&0		&0		&41 203		&0		&0		\\
			\hline
			Total & 790 551 & 5 324 & 6 099 & 12 971 243 & 132 886 & 56 717 \\
			\hline
			  & \multicolumn{6}{ l }{Soft coal}\\
			\hline
			  & \multicolumn{3}{ l|}{Reserves}  & \multicolumn{3}{ l }{Resources}\\
			\hline
			Region & Proven & Probable & Possible & Proven & Probable & Possible\\
			Name  & \cite{WEC2010} & \cite{BGR2010} & & \cite{WEC2010} & \cite{BGR2010} \\
			\hline
			\hline
			USA 			&30 851	&0		&0		&1 398 669	&0		&0	\\
			Canada  		&3 108	&0		&0		&17 371		&40 055	&108 995	\\
			EU-15 		&44 214	&0		&0		&89 158		&0		&0		\\
			Rest Europe 	&40 456 	&1 996	&3 124	&275 185		&14 961	&11 581	\\
			Russia 		&91 350	&0		&0		&1 371 030	&0		&0		\\
			China 		&52 300	&0		&0		&318 000		&0		&0	\\
			Japan	 	&10		&0		&0		&160		&1 132	&4 074	\\
			India 		&4 895	&0		&0		&38 647		&0		&0		\\
			Rest Asia 		&30 762	&7 086	&34 070	&387 263		&11 871	&57 198	\\
			Oceania 		&37 738	&62 840	&101 100	&46 973		&73 102	&112 300	\\
			Brazil	 	&4 559	&7 559	&4 575	&6 513		&10 799	&6 535	\\
			Rest America 	&5 633	&527	&790	&7 524		&0		&0	\\
			Africa	 	&180	&0		&0		&338		&0		&0	\\
			Middle East	&0		&0		&0		&0			&0		&0	\\
			\hline
			Total 		&346 056 &80 008 &143 659	&3 956 831 	&151 920 & 300 683\\
			\hline
			\multicolumn{7}{ l }{} \\
			\hline
			\multicolumn{1}{ l }{Costs} &\multicolumn{1}{ l }{USD2008 / t}&\multicolumn{1}{ l }{\cite{IEAWEO2008}}&\multicolumn{4}{ l }{}\\
			\hline
			  & \multicolumn{3}{ l|}{Reserves}  & \multicolumn{3}{ l }{Resources}\\
			\hline
			 & Proven & Probable & Possible & Proven & Probable & Possible\\
			 \hline
			 & 20 & 20 & 20 & 20 & 20 & 20\\
			 & 50 & 50 & 50 & 100 & 100 & 100\\
			\hline
		\end{tabular*}
	\caption{Coal resources by world region in units of Mt (million tonnes of coal). Hard coal includes anthracite and bituminous coal, while soft coal includes sub-bituminous coal and lignite. Since there is no clear demarcation between ranks of coal, the limit is put onto the calorific content, and thus coal resources with a calorific content higher than 16~500~kj/t belong to the hard coal category (as defined in \cite{BGR2010}), while coal resources with a lower calorific content belong to soft coal. Anthracite can have calorific contents of up to 35~000~kJ/t while lignite can have calorific values as low as 11~000~kJ/t.}
	\label{tab:Coal}
\end{center}
\end{table}

\begin{table}[p]\footnotesize
\begin{center}
		\begin{tabular*}{1\columnwidth}{@{\extracolsep{\fill}} l|r r r r|r r r r }
			\hline
			\multicolumn{1}{ l }{Uranium} &\multicolumn{1}{ l }{t}&\multicolumn{1}{ l }{\cite{IAEA2009}}&\multicolumn{6}{ l }{}\\
			\hline
			Region & \multicolumn{4}{ l|}{Reasonably Assured Reserves (RAR)}  & \multicolumn{4}{ l }{Inferred}\\
			\hline
			Name & $<$40\$/kg & $<$80\$/kg  & $<$130\$/kg  & $<$260\$/kg  & $<$40\$/kg  & $<$80\$/kg & $<$130\$/kg & $<$260\$/kg \\
			\hline
			\hline
			USA 			&0		&39 000		&207 400		&472 100		&0		&19 500	&103 700	&236 050	\\
			Canada  		&267 100	&336 800		&361 100		&387 400		&99 700	&110 600	&124 200	&157 200	\\
			EU-15 		&0		&7 000		&20 800		&33 800		&0		&0		&13 500	&110 400	\\
			Rest Europe 	&2 500   	&39 100		&88 500		&160 000		&3 200	&15 000	&43 850	&109 850	\\
			Russia 		&0		&100 400		&181 400		&181 400		&0		&57 700	&298 900	&384 900	\\
			China 		&52 000	&100 900		&115 900		&115 900		&15 400	&49 100	&55 500	&55 500	\\
			Japan	 	&0		&0			&6 600		&6 600		&0		&0		&3 300	&3 300	\\
			India 		&0		&0			&55 200		&55 200		&0		&0		&23 900	&24 900	\\
			Rest Asia 		&14 600	&326 600		&454 500		&533 500		&29 800	&276 800	&366 000	&474 900	\\
			Oceania 		&0		&1 163 000	&1 176 000	&1 179 000	&0		&449 000	&497 000	&500 000	\\
			Brazil	 	&139 900	&157 700		&157 700		&157 700		&0		&73 600	&121 000	&121 000	\\
			Rest America 	&0		&7 000		&11 700		&13 800		&0		&4 400	&10 100	&11 300	\\
			Africa	 	&93 800	&194 600		&644 100		&663 400		&78 500	&121 800	&260 000	&286 300	\\
			Middle East	&0		&44 000		&44 000		&44 700		&0		&67 800	&67 800	&69 200	\\
			\hline
			Total & 569 900 & 2 516 100 & 3 524 900 & 4 004 500 & 226 600 & 1 245 300 & 1 989 750 & 2 544 800\\
			\hline
			\hline
			Region & \multicolumn{4}{ l|}{Prognosticated}  & \multicolumn{4}{ l }{Speculative}\\
			\hline
			Name &  & $<$80\$/kg  & $<$130\$/kg  & $<$260\$/kg  & & $<$130\$/kg  & $<$260\$/kg & Unassigned \\
			\hline
			\hline
			USA 			&	&819 500	&1 169 300	&1 036 950	&	&858 000		&858 000		&482 000	\\
			Canada  		&	&50 000	&150 000		&150 000		&	&700 000		&700 000		&0		\\
			EU-15 		&	&7 000	&7 600		&7 600		&	&50 100		&50 100		&94 000	\\
			Rest Europe 	&   	&200	&41 650		&66 650		&	&6 650		&126 650		&314 000	\\
			Russia 		&	&0		&182 000		&182 000		&	&0			&0			&633 000	\\
			China 		&	&3600	&3 600		&3 600		&	&4 100		&4 100		&0		\\
			Japan	 	&	&0		&3 300		&3 300		&	&3 300		&3 300		&0		\\
			India 		&	&0		&0			&63 600		&	&0			&0			&17 000	\\
			Rest Asia 		&	&377 900	&591 400		&592 900		&	&1 776 600	&1 806 100	&264 700	\\
			Oceania 		&	&300 000	&300 000		&300 000		&	&0			&0			&500 000	\\
			Brazil	 	&	&73 600	&121 000		&121 000		&	&121 000		&121 000		&0		\\
			Rest America 	&	&6 600	&23 500		&23 500		&	&236 700		&236 700		&176 200	\\
			Africa	 	&	&49 400	&156 900		&156 900		&	&25 000		&25 500		&1 112 900	\\
			Middle East	&	&67 800	&89 000		&89 000		&	&84 800		&98 800		&0		\\
			\hline
			Total &  &  1 755 600 & 2 839 250 & 2 797 000 &  & 3 866 250 & 4 029 750 & 3 593 800\\
			\hline
		\end{tabular*}
	\caption{Uranium resources (in natural concentration) by world region in units of tonnes.}
	\label{tab:U}
\end{center}

\begin{center}
		\begin{tabular*}{1\columnwidth}{@{\extracolsep{\fill}} l|r r r r }
			\hline
			\multicolumn{1}{ l }{Thorium} &\multicolumn{1}{ l }{t}&\multicolumn{1}{ l }{\cite{IAEA2009}}&\multicolumn{2}{ l }{}\\
			\hline
			Region & RAR & Inferred & Indentified & Prognosticated\\
			Name  & $<$ 80 USD/kg & $<$ 80 USD/kg & $<$ 80 USD/kg & N/A \\
			\hline
			\hline
			USA 			&122 000	&278 000	&400 000	&274 000	\\
			Canada  		&0		&44 000	&44 000	&128 000	\\
			EU-15 		&0		&0		&0		&0		\\
			Rest Europe 	&54 000 	&213 000	&186 000	&164 000	\\
			Russia 		&75 000	&112 500	&75 000	&0	\\
			China 		&0		&0		&0		&0	\\
			Japan	 	&0		&0		&0		&0	\\
			India 		&319 000	&478 500	&319 000	&0	\\
			Rest Asia 		&0		&0		&0		&0	\\
			Oceania 		&46 000	&406 000	&452 000	&0	\\
			Brazil	 	&172 000	&130 000	&302 000	&330 000	\\
			Rest America 	&0		&300 000	&300 000	&0		\\
			Africa	 	&18 000	&127 000	&118 000	&410 000	\\
			Middle East	&0		&0		&0		&0	\\
			Unassigned	&23 000	&10 000	&33 000	&81 000	\\
			\hline
			Total &  829 000 & 2 099 000 & 2 229 000 & 1 387 000\\
			\hline
			\multicolumn{5}{ l }{} \\
			\hline
			\multicolumn{1}{ l }{Costs} &\multicolumn{2}{ l }{USD2008 / kg}&\multicolumn{2}{ l }{}\\
			\hline
			 & RAR & Inferred & Indentified & Prognosticated\\
			 \hline
			 Lower & 40 & 40 & 40 & 80\\
			 Upper & 80 & 80 & 80 & 260\\
			\hline
		\end{tabular*}
	\caption{Thorium resources (in natural concentration) by world region in units of tonnes.}
	\label{tab:Th}
\end{center}
\end{table}

\end{document}